\documentclass{amsbook}
\usepackage{amssymb,amscd,amsxtra,verbatim}
\usepackage[vcentermath,enableskew]{youngtab}
\usepackage[all]{xy}
\usepackage{pstricks}
\usepackage{hyperref}
\usepackage{setspace}

\newlength{\cellsz}
\newcounter{cellsize}
\newcommand{\setcellsize}[1]{%
  \setcounter{cellsize}{#1}%
  \setlength{\cellsz}{\value{cellsize}\unitlength}}%
\setcellsize{13}%
\newcommand\cellify[1]{\def\thearg{#1}\def\nothing{}%
\ifx\thearg\nothing \vrule width0pt height\cellsz depth0pt\else
\hbox to 0pt{{\begin{picture}(\value{cellsize},\value{cellsize})
  \put(0,0){\line(1,0){\value{cellsize}}}
  \put(0,0){\line(0,1){\value{cellsize}}}
  \put(\value{cellsize},0){\line(0,1){\value{cellsize}}}
  \put(0,\value{cellsize}){\line(1,0){\value{cellsize}}} \end{picture} \hss}}\fi%
\vbox to \cellsz{ \vss \hbox to \cellsz{\hss$#1$\hss} \vss}}
\newcommand\tableau[1]{\vcenter{\vbox{\let\\\cr
\baselineskip -16000pt \lineskiplimit 16000pt \lineskip 0pt
\ialign{&\cellify{##}\cr#1\crcr}}}}
\newcommand\tabl[1]{\vtop{\let\\\cr
\baselineskip -16000pt \lineskiplimit 16000pt \lineskip 0pt
\ialign{&\cellify{##}\cr#1\crcr}}}

\newtheorem{thm}{Theorem}[chapter]
\newtheorem{prop}[thm]{Proposition}
\newtheorem{lem}[thm]{Lemma}
\newtheorem{sublemma}[thm]{Sublemma}
\newtheorem{cor}[thm]{Corollary}
\newtheorem{conj}[thm]{Conjecture}

\theoremstyle{definition}
\newtheorem{definition}[thm]{Definition}
\newtheorem{ex}[thm]{Example}
\newtheorem{example}[thm]{Example}
\theoremstyle{remark}
\newtheorem{remark}[thm]{Remark}

\newtheorem{property}[thm]{Property}

\numberwithin{section}{chapter} \numberwithin{equation}{chapter}

\makeindex

\begin{document}
\frontmatter
\newcommand{\A}{\mathrm{Case A}}
\newcommand{\At}{\widetilde{A}}
\newcommand{\ba}[1]{\overline{#1}}
\newcommand{\B}{\mathrm{Case B}}
\newcommand{\Bo}{\mathcal{B}}
\newcommand{\bounded}{b}
\newcommand{\Bounded}{\mathcal{B}_n}
\newcommand{\colsums}{\mathrm{colsums}}
\newcommand{\content}{\mathrm{content}}
\newcommand{\CP}{\mathbb{P}}
\newcommand{\CC}{\mathbb{C}}
\newcommand{\ct}{\widetilde{c}}
\newcommand{\code}{\mathrm{code}}
\newcommand{\core}{c}
\newcommand{\Core}{\mathcal{C}_n}
\newcommand{\diag}{\mathrm{diag}}
\newcommand{\C}{\mathrm{Case C}}
\newcommand{\G}{\mathcal{G}}
\newcommand{\Gr}{\mathrm{Gr}}
\newcommand{\id}{\mathrm{id}}
\newcommand{\ipp}[1]{\langle #1 \rangle}
\newcommand{\inner}[2]{\langle #1\,,\,#2\rangle}
\newcommand{\ins}{\mathrm{inside}}
\newcommand{\Inv}{\mathrm{Inv}}
\newcommand{\II}{\mathcal{I}}
\newcommand{\IIL}{\II^\circ}
\newcommand{\La}{\Lambda}
\newcommand{\la}{\lambda}
\newcommand{\M}{\mathcal{M}}
\newcommand{\out}{\mathrm{outside}}
\newcommand{\OO}{\mathcal{O}}
\newcommand{\OOL}{\OO^\circ}
\newcommand{\pt}{\widetilde{p}}
\newcommand{\remind}[1]{\textbf{* #1 *}}
\newcommand{\rowsums}{\mathrm{rowsums}}
\newcommand{\res}{\mathrm{res}}
\newcommand{\Res}{\mathrm{Res}}
\newcommand{\SA}{\scriptsize{0}}
\newcommand{\SB}{\scriptsize{1}}
\newcommand{\say}[1]{{\bf*#1*}}
\newcommand{\scov}{\gtrdot}
\newcommand{\scovby}{\lessdot}
\newcommand{\sge}{\ge}
\newcommand{\sgtr}{>}
\newcommand{\sh}{\mathrm{shape}}
\newcommand{\size}{\mathrm{size}}
\newcommand{\sle}{\le}
\newcommand{\sless}{<}
\newcommand{\slh}{\widehat{\mathfrak{sl}}}
\newcommand{\smc}[1]{\scovby_{#1}}
\newcommand{\strong}{\mathrm{strong}}
\newcommand{\SSS}{\mathrm{Strong}}
\newcommand{\tabpair}[2]{\left(\,#1\,,#2\,\right)}
\newcommand{\tS}{\tilde{S}_n}
\newcommand{\tF}{\tilde{F}}
\newcommand{\vn}{\varnothing}
\newcommand{\wge}{\succeq}
\newcommand{\wgtr}{\succ}
\newcommand{\wle}{\preceq}
\newcommand{\wless}{\prec}
\newcommand{\wt}{\mathrm{wt}}
\newcommand{\weak}{\mathrm{weak}}
\newcommand{\WS}{\mathrm{Weak}}
\newcommand{\xt}{\widetilde{x}}
\newcommand{\yt}{\widetilde{y}}
\newcommand{\zz}{}
\newcommand{\Z}{\mathbb{Z}}

\newcommand{\T}{\rule{0pt}{2.6ex}}
\newcommand{\Bot}{\rule[-1.2ex]{0pt}{0pt}}

\newcommand{\Ah}{{A^\vee}}

\newcommand{\fixit}[1]{\texttt{*** #1 ***}}

\newcommand{\stc}[3]{#1\overset{#2}{\longrightarrow}#3}
\newcommand{\dws}[3]{#1 \overset{#2}{\tilde{\rightsquigarrow}}#3}
\newcommand{\ws}[3]{#1 \overset{#2}{\rightsquigarrow}#3}
\newcommand{\fc}{\mathrm{first}}
\newcommand{\lc}{\mathrm{last}}
\newcommand{\wci}[5]{\ws{#1}{#2}{#3},\,\,\stc{#1}{#4}{#5}}
\newcommand{\wcf}[5]{\ws{#1}{#2}{#3},\,\,\stc{#4}{#5}{#3}}

\newcommand{\um}{u^-}
\newcommand{\am}{a^-}
\newcommand{\bm}{b^-}
\newcommand{\im}{i^-}
\newcommand{\jm}{j^-}
\newcommand{\ip}{i^+}
\newcommand{\jp}{j^+}
\newcommand{\ap}{a^+}
\newcommand{\asp}{a_1}
\newcommand{\aspm}{a_1^-}
\newcommand{\bspm}{b_1^-}
\newcommand{\bp}{b^+}
\newcommand{\bsp}{b_1}
\newcommand{\Cm}{C^-}
\newcommand{\Cp}{C^+}
\newcommand{\Cpm}{{C'}^-}
\newcommand{\Cpp}{{C'}^+}
\newcommand{\is}{i_1}
\newcommand{\js}{j_1}
\newcommand{\xp}{x^+}
\newcommand{\xx}{{\hphantom{x}}}
\newcommand{\iins}[1]{\overset{#1}{\mapsto}}

\newcommand{\pushout}[3]{(#1)\searrow^{\!\!\!\!\mathcal{#2}}(#3)}

\newcommand{\pullback}[3]{(#1)\nwarrow^{\!\!\!\mathcal{#2}}(#3)}

\newcommand{\IA}{\mathcal{A}}
\newcommand{\IB}{\mathcal{B}}
\newcommand{\IS}{\mathcal{S}}
\newcommand{\IX}{\mathcal{X}}

\def\waa{$\tableau{\\{4}\\{3}\\{2}&{4}\\{1}&{3}}$}
\def\wba{$\tableau{\\{4}\\{3}\\{1}&{2}&{3}}$}
\def\wca{$\tableau{\\{3}&{4}\\{1}&{2}&{3}&{4}}$}
\def\wbb{$\tableau{\\{3}\\{2}\\{1}&{3}&{4}}$}

\def\strongaa{$\tableau{\\{4^*}\\{3}\\{3^*}&{4}\\{1^*}&{2^*}}$}
\def\strongba{$\tableau{\\{3}\\{3^*}\\{1^*}&{2^*}&{4^*}}$}
\def\strongca{$\tableau{\\{2^*}&{4}\\{1^*}&{3}&{3^*}&{4^*}}$}

\def\strongab{$\tableau{\\{4}\\{3}\\{3^*}&{4^*}\\{1^*}&{2^*}}$}
\def\strongbb{$\tableau{\\{4^*}\\{2^*}\\{1^*}&{3}&{3^*}}$}
\def\strongcb{$\tableau{\\{2^*}&{4^*}\\{1^*}&{3}&{3^*}&{4}}$}

\def\strongac{$ \tableau{\\{4}\\{3}\\{2^*}&{4^*}\\{1^*}&{3^*}}$}
\def\strongbc{$\tableau{\\{4^*}\\{3}\\{1^*}&{2^*}&{3^*}}$}
\def\strongcc{$\tableau{\\{3^*}&{4^*}\\{1^*}&{2^*}&{3}&{4}}$}

\def\strongad{$\tableau{\\{4^*}\\{3}\\{2^*}&{4}\\{1^*}&{3^*}}$}
\def\strongbd{$\tableau{\\{4^*}\\{3^*}\\{1^*}&{2^*}&{3}}$}
\def\strongcd{$\tableau{\\{3}&{4^*}\\{1^*}&{2^*}&{3^*}&{4}}$}

\def\strongae{$ \tableau{\\{4}\\{3^*}\\{2^*}&{4^*}\\{1^*}&{3}}$}
\def\strongbe{$\tableau{\\{3^*}\\{2^*}\\{1^*}&{3}&{4^*}}$}
\def\strongce{$\tableau{\\{3^*}&{4}\\{1^*}&{2^*}&{3}&{4^*}}$}

\def\strongaf{$ \tableau{\\{4^*}\\{3^*}\\{2^*}&{4}\\{1^*}&{3}}$}
\def\strongbf{$\tableau{\\{3}\\{2^*}\\{1^*}&{3^*}&{4^*}}$}
\def\strongcf{$\tableau{\\{3}&{4}\\{1^*}&{2^*}&{3^*}&{4^*}}$}

\title[Affine insertion and Pieri]{Affine insertion and Pieri rules for the affine Grassmannian}

\thanks{This project was partially supported by NSF grants DMS-0652641, DMS-0652652, DMS-0652668, and
DMS-0652648.}
\author{Thomas Lam}
\address{Department of Mathematics, Harvard University, Cambridge MA
02138 USA} \email{tfylam@math.harvard.edu}
\thanks{T. L. was partially supported by NSF DMS--0600677.}

\author{Luc Lapointe}
\address{Instituto de Matem\'{a}tica Y F\'{i}sica, Universidad de Talca, Casilla 747, Talca, Chile}
\email{lapointe@inst-mat.utalca.cl}
\thanks{L. L. was partially supported by the Anillo Ecuaciones Asociadas a Reticulados financed by the World
Bank through the Programa Bicentenario de Ciencia y Tecnolog\'{\i}a, and by the Programa Reticulados y Ecuaciones
of the Universidad de Talca}
\author{Jennifer Morse}
\address{Department of Mathematics, Drexel University,
Philadelphia, PA 19104 USA} \email{morsej@math.drexel.edu}
\thanks{J. M. was partially supported by NSF DMS--0638625.}
\author{Mark Shimozono}
\address{Department of Mathematics, Virginia Tech,
         Blacksburg, VA 24061 USA}
\email{mshimo@vt.edu} \thanks{M. S. was partially supported by NSF
DMS--0401012.}


\begin{abstract}
We study combinatorial aspects of the Schubert calculus of the
affine Grassmannian ${\rm Gr}$ associated with $SL(n,\mathbb{C})$.
Our main results are:
\begin{itemize}
\item Pieri rules for the Schubert bases of $H^*({\rm Gr})$ and
$H_*({\rm Gr})$, which expresses the product of a special Schubert
class and an arbitrary Schubert class in terms of Schubert classes.
\item A new combinatorial definition for $k$-Schur
functions, which represent the Schubert basis of $H_*({\rm Gr})$.
\item A combinatorial interpretation of the pairing $H^*({\rm Gr})\times
H_*({\rm Gr}) \rightarrow\Z$ induced by the cap product.
\end{itemize}
These results are obtained by interpreting the Schubert bases of
${\rm Gr}$ combinatorially as generating functions of objects we
call strong and weak tableaux, which are respectively defined using
the strong and weak orders on the affine symmetric group. We define
a bijection called affine insertion, generalizing the
Robinson-Schensted Knuth correspondence, which sends certain biwords
to pairs of tableaux of the same shape, one strong and one weak.  
Affine insertion offers a duality between the weak and strong orders
which does not seem to have been noticed previously.

Our cohomology Pieri rule conjecturally extends to the affine
flag manifold, and we give a series of related combinatorial
conjectures.
\end{abstract}

\maketitle \setcounter{page}{4}

\tableofcontents

\chapter*{Introduction}

Let $\Gr=G(\CC((t)))/G(\CC[[t]])$ denote the affine Grassmannian of
$G =SL(n,\CC)$, where $\CC[[t]]$ is the ring of formal power series
and $\CC((t))=\CC[[t]][t^{-1}]$ is the ring of formal Laurent
series. Since $\Gr \cong \mathcal{G}/\mathcal{P}$ for an affine
Kac-Moody group $\mathcal{G}$ and a maximal parabolic subgroup
$\mathcal{P}$, we may talk about the Schubert bases
\begin{align*}
  &\{\xi^w \in H^*(\Gr, \Z) \mid w \in \tS^0\} \\
  &\{\xi_w \in H_*(\Gr, \Z) \mid w \in \tS^0\}
\end{align*}
in the cohomology and homology of $\Gr$, where $\tS^0$ is the subset
of the affine symmetric group $\tS$ consisting of the affine
Grassmannian elements, which by definition are the elements of
minimal length in their cosets in $\tS/S_n$. Quillen (unpublished),
and Garland and Raghunathan~\cite{GR} showed that $\Gr$ is
homotopy-equivalent to the group $\Omega SU(n,\mathbb{C})$ of based
loops into $SU(n,\mathbb{C})$, and thus $H^*(\Gr)$ and $H_*(\Gr)$
acquire structures of dual Hopf-algebras. In~\cite{Bott}, Bott
calculated $H^*(\Gr)$ and $H_*(\Gr)$ explicitly -- they can be
identified with a quotient $\La^{(n)}$ and a subring $\La_{(n)}$ of
the ring $\Lambda$ of symmetric functions.  Using an algebraic
construction known as the {\it nilHecke ring}, Kostant and
Kumar~\cite{KK} studied the Schubert bases of $H^*(\Gr)$ (in fact
for flag varieties of Kac-Moody groups) and Peterson~\cite{Pet}
studied the Schubert bases of $H_*(\Gr)$.  Lam~\cite{Lam},
confirming a conjecture of Shimozono, identified the Schubert
classes $\xi^w$ and $\xi_w$ explicitly as symmetric functions in
$\La^{(n)}$ and $\La_{(n)}$.

In cohomology, the Schubert classes $\xi^w$ are given by the {\it
dual $k$-Schur functions} $$\{\tF_w \mid w \in \tS^0\} \subseteq
\La^{(n)},$$ introduced in~\cite{LMdual} by Lapointe and Morse. The
dual $k$-Schur functions are generating functions of objects called
$k$-tableaux \cite{LMcore}. In \cite{LamAS}, this construction was
generalized to the case of an arbitrary affine permutation
$w\in\tS$; $k$-tableaux are replaced by what we call \textit{weak
tableaux} and the generating function of weak tableaux $U$ is called
the affine Stanley symmetric function or \textit{weak Schur
function}
\begin{align*}
  \tF_w(x) = \WS_w(x) = \sum_U x^{\wt(U)}.
\end{align*}
When $w\in S_n\subset\tS$ is a usual Grassmannian permutation
(minimal length coset representative in $S_n/(S_r \times S_{n-r})$),
the affine Stanley symmetric function reduces to a usual Schur
function.

In homology, the Schubert classes $\xi_w$ are given by the {\it
$k$-Schur functions} $$\{s_\lambda^{(k)}(x)\mid \lambda_1 < n \}
\subseteq \Lambda_{(n)}$$ where the $k$ in $k$-Schur or dual
$k$-Schur function always means
\begin{align*}
  k = n-1.
\end{align*}
The $k$-Schur functions were first introduced by Lapointe, Lascoux,
and Morse~\cite{LLM} for the study of Macdonald
polynomials~\cite{Mac}, though so far a direct connection between
Macdonald polynomials and the affine Grassmannian has yet to be
established. A number of conjecturally equivalent definitions of
$k$-Schur functions have been presented
(see~\cite{LLM,LMfil,LMproofs,LMdual}).  In this article, a
$k$-Schur function will always refer to the definition of
\cite{LMproofs,LMdual} and we can thus view
$\{s_\lambda^{(k)}(x)\mid \lambda_1 < n = k+1\}$ as the basis of
$\Lambda_{(n)}$ dual to the basis $\{\tF_w \mid w \in \tS^0\}$ of
$\Lambda^{(n)}$ with respect to a bijection
\begin{equation} \label{EE:boundedtograss}
\begin{split}
  \{\la\mid \la_1<n\} &\to \tS^0 \\
  \la &\mapsto w
\end{split}
\end{equation}
(see Proposition \ref{P:shapetrans}). Given an interval $[v,w]$ in
the strong order (Bruhat order) on $\tS$, we introduce the notion of
a \textit{strong tableau} $T$ of shape $w/v$; it is a certain kind
of labeled chain from $v$ to $w$ in the strong order. We define a
\textit{strong Schur function} to be the generating function of
strong tableaux $T$ of shape $w/v$:
\begin{align*}
  \SSS_{w/v}(x) = \sum_T x^{\wt(T)}.
\end{align*}
One of our main results (Theorem~\ref{T:kSchurStrong}) is that
$k$-Schur functions are special cases of strong Schur functions:
\begin{align*}
  s_\la^{(k)}(x) = \SSS_{w/\id}(x)
\end{align*}
where $\id\in\tS$ is the identity and $\la\mapsto w\in\tS^0$ under
the bijection \eqref{EE:boundedtograss}. When $v,w \in S_n \subset
\tS$ are usual Grassmannian permutations, strong Schur functions
reduce to usual skew Schur functions. Strong tableaux, in the case
of $S_n \subset \tS$, are closely related to chains in the
$k$-Bruhat order (where here $k$ is unrelated to $n$) of Bergeron
and Sottile~\cite{BS}. An important difference is that our chains
are marked, reflecting the fact that affine Chevalley coefficients
are not multiplicity-free (see Remark \ref{R:mark}).


Our main result (Theorem~\ref{T:main}) is the construction of an
algorithmically defined bijection called {\it affine insertion}.  In
its simplest case, affine insertion establishes a bijection between
nonnegative integer matrices with row sums less than $n$, and pairs
$(P,Q)$ where $P$ is a strong tableau, $Q$ is a weak tableau, and
both tableaux start at $\id \in \tS$ and end at the same $v \in
\tS^0$.  This bijection reduces to the usual row-insertion
Robinson-Schensted-Knuth (RSK) algorithm (see~\cite{Ful}) as $n \to
\infty$. Affine insertion yields a combinatorial proof of the
following affine Cauchy identity, which is obtained from
Theorem~\ref{T:genCauchy} by taking $u=v=\id$:
\begin{align*}
\Omega_n(x,y) &= \sum_{w \in \tS^0} \SSS_w(x)\WS_w(y)
\end{align*}
and the affine Cauchy kernel is given by
\begin{align*}
\Omega_n(x,y) &= \prod_{i} \left(1 + y_ih_1(x) + y_i^2 h_2(x) +
\cdots + y_i^{n-1}h_{n-1}(x)\right) \\&= \sum_{\lambda \,:\,
\lambda_1 < n} h_\lambda(x) m_\lambda(y).
\end{align*}

This provides a combinatorial description of the reproducing kernel
of the perfect pairing
\begin{align*}
  H^*(\Gr) \times H_*(\Gr) \to \Z
\end{align*}
induced by the cap product.

The bijection also yields geometric information in the form of Pieri
rules, which are explicit formulae for certain structure constants
with respect to the Schubert bases of $H^*(\Gr)$ and $H_*(\Gr)$. We
use affine insertion to derive the Pieri rules
(Theorem~\ref{T:strongPieri} and~\ref{T:weakPieri}) for strong and
weak Schur functions.

Affine insertion exhibits a duality between the weak and strong
orders which does not seem to have been studied before, even in the
case of the finite symmetric group $S_n$. In particular one may show
that the number of pairs $(P,Q)$ of a standard strong tableau and a
standard weak tableau of the same ``shape'' (starting at the
identity $\id$ and ending at some permutation $w \in \tS^0$) and
size $m$ is given by $m!$. These issues are pursued in \cite{LS}
where such identities are established in Kac-Moody generality. One
may also generalize this duality in a different direction, to obtain
similar identities for combinatorial Hopf algebras.

The construction and proof of the affine insertion algorithm is
reduced to a ``local rule'' using the technology of Fomin's {\it
growth diagrams}~\cite{Fom}.  The local rule, which is constructed
directly on the level of affine permutations, represents the most
involved part of this paper.  Our local rule has many elements which
will be familiar to experts of Schensted insertion, including local
bijections analogous to boxes being bumped to the next row, or boxes
not interfering with each other. The strong covers $x \scovby y$ of
$\tS$ in this article roughly correspond to boxes in the traditional
language of Young tableaux.

As corollaries of the affine insertion theorem, one deduces Pieri
rules for the $H^*(\Gr)$ and $H_*(\Gr)$.  Let $c_{0,m} =
s_{m-1}\cdots s_1s_0 \in \tS$ where $\{s_i \mid i \in \Z/n\Z\}$
denote the simple generators of $\tS$.  We obtain an affine
homology Pieri rule in $H_*(\Gr)$ (Theorem~\ref{T:homPieri}):
$$
\xi_{c_{0,m}}\, \xi_w = \sum_{\ws{w}{}{z}} \xi_z,
$$
where the sum runs over weak strips $\ws{w}{}{z}$ of size $m$ and
an affine cohomology Pieri rule for $H^*(\Gr)$
(Theorem~\ref{T:cohomPieri}):
$$
\xi^{c_{0,m}}\, \xi^w = \sum_{S} \xi^{\out(S)},
$$
where the sum runs over strong strips $S$ of size $m$ with $\ins(S)
= w$.  The affine cohomology Pieri rule is an affine analogue of the
usual Pieri rule for the Schubert calculus of the flag
manifold~\cite{Sot}.  The affine homology Pieri rule does not appear
to have a classical geometric counterpart.  It can be deduced
directly from the identification of the affine homology Schubert
basis $\{\xi_w \mid w \in \tS^0\}$ as $k$-Schur functions
in~\cite{Lam}, together with the Pieri rule for $k$-Schur functions
first stated in~\cite{LMproofs} and described in the notation of
this article in \cite{LMW}.

Conjecture~\ref{conj:flags} asserts that the natural analogue of the
affine cohomology Pieri rule holds for the cohomology $H^*(\G/\Bo)$
of the affine flag variety.  This conjecture is related to a series
of combinatorial conjectures (Conjecture~\ref{conj:strongSym})
concerning the strong Schur functions $\SSS_{v/w}(x)$, including
symmetry and positivity when expressed in terms of $k$-Schur
functions.  The analogous properties for weak Schur functions were
established in~\cite{LamAS,Lam}.

In the last part of our paper we translate weak and strong tableaux,
together with the affine insertion bijection, from permutations into
the more traditional language of partitions. This is performed using
a classical bijection \cite{L,MM} between $\tS^0$ and the set of
partitions which are {\it $n$-cores}. For weak strips and weak
tableaux, the corresponding combinatorics involving cores was worked
out in~\cite{LMcore}.  Our main result here
(Proposition~\ref{P:scover}) gives a purely partition-theoretic
description of marked strong covers, and hence strong strips and
strong tableaux.  As a consequence the affine cohomology Pieri rule
acquires a form similar to that of the Pieri rule for Schur
functions, with horizontal strips replaced by ``strong strips built
on cores''.  We also use the combinatorial description of strong
covers to define a ``spin''-statistic on strong tableaux and
conjecture (Conjecture~\ref{conj:kSchurt}) that the original
$k$-Schur functions (depending on a parameter $t$)
of~\cite{LLM,LMfil} are spin-weight generating functions of strong
tableaux of fixed shape.

Our work poses further challenges for both geometers and
combinatorialists. The two Pieri rules beg for a more geometric
proof; in the cohomology case there should be a geometric proof
similar to that of Sottile~\cite{Sot}, and alternatively a more
algebraic derivation might be possible using the recursive machinery
of Kostant and Kumar's nilHecke ring~\cite{KK}.  The ``monomial
positivity'' of the cohomology classes $\xi^w = \WS_w(x)$ can be
interpreted geometrically as arising from Bott's map
$(\mathbb{CP}^n)^\infty \to \Gr$.  It would be interesting to obtain
a geometric explanation of the monomial positivity of $\xi_w =
\SSS_w(x)$.

The positivity of structure coefficients for both weak and strong
Schur functions have yet to be given a combinatorial interpretation.
The structure constants for the strong Schur functions yields as a
special case the WZW fusion coefficients (or equivalently the
structure constants of the quantum cohomology $QH^*(G/P)$ of the
Grassmannian) as proved by Lapointe and Morse~\cite{LMdual}. More
generally, Peterson has shown that the structure constants of the
quantum cohomology of any (partial) flag manifold can be obtained
from the structure constants of the homology $H_*(\Gr)$ of the
affine Grassmannian. Obtaining a combinatorial interpretation for
these structure constants is likely to be a challenging problem.
Another interesting problem is to give a direct combinatorial proof
of the symmetry of strong Schur functions.

The literature contains other ``affine'' or ``infinite''
generalizations of the Schensted algorithm, see for example
\cite{Pak, Shi}. However, as far as we are aware, these algorithms
biject affine or infinite permutations with certain tableaux, while
our insertion pairs {\it usual} permutations with affine tableaux.

\medskip

{\bf Organization.} This paper is roughly divided into three parts.
In the first part (Chapters~\ref{ch:notation}--\ref{ch:main}) we
give the necessary definitions and present our main theorems.
Chapter~\ref{ch:notation} contains notation for symmetric functions
and Schubert bases of the affine Grassmannian.  Our two main
objects, strong and weak tableaux, are introduced in
Chapters~\ref{ch:strongtableaux} and~\ref{ch:weaktableaux}
respectively. In Chapter~\ref{ch:main}, we present and prove our
main results modulo the proof of affine insertion.  The proof of the
affine insertion algorithm is reduced to properties of local rules
using Fomin's growth diagram machinery.

The local rules for the affine insertion algorithm are defined and
studied in the second, and most technical, part of our paper
(Chapters~\ref{ch:forward}--\ref{ch:bijectivity}). In
Chapters~\ref{ch:forward} and~\ref{ch:inverse} we define the forward
and reverse local rules, and show that they are well-defined. In
Chapter~\ref{ch:bijectivity} we prove that affine insertion is
bijective.

The last part of our paper
(Chapters~\ref{ch:coresbijection}--\ref{ch:coresinsertion}) contains
translations of our combinatorial constructions into the language of
partitions and cores.  In Chapter~\ref{ch:coresbijection}, we
explain a number of bijections between the coroot lattice, affine
Grassmannian permutations, cores, offset sequences and $k$-bounded
partitions. In Chapters~\ref{ch:corestableaux}
and~\ref{ch:coresinsertion}, we explain weak and strong tableaux and
affine insertion using the combinatorial language of cores.

\mainmatter

\chapter{Schubert Bases of $\Gr$ and Symmetric Functions}
\label{ch:notation} \section{Symmetric functions} Here we introduce
notation for symmetric functions, which can be found in greater
detail in~\cite{Mac}.  Let $\La = \La(x)$ denote the ring of
symmetric functions in infinitely many variables $x_1,x_2,\ldots$
over $\Z$. It is generated over $\Z$ by the algebraically
independent {\it homogeneous symmetric functions} $h_1,h_2,\ldots$,
where ${\rm deg}\,h_i = i$.  The ring $\La$ is equipped with an
algebra involution $\omega: \La \rightarrow \La$ given by
$\omega(h_i) = e_i$ where $e_i$ denotes the {\it elementary
symmetric functions}.  For a partition $\lambda = (\lambda_1 \ge
\lambda_2 \ge \dotsm)$ let $h_\lambda: =
h_{\lambda_1}h_{\lambda_2}\cdots$. The {\it Hall inner product}
$\inner{}{}_\La: \La \times \La \to \Z$ is defined by
$\inner{h_\lambda}{m_\mu}_\La = \delta_{\lambda\mu}$ where $m_\mu$
denotes the {\it monomial symmetric functions}.

The ring $\La$ has a coproduct $\Delta:\La \rightarrow \La \otimes_{\Z}
\La$ given by $\Delta(h_i) = \sum_{0 \leq j \leq i} h_{j} \otimes
h_{i-j}$ where $h_0:=1$.  Together with the Hall inner product, this
gives $\La$ the structure of a self-dual commutative and
cocommutative Hopf algebra.  The antipode $c$ is given by $c(h_i) =
(-1)^ie_i$; thus if $f \in \Lambda$ is homogeneous of degree $d$
then $c(f) = (-1)^d \omega(f)$.

Now let $\La^{(n)}:= \La/(m_\la\mid \la_1 \ge n)$.  This is a
quotient Hopf algebra of $\La$.  Let $\La_{(n)} :=
\Z[h_1,h_2,\dotsc,h_{n-1}]$.  This is a sub-Hopf algebra of $\La$.
The Hall inner product gives $\La^{(n)}$ and $\La_{(n)}$ the
structures of dual Hopf algebras.  One possible choice of dual bases
is $\{m_\la \mid \la_1 \leq n-1\}$ for $\La^{(n)}$ and $\{h_\lambda
\mid \la_1 \leq n-1 \}$ for $\La_{(n)}$.  The algebra involution
$\omega$ of $\La$ restricts to an involution of $\La_{(n)}$.  By
duality we also obtain an involution $\omega^+$ of $\La^{(n)}$
characterized by the property $\inner{f}{g}_\La =
\inner{\omega(f)}{\omega^+(g)}_\La$ for $f \in \La_{(n)}$ and $g \in
\La^{(n)}$.  For $f \in \La$ let $\ba{f} \in \Lambda^{(n)}$ denote
its image in the quotient.  Then $\omega^+(\ba{f}) =
\ba{\omega(f)}$. If $f \in \La$, when the context makes it clear we
will just write $f$ to denote its image in $\Lambda^{(n)}$.

\section{Schubert bases of $\Gr$} \label{sec:SchubertGr} Let $\Gr$
denote the affine Grassmannian of $SL(n,{\mathbb C})$. It is an
ind-scheme equipped with a Schubert-decomposition
$$
\Gr = \bigsqcup_{w \in \tS^0} \Omega_w = \bigcup_{w \in \tS^0} X_w
$$
where the unions are taken over the set $\tS^0$ of all
$0$-Grassmannian permutations in the affine symmetric group $\tS$
(see Chapter~\ref{ch:strongtableaux}); $\Omega_w$ denotes the
Schubert cell indexed by $w$ and $X_w$ denotes the Schubert variety.
Let $\xi^w \in H^*(\Gr)$ and $\xi_w \in H_*(\Gr)$ denote the
corresponding Schubert classes in cohomology and homology;
see~\cite{Gra, Lam, Kum}.  The cap product yields a pairing
$$\inner{\cdot}{\cdot}_{\Gr}: H^*(\Gr) \times H_*(\Gr) \rightarrow
\Z$$ under which the Schubert bases $\{\xi^w \mid w \in \tS^0\}$ and
$\{\xi_w \mid w \in \tS^0 \}$ are dual. Throughout this paper, all
(co)homology rings have coefficients in $\Z$. For $u,v,w\in\tS^0$
define $c^w_{uv} \in \Z$ and $d^w_{uv} \in \Z$ by
\begin{eqnarray}
\xi_u \xi_v &=& \sum_w d^w_{uv} \xi_w, \label{E:struc1} \\
\xi^u \xi^v &=& \sum_w c^w_{uv} \xi^w. \label{E:struc2}
\end{eqnarray}  The
structure constants $c^w_{uv}$ were studied in~\cite{KK} using the
{\it nilHecke ring}.  It follows from work of Graham~\cite{Gra}
and Kumar~\cite{Kum} that $c^w_{uv}\in\Z_{\ge0}$ and from work of
Peterson~\cite{Pet} that $d^w_{uv}\in\Z_{\ge0}$. Our work yields
combinatorial interpretations for some of these numbers.

The space $\Gr$ is homotopy-equivalent to the based loop space
$\Omega SU(n)$; see~\cite{GR,PS}. Thus $H_*(\Gr)$ and $H^*(\Gr)$ are
endowed with the structures of dual commutative and co-commutative
Hopf algebras. In~\cite{Bott}, Bott calculated these Hopf algebras
explicitly. By identifying the generators explicitly one obtains
isomorphisms $H^*(\Gr) \cong \La^{(n)}$ and $H_*(\Gr)\cong
\La_{(n)}$ such that the diagram
\begin{align*}
\xymatrix{%
{H^*(\Gr) \times H_*(\Gr)} \ar[dd] \ar[dr]^{\inner{\cdot}{\cdot}_{\Gr}} & \\
& {\Z} \\
{\La^{(n)} \times \La_{(n)}} \ar[ur]_{\inner{\cdot}{\cdot}_{\La}} &
}%
\end{align*}
commutes.

A natural problem is the identification of the Schubert classes
$\xi_w$ and $\xi^w$ as symmetric functions.  Confirming a conjecture
of Shimozono, in~\cite{Lam} Lam showed that the Schubert classes
$\xi_w$ and $\xi^w$ are represented respectively by the {\it
$k$-Schur functions} \cite{LLM,LMfil,LMproofs} and {\it affine Schur
functions} (also called {\it dual $k$-Schur functions}) \cite{LamAS}
\cite{LMdual}.

\begin{thm}
\label{T:Lam} \cite{Lam} Under the isomorphism $H^*(\Gr) \cong
\La^{(n)}$, the Schubert class $\xi^w$ is sent to the affine Schur
function $\tF_w \in \La^{(n)}$. Under the isomorphism $H_*(\Gr)
\cong \La_{(n)}$, the Schubert class $\xi_w$ is sent to the
$k$-Schur function $s^{(k)}_{w} \in \La_{(n)}$, where $k = n-1$.
\end{thm}


The affine Schur functions $\tF_w$ are generating functions of
combinatorial objects known as {\it $k$-tableaux} and were first
introduced by Lapointe and Morse in~\cite{LMcore}.  We shall define
these objects in Chapter~\ref{ch:weaktableaux}, following the
approach of~\cite{LamAS}. It is shown in~\cite{LamAS,LMdual} that
$\{\tF_w \mid w \in \tS^0\}$ forms a basis of $\La^{(n)}$.  The
$k$-Schur functions $\{s^{(k)}_{w} \mid w \in \tS^0\}$ form the dual
basis of $\La_{(n)}$ to the affine Schur functions.  The $k$-Schur
functions $s^{(k)}_\la(x)$ or $s_w^{(k)}(x)$ used here are
conjecturally (see~\cite{LMproofs}) the $t=1$ specializations of the
$k$-Schur functions $s^{(k)}_\la(x;t)$ first introduced by Lascoux,
Lapointe, and Morse in~\cite{LLM} to study Macdonald polynomials.
The $k$-Schur functions (and affine Schur functions) are usually
indexed by partitions $\la$ such that $\la_1 \leq k$. We will
describe the bijection between the sets $\{\la \mid \la_1 \leq k\}$
and $\{ w \in \tS^0\}$ in Chapter~\ref{ch:coresbijection}. For the
first portion of this paper we will use affine permutations as
indices for affine Schur and $k$-Schur functions.

\section{Schubert basis of the affine flag variety} Let $\G/\Bo$
denote the flag variety for the affine Kac-Moody group
$\widehat{SL}(n,{\mathbb C})$. Again we omit the construction of
this ind-scheme and refer the reader to~\cite{Kum}. The space
$\G/\Bo$ has a decomposition into Schubert varieties indexed by
affine permutations $w \in \tS$. We let $\xi_B^w \in H^*(\G/\Bo)$
for $w \in \tS$ denote the cohomology Schubert basis of the affine
flag variety.  There is a (surjective) morphism $\G/\Bo \to \Gr$
which induces an algebra inclusion $\iota: H^*(\Gr) \hookrightarrow
H^*(\G/\Bo)$.  The Schubert classes are sent to Schubert classes
under $\iota$, so that $\iota(\xi^w) = \xi^w_B$ for $w \in \tS^0$.
Thus we may define integers $c^w_{uv} \in \Z$ as the structure
constants of $H^*(\G/\Bo)$:
\begin{align} \label{E:affineflagconstants}
\xi^u_B \xi^v_B = \sum_w c^w_{uv} \xi^w_B
\end{align}
and when $w,u,v \in \tS^0$ this agrees with the definition in
Section~\ref{sec:SchubertGr}.  Again, by general results
of~\cite{Gra,Kum}, the integers $c^w_{uv}$ are nonnegative.

\chapter{Strong Tableaux} \label{ch:strongtableaux} We introduce
some combinatorial constructions involving the affine symmetric
group $\tS$. For $a\in\Z$ let $\overline{a}$ denote the coset
$a+n\Z\in \Z/n\Z$.

\section{$\tS$ as a Coxeter group} \label{SS:coxeter}

The affine symmetric group $\tS$ is an infinite Coxeter group, with
generators $\{s_0,s_1,\dotsc,s_{n-1}\}$ of \textit{simple
reflections} and relations $s_i^2=\id$, $s_is_j=s_js_i$ if $\ba{i}$
and $\ba{j}$ are not adjacent in $\Z/n\Z$, and $s_is_js_i=s_js_is_j$
if $\ba{i}$ and $\ba{j}$ are adjacent in $\Z/n\Z$. Here pairs of
elements of the form $\{\ba{i},\ba{i+1}\}$ are adjacent in $\Z/n\Z$
and other pairs of elements are not. The length $\ell(w)$ of
$w\in\tS$ is the number $\ell$ of simple reflections in a reduced
decomposition of $w$, which by definition is a factorization
$w=s_{i_1}s_{i_2}\dotsm s_{i_\ell}$ of $w$ into a minimum number of
simple reflections.

For $a,b\in\Z$ we write $s_a=s_b$ if $\ba{a}=\ba{b}$. A
\textit{reflection} is an element that is conjugate to a simple
reflection.

\subsection{$\tS$ as periodic permutations of $\Z$}
The affine symmetric group $\tS$ can be realized as the set of
permutations $w$ of $\Z$ such that $w(i+n)=w(i)+n$ for all
$i\in\Z$ and $\sum_{i=1}^n (w(i)-i)=0$ (see \cite{Lus}). We
sometimes specify an element $w\in\tS$ by ``window" notation
\[
w = [w(1),w(2),\dotsc,w(n)]
\]
as this uniquely determines $w$. Multiplication of elements
$v,w\in\tS$ is given by function composition: $(vw)(i)=v(w(i))$
for all $i\in \Z$. We recall Shi's length formula~\cite{Shi}
\begin{equation}\label{E:Shi}
\ell(w) = \sum_{1\le i<j\le n} \left| \left\lfloor
\dfrac{w(j)-w(i)}{n} \right\rfloor \right|.
\end{equation}
An \textit{inversion} of $w\in\tS$ is a pair $(i,j)\in\Z^2$ such
that $\ba{i}\ne\ba{j}$, $i<j$, $1\le i\le n$, and $w(i)>w(j)$. Let
$\Inv(w)$ denote the set of inversions of $w$. Then we have (see
\cite{H})
\begin{equation} \label{E:leninv}
  \ell(w) = |\Inv(w)|.
\end{equation}

For $r,s\in\Z$ with $\ba{r}\not=\ba{s}$ let $t_{r,s}$ be the unique
element of $\tS$ defined by $t_{r,s}(r)=s$, $t_{r,s}(s)=r$, and
$t_{r,s}(i)=i$ for $\ba{i}\not\in\{\ba{r},\ba{s}\}$. We have
\begin{align}
\label{E:tconj} w\,t_{ij}w^{-1} =t_{w(i),w(j)} \qquad\text{for all
$w\in\tS$.}
\end{align}
The simple reflections are given by $s_i = t_{i,i+1}$ for $i\in\Z$.
By \eqref{E:tconj} the reflections in $\tS$ are precisely the
elements of the form $t_{r,s}$.

\begin{ex} For $n=3$, $t_{0,4}=[-3,2,7]$ since
$t_{0,4}(1)=t_{0,4}(4)-3=0-3=-3$ and
$t_{0,4}(3)=t_{0,4}(0)+3=4+3=7$.
\end{ex}

One may also specify $w\in\tS$ in ``two-line notation" by expanding
the window to include all values of $w:\Z\rightarrow\Z$:
\begin{equation*}
\begin{matrix}
\dotsm & -2 & -1 & 0 & 1 & 2 & \dotsm \\ %
 \dotsm & w(-2)&w(-1)&w(0)&w(1)&w(2)&\dotsm
\end{matrix}
\end{equation*}
In this notation $t_{ij} \, w $ is obtained from $w$ by the
exchanging the values $i+kn$ and $j+kn$ in the lower row for all
$k\in\Z$. By \eqref{E:tconj} the permutation $w \,t_{ij}$ is
obtained from $w$ by exchanging the elements in the
\textit{positions} $i+kn$ and $j+kn$ for all $k\in\Z$.

\section{Fixing a maximal parabolic} \label{SS:maxparabolic} For
the duration of the paper, fix $l \in \Z$. Let $S^{\ba{l}}_n \subset
\tS$ be the maximal parabolic subgroup generated by $s_i$ for
$\ba{i}\not=\ba{l}$. It is isomorphic to $S_n$. The definitions of
strong strip and strong tableau in Section \ref{SS:strong} depend on
the choice of $l$. We say that an element $w\in\tS$ is
\textit{$l$-Grassmannian} if it is of minimum length in its coset in
$\tS/S^{\ba{l}}_n$, that is, $w(l+1) < w(l+2) < \dotsm < w(l+n)$. We
say that $w$ is \textit{Grassmannian} to mean that it is
$0$-Grassmannian.

\section{Strong order and strong tableaux} \label{SS:strong}

The strong order on $\tS$ is by definition the partial order with
covering relation $w\lessdot u$, which holds exactly when
$\ell(u)=\ell(w)+1$ and $wt_{ij}=u$ for some reflection $t_{ij}$.

For $a,b\in\Z$ we sometimes use the notation $[a,b]=\{c\in\Z \mid
a\le c\le b\}$ and $(a,b)=\{c\in\Z\mid a < c < b\}$ for intervals of
integers.

\begin{lem}
\label{L:strongcover} Let $w \in \tS$ and let $i < j$ be integers
such that $\ba{i}\not=\ba{j}$. Then
\begin{enumerate}
\item
$w \lessdot w\, t_{ij}$ if and only if $w(i) < w(j)$ and for each $k
\in (i,j)$, $w(k)\not\in[w(i),w(j)]$.
\item
$wt_{ij}\lessdot w$ only if $w(j)<w(i)$ and for each $k \in (i,j)$,
$w(k)\not\in[w(j),w(i)]$.
\end{enumerate}
Moreover if the strong cover holds then either $j-i<n$ or
$|w(i)-w(j)|<n$.
\end{lem}
\begin{proof} Suppose $w(i)<w(j)$. We have
$$t_{ij} = (s_is_{i+1} \dotsm s_{i+a-1})s_{i+a}(s_{i+a-1}
\dotsm s_{i+1}s_i)$$ for $i+a=j-1$. Let
\begin{align*}
w^{(r)} &= w \, s_i s_{i+1} \dotsm s_{i+r-1} &\qquad&\text{for $r \in [0,a+1]$} \\
w_{(r)} &= w^{(a+1)} \, s_{i+a-1} \dotsm s_{i+a-r} && \text{for $r
\in [0,a]$.}
\end{align*}
We have  $\ell(w^{(r)}) = \ell(w^{(r-1)}) \pm 1$ for $r \in [1,a+1]$ and
$\ell(w_{(r)}) = \ell(w_{(r-1)}) \pm 1$ for $r \in [1,a]$.
Note that in
passing from $w = w^{(0)}$ to $w^{(a)}$ (resp. $w^{(a+1)}=w_{(0)}$
to $w_{(a)}=wt_{ij}$), given $A,B$ and $C$ as defined below,
we compare $w(i)$ (resp. $w(j)$) with all
$w(k)$ where $k \in A\cup B$  (resp. $k\in A \cup C$) and count the number of inversions that have been gained or lost.
Note that
$w(j)>w(i)>w(k)$ for all $k\in C$. Including the ``middle" step
$w^{(a)}$ to $w^{(a+1)}$ that exchanges the elements $w(i)$ and
$w(j)$ in adjacent positions $j-1$ and $j$ respectively, and writing
\begin{align*}%
A &= \{k \in (i,j) \mid \ba{k} \notin \{\ba{i},\ba{j}\}\} \\
A_- &= \{k \in A \mid w(k) < w(i) \}\\
A_+ &= \{k \in A \mid w(k) > w(j)\} \\
A_= &= \{k \in A \mid w(i) < w(k) < w(j)\} \\
B &= \{k \in (i,j) \mid \ba{k} =
\ba{j} \} \\
B_+ &= \{k \in B \mid w(k) > w(i) \}
\\
B_- &= \{k \in B \mid w(k) < w(i)\} \\
C& =\{i-n,i-2n,\dotsc,i-\lfloor(w(j)-w(i))/n\rfloor n\},
\end{align*}
we have
\begin{align*}
\ell(w^{(a+1)} ) &= \ell(w) + |A_=| + |A_+| - |A_-| + |B_+|-|B_-|+1; \\
\ell(w \, t_{ij})  & =\ell(w_{(a)})  = \ell(w^{(a+1)}) +|A_=| -|A_+|+|A_-| +|C|.
\end{align*}
But $|B|=|B_+|+|B_-|=|C|$. So $\ell(w \, t_{ij})-\ell(w) = 2|A_=| +
2|B_+|+1$. Therefore $w \lessdot w \, t_{ij}$ if and only if $|A_=|
= 0$ and $|B_+| = 0$. But the following are equivalent: $B_+$ is
empty; $j-i<n$ or $w(j)-w(i)<n$; for every $k\in(i,j)$ with
$\ba{k}\in\{\ba{i},\ba{j}\}$, $w(k)$ is not in the interval
$(w(i),w(j))$. Part (1) follows.

The statement in Part (2) follows immediately from Part (1).
\end{proof}

\begin{example}
Even if $w(j)<w(i)$ and for each $k \in (i,j)$ we have
$w(k)\not\in[w(j),w(i)]$, it is not always true that
$wt_{ij}\lessdot w$.  For example, let $n = 3$, $w = [10,2,-6]$, $i
= 1$, $j = 5$.  Then $w$ is not a strong cover of $wt_{ij} =
[5,7,-6]$.
\end{example}


A {\it marked strong cover} $C=(\stc{w}{i,j}{u})$ consists of
$w,u\in \tS$ and an ordered pair $(i,j)\in \Z^2$ such that
\begin{enumerate}
\item
$w \lessdot u$ is a strong cover with $w\, t_{ij} = u$.
\item
The reflection $t_{ij}$ \textit{straddles} $l$, that is, $i \le l <
j$, where $l$ is defined in Section \ref{SS:maxparabolic}.
\end{enumerate}
We use the notation $\ins(C)=w$ and $\out(C)=u$. We say that $C$ is
{\it marked} at the integer
\begin{align} \label{E:mark}
  m(C) = w(j) = u(i).
\end{align}

\begin{remark} \label{R:mark} Let $w \lessdot u$ be a strong cover.
The number of pairs $(i,j)$ such that $\stc{w}{i,j}{u}$ is a marked
strong cover, is equal to the affine Chevalley multiplicity, which
by definition is the structure constant $c_{w,s_l}^u$ for the
Schubert basis of the cohomology $H^*(\G/\Bo)$ of the affine flag
variety. This is merely a translation of the Chevalley rule in
\cite{KK} for a Kac-Moody flag manifold, in the special case of the
affine flag variety.
\end{remark}

\begin{remark}
Our notion of strong marked cover when restricted to the symmetric
group $S_n \subset \tS$ essentially produces the $k$-Bruhat order
studied by Sottile~\cite{Sot} and Bergeron and Sottile~\cite{BS}.
\end{remark}

\begin{prop} \label{P:strongGrass} If
$\stc{w}{i,j}{u}$ is a marked strong cover and $w$ is
$l$-Grass-mannian then $u$ is $l$-Grassmannian as well.
\end{prop}
\begin{proof}
Suppose $w$ is $l$-Grassmannian so that $w(l+1) < w(l+2) < \cdots <
w(l+n)$.  We may pick $i',j'$ so that $t_{i,j} = t_{i',j'}$ and $i'
\leq l <j' \leq l+n$.  By Lemma~\ref{L:strongcover}, we must have
$w(i') < w(j')$ and for each $k$ satisfying $l+1 \leq k < j'$ we
must have $w(k) < w(i')$ (since $w(k) < w(j')$).  In particular,
since $k > i'$, we must have $\ba{k} \neq \ba{i'}$.  Thus $u(l+1) <
u(l+2) < \cdots < u(j')$.

Now suppose $\ba{i'} = \ba{a}$ where $j' < a \leq l+n$.  Then $u(a)
> w(a)$ so we must have $u(a) > u(a-1) > \dotsm > u(j')$.  Now for
each $k$ satisfying $a < k \leq l+n$, we have $i' < k-bn < j'$ where
$i' = a - bn$ for a positive integer $b$.  Thus $u(a) = w(j') + bn$.
Again by Lemma~\ref{L:strongcover}, we have either $w(k-bn) < w(i')
= w(a) - bn$ or $w(k-bn) > w(j')$.  Since $w(k) > w(a)$ the first
situation does not occur.  Thus $u(k) = w(k)
> w(j') + bn = u(a)$.  Combining the inequalities, we conclude that
$u$ is $l$-Grassmannian.
\end{proof}

A \textit{strong tuple} $S=[w;(C_1,C_2,\dotsc,C_r);u]$ consists of a
``sequence of marked strong covers from $w$ to $u$", that is,
elements $w,u\in\tS$ and a sequence of marked strong covers $C_k$
such that $\out(C_k)=\ins(C_{k+1})$ for each $0\le k \le r$, where
by convention $\out(C_0)=w$ and $\ins(C_{r+1})=u$. This is a certain
kind of chain in the strong order from $w$ to $u$ with data attached
to the covers. Sometimes $w$ and $u$ are suppressed in the notation.
We write $\ins(S)=w$ and $\out(S)=u$. The size of $S$ is the number
$r$ of covers in $S$. By definition $\size(S)=\ell(u)-\ell(w)$. If
$\size(S)>0$, we refer to the first and last covers in $S$ by
$\fc(S)=C_1$ and $\lc(S)=C_r$. If $\size(S)=0$ then the empty strong
tuple $S$ is determined only by the element $\ins(S)=\out(S)$. We
also use the notation $C^-$ and $C^+$ to refer respectively to the cover
before and after the cover $C$ within a given strong tuple
$S=(\dotsc,C^-,C,C^+,\dotsc)$. If $C=(\stc{w}{i,j}{u})$ then we use
the notation $C^-=(\stc{w^-}{i^-,j^-}{u^-})$ and
$C^+=(\stc{w^+}{i^+,j^+}{u^+})$.

A \textit{strong strip} is a strong tuple $S$ with increasing
markings, that is, $m(C_1) < m(C_2) < \dotsm < m(C_r)$. A strong
strip of size $1$ is the same thing as a strong cover. If we wish to
emphasize the inside and outside permutations then we use the
notation $\stc{w}{S}{u}$.

\begin{lem}
\label{L:markneq} Let $S=(C_1,\dotsc,C_r)$ be a strong tuple. Then
$m(C_k) \ne m(C_{k+1})$ for $1\le k < r$.
\end{lem}
\begin{proof} Let $S=(\dotsc,C^-,C,\dotsc)$ and
$C=(\stc{w}{i,j}{u})$. We have $m(C^-)=u^-(i^-)=w(i^-)$ and
$m(C)=w(j)$. So $m(C^-)=m(C)$ implies that $i^-=j$, which
contradicts the straddling inequalities $i^-\le l<j$.
\end{proof}

Given strong strips $S_1$ and $S_2$ such that $\out(S_1)=\ins(S_2)$
let $S_1\cup S_2$ denote the strong tuple obtained by concatenating
the sequences of strong covers defining $S_1$ and $S_2$. Then $S_1
\cup S_2$ is a strong strip if and only if one of the $S_i$ is
empty, or if both are nonempty and $m(\lc(S_1))<m(\fc(S_2))$.

A \textit{strong tableau} is a sequence $T=(S_1,S_2,\dotsc)$ of
strong strips $S_k$ such that $\out(S_k)=\ins(S_{k+1})$ for all
$i\in\Z_{>0}$ and $\size(S_k)=0$ for all sufficiently large $k$.
We define $\ins(T)=\ins(S_1)$ and $\out(T)=\out(S_k)$ for $k$
large. The \textit{weight} $\wt(T)$ of $T$ is the sequence
$$\wt(T)=(\size(S_1),\size(S_2),\dotsc).$$

We say that $T$ has \textit{shape} $u/v$ where $u=\out(T)$ and
$v=\ins(T)$.  If $T$ has shape $u/\id$ we simply say that $T$ has
shape $u$.

\section{Strong Schur functions}

\begin{definition}
For fixed $u,v\in \tS$, define the {\it strong Schur function}
\begin{align} \label{E:strongtabGF}
  \SSS_{u/v}(x) = \sum_T x^{\wt(T)}
\end{align}
where $T$ runs over the strong tableaux of shape $u/v$.
\end{definition}

We will use the convention that $\SSS_u(x) = \SSS_{u/\id}(x)$.  By
Proposition~\ref{P:strongGrass}, $\SSS_{u/v}(x) = 0$ if $v$ is
$l$-Grassmannian and $u$ is not.  We shall show later in
Theorem~\ref{T:kSchurStrong} that when $u$ is $l$-Grassmannian,
$\SSS_u(x)$ is a $k$-Schur function and thus possess remarkable
properties as shown in~\cite{LMproofs}.  However, for general $u, v
\in \tS$, the generating function $\SSS_{u/v}(x)$ is not well
understood, especially compared to the weak Schur functions to be
introduced in Chapter~\ref{ch:weaktableaux}. See
Section~\ref{sec:strongProperty}.

\begin{example}
\label{ex:strongc} Let $u = c_{l,m} := s_{l+m-1}\cdots s_{l+1}s_l$
for an integer $m$ satisfying $0 \leq m \leq n-1$.  For example, if
$l = 0$ then in window notation, $u =
[0,1,\ldots,m-1,m+1,m+2,\ldots,n-1,n+m]$. Let us calculate
$\SSS_u(x)$. The only strong cover $C = (\stc{w}{i,j}{u})$ where $w$
is $l$-Grassmannian, is given by $w = c_{l,m-1}$ and $(i,j) = (l,
l+m)$. This strong cover $C_m$ is marked at $m(C) = l+m$. Since
$\id\in\tS$ is $l$-Grassmannian for any $l$, by
Proposition~\ref{P:strongGrass} we see that a strong tableau $T =
(S_1,S_2,\ldots,S_r)$ with shape $u/\id$ is determined by specifying
integers $0 = m_0 \leq m_1 \leq m_2 \leq \cdots \leq m_r = m$ such
that $S_k = [c_{l,m_{k-1}};C_{m_{k-1}+1},\ldots,C_{m_k};c_{l,m_k}]$.
Thus $\SSS_u(x) = h_m(x)$.
\end{example}

\chapter{Weak Tableaux} \label{ch:weaktableaux}
\section{Cyclically decreasing permutations and weak tableaux} This
section follows \cite{LamAS}, which builds on earlier work in the
special case of $0$-Grassmannian elements (or $n$-cores)
\cite{LMcore,LMdual}. The intervals $I=\{a,a+1,\dotsc,b-1,b\}
\subsetneq \Z/n\Z$ considered in cyclical fashion will be denoted
with the interval notation $[a,b]$.

The \textit{left weak order} $\wle$ on $\tS$ (sometimes also called
the left weak Bruhat order) is defined by $w \wle v$ if and only if
there is a $u\in\tS$ such that $uw=v$ with
$\ell(u)+\ell(w)=\ell(v)$.

Given a proper cyclic interval $I=[a,b]\subsetneq\Z/n\Z$, let
$c_I=s_b s_{b-1}\dotsm s_a\in\tS$ be the product of the reflections
indexed by $I$, appearing in decreasing order. Given any proper
subset $A\subsetneq\Z/n\Z$ let $c_A = c_{I_1}\dotsm c_{I_t}$ where
$A=I_1\cup\dotsm\cup I_t$ is the decomposition of $A$ into maximal
cyclic intervals $I_k$ which are called the \textit{cyclic
components} of $A$. The element $c_A$ is well-defined since
$c_{I_i}$ and $c_{I_j}$ commute for $i\not=j$. Say that $c\in \tS$
is \textit{cyclically decreasing} if $c=c_A$ for some
$A\subsetneq\Z/n\Z$. Write $A(c)$ for this subset $A$.

\begin{example}
Let $n = 10$ and $A = \{0,1,3,4,6,9\}$. The cyclic components of
$A$ are $[9,1]$, $[3,4]$ and $[6,6]$ and we have $$c_A =
(s_1s_0s_9)(s_4s_3)(s_6) = (s_4s_3)(s_1s_0s_9)(s_6) = \dotsm.$$
Since $\ba{2}\not\in A$, the action of $c_A$ on $\Z$ acts on the
``window" $[3,12]\subset\Z$ by the cycles
$12\mapsto11\mapsto10\mapsto9\mapsto12$,
$5\mapsto4\mapsto3\mapsto5$, $7\mapsto6\mapsto7$, $8\mapsto8$, and
on all of $\Z$ periodically.
\end{example}

A \textit{weak strip} $W=(\ws{w}{}{v})$ consists of a pair
$w,v\in\tS$ such that $w \wle v$ and $vw^{-1}$ is cyclically
decreasing; it is a certain kind of interval in the left weak order.
If we wish to emphasize the cyclically decreasing element
$c=vw^{-1}$ then we write $W=(\ws{w}{A}{v})$ where $A=A(c)$. The
definition of $\wle$ implies that $\ell(c)+\ell(w)=\ell(v)$. The
\textit{size} of $W$, denoted $\size(W)$, is the integer
$\ell(v)-\ell(w)=\ell(c)=|A|$. We write $\ins(W)=w$ and $\out(W)=v$.

A \textit{weak tableau} $U$ is a sequence $U=(W_1,W_2,\dotsc)$ of
weak strips such that $\out(W_k)=\ins(W_{k+1})$ for all
$k\in\Z_{>0}$ and $\size(W_k)=0$ for $k$ large. Let
$\ins(U)=\ins(W_1)$ and $\out(U)=\out(W_k)$ for $k$ large. The
tableau $U$ gives precisely the data of a chain
$(\ins(W_1),\ins(W_2),\dotsc)$ in the left weak order on $\tS$ such
that consecutive elements define a weak strip. The \textit{weight}
$\wt(U)$ of a weak tableau $U$ is the composition
$\wt(U)=(\size(W_1),\size(W_2),\dotsc)$. We say that $U$ is a
tableau of \textit{shape} $u/v$ where $u=\out(U)$ and $v=\ins(U)$.
If $U$ has shape $u/\id$ we say that $U$ has shape $u$.

Note that unlike strong tableaux, weak tableaux do not depend on the
choice $l \in \Z$ of maximal parabolic.

\section{Weak Schur functions}
\begin{definition}
For $u,v\in\tS$ define the {\it weak Schur function} or {\it (skew)
affine Stanley symmetric function} by
\begin{align} \label{E:weakGF}
  \WS_{u/v}(x) = \sum_U x^{\wt(U)}
\end{align}
where $U$ runs over the weak tableaux of shape $u/v$.
\end{definition}

We use the shorthand $\WS_u(x)$ to mean $\WS_{u/\id}(x)$.  The
generating functions $\WS_{u/v}(x)$ were first introduced
in~\cite{LamAS} where they were called skew affine Stanley symmetric
functions, though weak tableaux were not explicitly used.  It is not
difficult to see that if $v \wle u$ then $\WS_{u/v}(x) =
\WS_{uv^{-1}}(x)$. Thus each weak tableaux generating function is in
fact an {\it affine Stanley symmetric function} (denoted $\tF_w(x)$
in~\cite{LamAS}).

If $u\in \tS$ is $0$-Grassmannian then weak tableaux of shape $u$
are the {\it $k$-tableaux} (with $k=n-1$) first defined
in~\cite{LMcore}. We shall translate strong and weak tableaux into
the context of $k$-tableaux in Chapter~\ref{ch:corestableaux}.  In
the case that $u \in \tS^0$, $\WS_{u}(x)$ is the dual $k$-Schur
function introduced in~\cite{LMdual} (called an affine Schur
function in~\cite{LamAS}).

The basic theorem for $\WS_{u/v}(x)$ is the following.
\begin{thm}[Symmetry of weak Schurs~\cite{LamAS}]
\label{T:weakSym} For $u, v \in \tS$, the generating function
$\WS_{u/v}(x)$ is a symmetric function in $x_1,x_2,\ldots$.
\end{thm}

We shall later also need the following properties of $\WS_{u/v}(x)$.
Let
\begin{align} \label{E:Dynkinflip}
  w &\mapsto w^*
\end{align}
be the unique involutive automorphism of the group $\tS$ such that
$s_i^* = s_{l-i}$ for all $i$.

\begin{thm}[Conjugacy of weak Schurs~\cite{LamAS}]
\label{T:weakconjugate} Let $w \in \tS$.  Then
$$
\omega^+(\WS_{w}(x)) = \WS_{w^{-1}}(x) = \WS_{w^*}(x).
$$
\end{thm}

\begin{thm}[Affine Schurs form a basis~\cite{LamAS, LMdual}]
\label{T:affineSchurbasis} The set $\{\WS_w(x) \mid w \in \tS^0\}$
forms a basis of $\La^{(n)}$.
\end{thm}

\begin{remark}
If $\sigma$ denotes the map on $\tS$ induced by a rotational
automorphism of the affine Dynkin diagram, then $\WS_w(x)=
\WS_{\sigma(w)}(x)$ (see also \cite[Propostion~18]{LamAS}). Thus
$\tS^0$ can be replaced by $\tS^l$ in
Theorem~\ref{T:affineSchurbasis}.
\end{remark}

\begin{example}
\label{ex:weakc} Let $w = c_{l,m} := s_{l+m-1}s_{l+m-2}\cdots s_{l}$
for $m \geq 0$ a non-negative integer. Then $w$ has a unique
reduced decomposition. Note that $c_{l,m}$ is always cyclically
decreasing if $m \leq n-1$ and that $c_{l,m}c_{l,m'}^{-1} =
c_{l+m',m-m'}$ if $m > m' \geq 1$.  The weak tableaux $U$ with shape $w$ are of the form $U =
(W_1,W_2,\ldots,W_r)$ where $W_k =
(\ws{c_{l,m_{k-1}}}{}{c_{l,m_k}})$ and $1 = m_0 \leq m_1 \leq m_2
\leq \cdots \leq m_k = m$ satisfies $m_k - m_{k-1} \leq n-1$ for $1
\leq k \leq r$.  The weight of $U$ is $\wt(U) =
(m_1-m_0,m_2-m_1,\ldots,m_r-m_{r-1})$.  The weak tableau generating
function, or affine Stanley symmetric function labeled by $w$ is
thus $\WS_w(x) = \tF_w(x)$, %
which is given by the image of $h_m(x)$ in $\La^{(n)}$.
\end{example}

\section{Properties of weak strips} For later use we collect some
properties of cyclically decreasing elements and weak strips. Let $A
\subsetneq \Z/n\Z$. We say that $i\in\Z$ is \textit{$A$-nice} if
$\ba{i-1}\notin A$ and that $i$ is \textit{$A$-bad} if $\ba{i}\notin
A$.

Given $i\in\Z$, let $j$ be the smallest $A$-nice integer
such that $j>i$.  We have
\begin{equation} \label{E:cvalue2}
  i+n > c_A(i) = \begin{cases} j-1 & \text{if $i$ is $A$-nice } \\
  i-1 &\text{otherwise.}
  \end{cases}
\end{equation}
In other words, if $i<j$ are consecutive $A$-nice integers for
$A\subsetneq\Z/n\Z$, then $j-i \le n$ and $c_A$ acts on the set of
integers $[i,j-1]$ by the cycle $j-1\mapsto j-2\mapsto\dotsm\mapsto
i+1\mapsto i\mapsto j-1$ and $\ba{i}$ is the cyclic minimum of the
cyclic component $[\ba{i},\ba{j-2}]$ of $A$. In particular
\begin{equation}
\label{E:cvalue}
\begin{alignedat}{2}
  c_A(i) &= i-1 &\qquad&\text{if $i$ is not $A$-nice} \\
  c_A(i) &= i &\qquad&\text{if $i$ is $A$-nice and $A$-bad} \\
  c_A(i) &> i &\qquad&\text{if $i$ is $A$-nice and not $A$-bad.}
\end{alignedat}
\end{equation}

We obtain the following Lemmata.

\begin{lem} \label{L:nicebad}
 Let $A\subsetneq \Z/n\Z$. Then $c_A$ restricts to
an order-preserving bijection from the set of $A$-nice integers to
the set of $A$-bad integers.
\end{lem}

\begin{lem} \label{L:cycdecrinv} Let $A\subsetneq\Z/n\Z$ and $i<j$.
Then $c_A(i)>c_A(j)$ if and only if $j\le c_A(i)$, $i$ is $A$-nice
and not $A$-bad. Equivalently, $c_A(i)>c_A(j)$ if and only if $i$ is $A$-nice and the
next $A$-nice integer after $i$ is larger than $j$.
In this case $j-i<n$ and $j$ is not $A$-nice.
\end{lem}

\begin{lem}
\label{L:weakcriterion} Let $c=c_A$ be cyclically decreasing and
$v=cw$. The following are equivalent:
\begin{enumerate}
\item $\ws{w}{A}{v}$ is a weak strip.
\item For every pair of consecutive  $A$-nice integers $a<b$ we have
$$w^{-1}(a)< \min(w^{-1}(a+1),w^{-1}(a+2),\dotsc,w^{-1}(b-1)).$$
\item For every pair of consecutive $A$-bad integers $a<b$ we have
$$v^{-1}(b) < \min(v^{-1}(a+1),\dotsc, v^{-1}(b-1)).$$
\end{enumerate}
\end{lem}

\begin{lem} \label{L:weaknoinv}
Let $\ws{w}{A}{v}$ be a weak strip. Then there do not exist integers
$i<j$ such that $w(i)>w(j)$ and $v(i)<v(j)$.
\end{lem}
\begin{proof}
Suppose $i$ and $j$ exist. By Lemma \ref{L:cycdecrinv} $w(j)$ is
$A$-nice and $w(j)<w(i)<b$ where $b$ is the smallest $A$-nice
integer greater than $w(j)$. Applying Lemma \ref{L:weakcriterion}(2)
we have the contradiction $j<i$.
\end{proof}

\begin{lem} \label{L:weakGrass}
Let $\ws{w}{A}{v}$ be a weak strip. If $v$ is $l$-Grassmannian, then $w$
is $l$-Grassmannian as well.
\end{lem}
\begin{proof}
Suppose $v(l+1)<v(l+2)< \cdots < v(l+n)$.  Then Lemma~\ref{L:weaknoinv}
implies that $w(l+1)<w(l+2)< \cdots < w(l+n)$.
\end{proof}

\begin{lem} \label{L:Aadd} Let $A\subsetneq\Z/n\Z$,
$q<p$ consecutive $A$-nice integers such that $\ba{q}\ne\ba{p}$ and
$B=A\cup\{\ba{p-1}\}$. Then
\begin{align} \label{E:changeA}
  c_B = c_A t_{qp}
\end{align}
\end{lem}
\begin{proof} Let $r>p$ be the next $A$-nice integer.
After canceling $c_I$ for each of the common cyclic components $I$
of $A$ and $B$, we may assume that $A=[\ba{q},\ba{p-2}]\cup
[\ba{p},\ba{r-2}]$ and $B=[\ba{q},\ba{r-2}]$. Then \eqref{E:changeA}
is equivalent to
\begin{align*}
 s_{p-1} (s_{p-2} \dotsm s_q) = (s_{p-2}\dotsm s_q) t_{qp}
\end{align*}
which holds by \eqref{E:tconj}.
\end{proof}

\section{Commutation of weak strips and strong covers}
\label{SS:comm}

An \textit{initial pair} $(W,S)$ consists of a weak strip $W$ and
strong strip $S$ with $\ins(W)=\ins(S)$. A special case is an
initial pair $(W,C)$ where $C$ is a marked strong cover. A
\textit{final pair} $(W',S')$ consists of a weak strip $W'$ and
strong strip $S'$ such that $\out(W')=\out(S')$. Since a marked
strong cover $C'$ is a special case of a strong strip we can refer
to a final pair $(W',C')$.
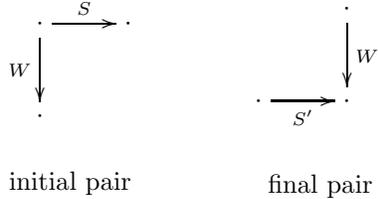
\begin{align*}
\begin{matrix}
\xymatrix{%
{\cdot} \ar[r]^S \ar[d]_W & {\cdot} \\
{\cdot}} \\ \\ \text{initial pair}
\end{matrix}\qquad\qquad
\begin{matrix}
\xymatrix{%
{} & {\cdot} \ar[d]^{W'} \\
{\cdot} \ar[r]_{S'} & {\cdot}} \\ \\ \text{final pair}
\end{matrix}
\end{align*}

Let $(W,C)=(\wci{w}{A}{v}{i,j}{u})$ be an initial pair. By
associativity in $\tS$ one always has a commutative diagram
\begin{align*}
\begin{CD}
w @>{i,j}>> u \\
@V{A}VV @VV{A}V \\
v @>>{i,j}> x
\end{CD}
\end{align*}
where $x=v\,t_{ij}=c_A u$, in which the arrows labeled by $A$
represent left multiplication by $c_A$ and those labeled by $i,j$
represent right multiplication by $t_{ij}$. Commutativity of such a
diagram means that the corresponding products of elements in $\tS$
result in the same element.

We say that the initial pair $(W,C)$ \textit{commutes} (or that $W$
and $C$ commute, or that $c_A$ and $C$ commute, and so on) if $v<x$,
that is, $v(i)<v(j)$ or equivalently $x(i)>x(j)$. This should not be
confused with the commutation of a diagram.

Let $(W',C')=(\wcf{u}{A'}{x}{v}{a,b})$ be a final pair. We say that
$(W',C')$ commutes if $w < u$ (that is, $w(a)<w(b)$ or equivalently
$u(a)>u(b)$) where $w=u\,t_{ab}=c_{A'}^{-1}v$.

\begin{lem} \label{L:initfincommute}
Let $w,u,v,x\in\tS$, $A\subsetneq \Z/n\Z$ and $t_{ij}$ be
such that $u=w\,t_{ij}$, $v=c_A w$, and $x=v\,t_{ij}=c_A u$. Then
the following are equivalent:
\begin{enumerate}
\item $(\wci{w}{A}{v}{i,j}{u})$ is a commuting initial
pair.
\item $(\wcf{u}{A}{x}{v}{i,j})$ is a commuting final
pair.
\end{enumerate}
\end{lem}
\begin{proof} For the forward direction we have
\begin{align*}
|A|+\ell(u) \ge \ell(x) \ge \ell(v)+1 =|A|+\ell(w)+1=|A|+\ell(u)
\end{align*}
since $c_A u=x$, $v<x$, $\ws{w}{A}{v}$ is a weak strip, and
$w\lessdot u$. All inequalities must be equalities, which proves
that $\ws{u}{A}{x}$ is a weak strip and $v\lessdot x$. This proves
the forward direction. The reverse direction is similar.
\end{proof}

We list some Lemmata regarding noncommuting initial pairs.

\begin{lem}
\label{L:Fcyclicobstruction} Let $(W,C)=(\wci{w}{A}{v}{i,j}{\cdot})$
be an initial pair and $c=c_A$.
\begin{enumerate}
\item The following are equivalent: (i) $(W,C)$ does not commute;
(ii) $c(w(i))>c(w(j))$; (iii) $w(j)\le c(w(i))$, $\ba{w(i)}\in A$
and $w(i)$ is $A$-nice; (iv) $w(i)$ is $A$-nice and the next $A$-nice integer after
$w(i)$ is larger than $w(j)$.
\item If $(W,C)$ does not commute then
$0 < w(j)-w(i)<n$ and $w(j)$ is not $A$-nice.
\end{enumerate}
\end{lem}
\begin{proof} By definition, $(W,C)$ commutes if and only if $v(i)<v(j)$.
The lemma then follows from Lemma~\ref{L:cycdecrinv}.
\end{proof}

\begin{lem}
\label{L:Fremoveresidue} Let $(\wci{w}{A}{v}{i,j}{u})$ be a
noncommuting initial pair and $\Ah = A-\{\ba{w(j)-1}\}$. Then
$\ws{u}{\Ah}{v}$ is a weak strip.
\end{lem}
\begin{proof} By a length calculation it suffices to show that
$c_{\Ah}  u = v$. By Lemma \ref{L:Fcyclicobstruction}(iv),
$w(i)$ and $w(j)$ are consecutive $A^\vee$-nice integers.
Applying Lemma~\ref{L:Aadd} with the consecutive $\Ah$-nice integers $w(i)<w(j)$, we
get $c_{A^\vee}= c_A t_{w(i) w(j)}=c_A t_{u(i) u(j)}$ which easily gives the result.
\end{proof}

The picture for Lemma \ref{L:Fremoveresidue} is given below.
\begin{align*}
\xymatrix{%
{w} \ar[r]^{i,j} \ar[d]_{A} & {u} \ar@{-->}[dl]^{\Ah} \\
{v} & {}
}%
\end{align*}

\begin{lem}
\label{L:Faddresidue} Let $\ws{u}{\Ah}{v}$ be a weak strip with
$|\Ah|<n-1$, $q < p$ consecutive $\Ah$-nice integers such that
$u^{-1}(q) < u^{-1}(p)$, $A'=\Ah\cup\{\ba{p-1}\}$, and $x=c_{A'}u$.
Then
\begin{equation}\label{E:weakadd}
  x = c_{A'}\, u =  v\, t_{u^{-1}(q),u^{-1}(p)},
\end{equation}
$\ws{u}{A'}{x}$ is a weak strip, and $v\lessdot x$ is a strong
cover. If in addition $u^{-1}(q) \le l < u^{-1}(p)$ then
$(\wcf{u}{A'}{x}{v}{a,b})$ is a noncommuting final pair where
$(a,b)=(u^{-1}(q),u^{-1}(p))$.
\end{lem}
\begin{proof} Let $r$ be the next $\Ah$-nice integer after $p$.
Using the hypothesis $u^{-1}(q)<u^{-1}(p)$ and applying Lemma
\ref{L:weakcriterion} to the pairs of consecutive $\Ah$-nice
integers $q<p$ and $p<r$, and to the pair of consecutive $A'$-nice
integers $q<r$, we conclude that $\ws{u}{}{c_{A'}u=x}$ is a weak
strip. Equation \eqref{E:weakadd} follows from Lemma \ref{L:Aadd}.
The strong cover assertion follows from \eqref{E:weakadd} and a
length computation. Noncommutativity holds since $q<p$.
\end{proof}

\begin{lem} \label{L:Faddresidue2} Let $\ws{u}{\Ah}{v}$ be a weak
strip with $|\Ah|<n-1$, $q < p$ consecutive $\Ah$-bad integers such
that $v^{-1}(q) < v^{-1}(p)$, $A' = \Ah \cup \{\ba{q}\}$, and
$x=c_{A'}u$. Then
\begin{equation}\label{E:weakadd2}
x = c_{A'} \, u = v \, t_{v^{-1}(q),v^{-1}(p)},
\end{equation}
$\ws{u}{A'}{x}$ is a weak strip, and $v \lessdot x$ is a strong
cover. Moreover if $v^{-1}(q)\le l < v^{-1}(p)$ then
$(\wcf{u}{A'}{x}{v}{a,b})$ is a noncommuting final pair where
$(a,b)=(v^{-1}(q),v^{-1}(p))$.
\end{lem}
\begin{proof}By Lemma~\ref{L:nicebad} and \eqref{E:cvalue2}, $q_1=c_{\Ah}^{-1}(q)$ and $p_1=c_{\Ah}^{-1}(p)$ are
consecutive $\Ah$-nice integers with $p_1-1=q$. We have
$v^{-1}(q)=u^{-1}(c_{\Ah}^{-1}(q))=u^{-1}(q_1)$ and
$v^{-1}(p)=u^{-1}(p_1)$. The Lemma follows by an application of
Lemma \ref{L:Faddresidue} to the weak strip $\ws{u}{\Ah}{v}$ and the
$\Ah$-nice integers $q_1<p_1$.
\end{proof}

The diagram for Lemmata \ref{L:Faddresidue} and \ref{L:Faddresidue2}
is given as follows.
\begin{align*}
\xymatrix{%
{} & {u} \ar[dl]_{\Ah} \ar@{-->}[d]^{A'} \\
{v} \ar@{-->}[r]_{a,b} & {x}
}%
\end{align*}

We now give the corresponding lemmata for final pairs.

\begin{lem}
\label{L:Rcyclicobstruction} Let $(W',C')=(\wcf{u}{A'}{x}{v}{a,b})$
be a final pair.
\begin{enumerate}
\item $(W',C')$ does not commute if and only if $u(a)<u(b)$. \item
If $(W',C')$ does not commute then $\ba{u(a)}\in A'$, $u(a)$ is
$A'$-nice, $u(b)$ is not $A'$-nice, $c_{A'}u(a) \ge u(b)$, and $0
< u(b)-u(a)<n$.  In particular, the next $A$-nice integer after
$u(a)$ is in that case larger than $u(b)$.
\end{enumerate}
\end{lem}
\begin{proof} The proof is similar to that of Lemma
\ref{L:Fcyclicobstruction}.
\end{proof}

\begin{lem}
\label{L:Rremoveresidue} Let $(\wcf{u}{A'}{x}{v}{a,b})$ be a
noncommuting final pair and $\Ah = A'-\{\ba{u(b)-1}\}$. Then
$\ws{u}{\Ah}{v}$ is a weak strip.
\end{lem}
\begin{proof}  By a length calculation it suffices to show that $c_{\Ah}  u =
v$. By Lemma \ref{L:Rcyclicobstruction},
$u(a)<u(b)$ are consecutive $\Ah$-nice integers. The result then follows by an application of
Lemma~\ref{L:Aadd} with the consecutive $\Ah$-nice integers $u(a)<u(b)$.
\end{proof}

The picture for Lemma \ref{L:Rremoveresidue} is given below.
\begin{align*}
\xymatrix{%
{} & {u} \ar@{-->}[dl]_{\Ah} \ar[d]^{A'} \\
{v} \ar[r]_{a,b} & {x}
}%
\end{align*}

\begin{lem} \label{L:Raddresidue} Let $\ws{u}{\Ah}{v}$ be a weak
strip with $|\Ah|<n-1$ and $p<q$ consecutive $\Ah$-nice integers
such that $u^{-1}(q)<u^{-1}(p)$, $A=\Ah\cup\{\overline{q-1}\}$, and
$w=c_A^{-1}v$. Then
\begin{align*}
  w = c_A^{-1} v = u t_{u^{-1}(q),u^{-1}(p)},
\end{align*}
$\ws{w}{A}{v}$ is a weak strip, and $\stc{w}{i,j}{u}$ is a strong
strip where $(i,j)=(u^{-1}(q),u^{-1}(p))$. If also $u^{-1}(q)\le
l<u^{-1}(p)$ then $(\wci{w}{A}{v}{i,j}{u})$ is a noncommuting
initial pair.
\end{lem}
\begin{proof} The proof of the weak and strong strip properties is
similar to the proof of Lemma \ref{L:Faddresidue}. The noncommuting
property is equivalent to $c_A(p)>c_A(q)$, which holds by the
assumptions on $p$ and $q$.
\end{proof}

The diagram for Lemma \ref{L:Raddresidue} is given below.
\begin{align*}
\xymatrix{%
{w} \ar@{-->}[r]^{i,j} \ar@{-->}[d]_{A} & {u} \ar[dl]^{\Ah} \\
{v} & {}
}%
\end{align*}

\chapter{Affine Insertion and Affine Pieri} \label{ch:main}

We employ Fomin's general method of growth diagrams \cite{Fom} to
help define a bijection that we call affine insertion. We deduce our
main results from the properties of affine insertion.

\section{The local rule $\phi_{u,v}$} For $u,v\in\tS$, let
$\IIL_{u,v}$ be the set of triples $(W,S,e)$ where $(W,S)$ is an
initial pair with $\out(W)=v$ and $\out(S)=u$, and $e\in\Z_{\ge0}$
is such that $\size(W)+e<n$. Let $\OOL_{u,v}$ be the set of final
pairs $(W',S')$ such that $\ins(W')=u$ and $\ins(S')=v$.

\begin{prop} \label{P:local}
For each $u,v\in \tS$ there is a bijection
\begin{align*}
  \phi_{u,v}: \IIL_{u,v} &\rightarrow \OOL_{u,v} \\
  (W,S,e) &\mapsto (W',S')
\end{align*}
such that
\begin{equation} \label{E:locallength}
  \size(S')=\size(S)+e.
\end{equation}
\end{prop}
The weak and strong strips form the edges of a diagram
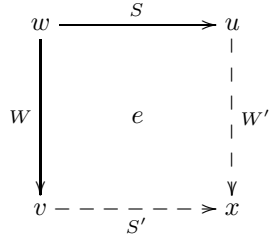
\begin{align*}
\xymatrix{ %
{w} \ar[rr]^{S} \ar[dd]_{W} & {} & {u} \ar@{-->}[dd]^{W'} \\
& {e} & \\
 {v} \ar@{-->}[rr]_{S'} & & {x}
}
\end{align*}
We call $\phi_{u,v}$ the local rule because of its role in the
affine insertion bijection defined in the following section. The
formulation of the local rule and the proof of its bijectivity form
the core of our arguments and occupy Chapters \ref{ch:forward}
through \ref{ch:bijectivity}.

\section{The affine insertion bijection $\Phi_{u,v}$} A {\it
$n$-bounded matrix} is a matrix $m=(m_{ij})_{i,j>0}$ with
nonnegative integer entries, only finitely many of which are
nonzero, all of whose row sums are strictly less than $n$.  We write
$\rowsums(m)$ and $\colsums(m)$ to indicate the sequences of
integers given by the row sums and column sums of $m$ respectively.
Let $\M$ denote the set of $n$-bounded matrices.

Fix $u,v\in \tS$. Let $\II_{u,v}$ be the set of triples $(T,U,m)$
where $T$ is a skew strong tableau, $U$ is a skew weak tableau, and
$m\in\M$, such that
\begin{align*}
  \ins(T) &= \ins(U) \\
  \out(T) &= u \\
  \out(U) &= v \\
  \wt(U) + \rowsums(m) &\le (n-1,n-1,\dotsc)
\end{align*}
When $u=v=\id\in\tS$ then $\II_{\id,\id} \cong \M$ since $T$ and $U$
must respectively be the empty strong strip and weak strip from
$\id$ to $\id$.

Let $\OO_{u,v}$ be the set of pairs $(P,Q)$ where $P$ is a skew strong
tableau and $Q$ a skew weak tableau with $\ins(P)=v$, $\ins(Q)=u$, and
$\out(P)=\out(Q)$.

The following is our main theorem. The reduction of its proof to
that of Proposition \ref{P:local}, is an instance of Fomin's theory
\cite{Fom}.

\begin{thm} \label{T:main}
There is a bijection
\begin{equation}
\label{E:affins}
\begin{split}
  \Phi_{u,v}: \II_{u,v}  &\rightarrow \OO_{u,v} \\
    (T,U,m) &\mapsto (P,Q)
\end{split}
\end{equation}
such that %
\begin{align}
\label{E:cols} \wt(T)+\colsums(m)&=\wt(P) \\ %
\label{E:rows} \wt(U)+\rowsums(m)&=\wt(Q).
\end{align}
\end{thm}
We picture the input and output of the bijection by the diagram
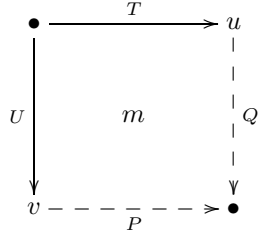
\begin{align*}
\xymatrix{ %
{\bullet} \ar[rr]^{T} \ar[dd]_{U} & {} & {u} \ar@{-->}[dd]^{Q} \\
& {m} & \\
 {v} \ar@{-->}[rr]_{P} & & {\bullet}
}
\end{align*}
where solid and dotted arrows indicate input and output data
respectively.
\begin{proof}[Proof of Theorem \ref{T:main}]
For each $u',v'\in\tS$, Proposition \ref{P:local} specifies a
bijection $\phi_{u',v'}:\IIL_{u',v'} \rightarrow \OOL_{u',v'}$. Let
$(T,U,m)\in\II_{u,v}$. The growth diagram of $(T,U,m)$ (defined by
the local rules $\phi_{u',v'}$) is by definition the directed graph
with vertices $G_{ij}\in \tS$ for $i,j\in\Z_{\ge0}$ indexed
matrix-style, with horizontal edges $G_{i,j-1}\to G_{i,j}$ endowed
with the structure of strong strips, and vertical edges
$G_{i-1,j}\to G_{i,j}$ given by weak strips, such that
\begin{enumerate}
\item[(G1)] The zero-th row of $G$ is the strong tableau $T$.
\item[(G2)] The zero-th column of $G$ is the weak tableau $U$.
\item[(G3)] For every $(i,j)\in \Z_{>0}^2$, the two-by-two subgraph
\begin{align*}
\xymatrix{ %
{G_{i-1,j-1}} \ar[rr]^{S} \ar[dd]_{W} & {} & {G_{i-1,j}} \ar[dd]^{W'} \\
& {m_{ij}} & \\
 {G_{i,j-1}} \ar[rr]_{S'} & & {G_{i,j}}
}
\end{align*}
satisfies $$\phi_{u',v'}(W,S,m_{ij})=(W',S')$$ where $u'=G_{i-1,j}$
and $v'=G_{i,j-1}$.
\end{enumerate}
The north and west boundaries of the growth diagram are specified by
(G1) and (G2) above. Here by convention a strong (resp. weak)
tableau is defined by an infinite sequence of strong (resp. weak)
strips, but only finitely many of these strips have positive size.

If $(i,j)\in\Z_{>0}^2$, by induction we may assume that the part of
the growth diagram northwest of $G_{i,j}$, is uniquely defined. Then
by (G3), $G_{i,j}$ and its incoming edges are specified by
$\phi_{u',v'}$. It follows that the growth diagram is well-defined.

The growth diagram satisfies two additional properties.
\begin{enumerate}
\item Let $G_{i,\bullet}$ denote the strong tableau given by the
$i$-th row of $G$. Then for $i\gg0$, $G_{i,\bullet}$ stabilizes;
call the limiting tableau $G_{\infty,\bullet}$.
\item Let $G_{\bullet,j}$ denote the weak tableau given by the
$j$-th column of $G$. Then for $j\gg0$, $G_{\bullet,j}$ stabilizes;
call the limiting tableau $G_{\bullet,\infty}$.
\end{enumerate}
To see this, there is an $N\gg0$ such that $U$ stabilizes after $N$
steps and $m_{ij}=0$ for $i\ge N$. By \eqref{E:locallength} every
column $G_{\bullet,j}$ stabilizes after $N$ steps. In other words,
for all $i\ge N$, the $i$-th row $G_{i,\bullet}$ is the same strong
tableau. This proves the existence of $G_{\infty,\bullet}$. In a
similar manner one may show that $G_{\bullet,\infty}$ exists.

The affine insertion bijection is defined by
\begin{align} \label{E:Phiuvdef}
  \Phi_{u,v}(T,U,m) = (P,Q) := (G_{\infty,\bullet},G_{\bullet,\infty}).
\end{align}
By construction $(P,Q)\in \OO_{u,v}$ and $\Phi_{u,v}$ is
well-defined.

To show that $\Phi_{u,v}$ is a bijection we define the inverse map
$\Psi_{u,v}$. Given $(P,Q)\in\OO_{u,v}$, let $x=\out(P)=\out(Q)$.
Let $N_1$ and $N_2$ be sufficiently large such that for all $i\ge
N_1$ and $j\ge N_2$, the $i$-th element of $P$ and the $j$-th
element of $Q$ are equal to $x$. We define a growth diagram $G$ as
follows. We set $G_{i,\bullet}=P$ for $i\ge N_1$ and
$G_{\bullet,j}=Q$ for $j\ge N_2$. This is consistent: these
definitions overlap in the region $i\ge N_1$ and $j\ge N_2$, where
the entries $G_{i,j}$ are all equal to $x$ and all strong and weak
strips are empty. In the middle of each two-by-two subdiagram in
this region we place the integer $0$. We now use the inverse
$\psi_{u',v'}$ of the local rule $\phi_{u',v'}$ to fill in each
two-by-two subdiagram of $G$ (including the ``excitation integers"
in the middle of the subdiagram) given its south and east borders.
Then all of $G$ may be computed, as well as a matrix $m$. Letting
$T=G_{0,\bullet}$ and $U=G_{\bullet,0}$ be the north and west
boundaries of $G$, we define $\Psi(P,Q)=(T,U,m)$. It is easy to show
that $\Psi_{u,v}\circ \Phi_{u,v}=\id_{\II_{u,v}}$ and
$\Phi_{u,v}\circ\Psi_{u,v} = \id_{\OO_{u,v}}$.
\end{proof}

As a special case of affine insertion (and using
Proposition~\ref{P:strongGrass}) we obtain an RSK bijection for the
affine Grassmannian.

\begin{thm} \label{T:affGrRSK} For $l=0$,
$\Phi_{\id,\id}$ gives a bijection $\M\rightarrow \OO_{\id,\id}$
from $n$-bounded matrices to pairs $(P,Q)$ of strong and weak
tableaux from $\id$ to a common $0$-Grassmannian element
$\out(P)=\out(Q)$.
\end{thm}

We shall prove later (Theorem~\ref{T:limitRSK}) that when
$n\to\infty$ the affine insertion bijection $\Phi_{\id,\id}$
coincides with the classical row insertion RSK correspondence.

\subsection{Cauchy identity and Pieri rules for strong and weak tableaux}
Define the {\it affine Cauchy kernel} $\Omega_n(x,y)$ by
\begin{align*}
\Omega_n(x,y) &= \prod_{i} \left(1 + y_ih_1(x) + y_i^2 h_2(x) +
\cdots + y_i^{n-1}h_{n-1}(x)\right) \\&= \sum_{\lambda \,:\,
\lambda_1 < n} h_\lambda(x) m_\lambda(y).\end{align*} It is an
element of a completion $\Lambda_{(n)}(x) \hat{\otimes}
\Lambda^{(n)}(y)$ of $\La_{(n)}(x) \otimes \La^{(n)}(y)$.

The following is an immediate enumerative consequence of
Theorem~\ref{T:main}.
\begin{thm}[Generalized Affine Cauchy Identity] \label{T:genCauchy}
Let $u,v\in\tS$. Then the following identity holds in the quotient
$\Z[[x_1,x_2,\ldots,y_1,y_2,\ldots]]/(y_1^n,y_2^n,\ldots)$ of the
formal power series ring in two infinite sets of variables:
$$
\Omega_n(x,y)\,\sum_{w \in \tS} \SSS_{u/w}(x)\WS_{v/w}(y)  = \sum_{z
\in \tS} \SSS_{z/v}(x)\WS_{z/u}(y).
$$
\end{thm}
\begin{cor}\label{C:skewCauchy}
Let $u,v\in\tS^0$. Then the following identity holds in the quotient
of the formal power series ring
$\Z[[x_1,x_2,\ldots,y_1,y_2,\ldots]]/(y_1^n,y_2^n,\ldots)$:
$$
\Omega_n(x,y)\,\sum_{w \in \tS^0} \SSS_{u/w}(x)\WS_{v/w}(y)  =
\sum_{z \in \tS^0} \SSS_{z/v}(x)\WS_{z/u}(y).
$$
\end{cor}
\begin{proof}
In Theorem~\ref{T:genCauchy}, use Proposition~\ref{P:strongGrass}
and Lemma~\ref{L:weakGrass} to deduce that $\SSS_{z/v}(x) = 0$
unless $z \in \tS^0$ and $\WS_{v/w}(y) = 0$ unless $w \in \tS^0$.
\end{proof}

\begin{cor}[Affine Cauchy Identity] \label{C:cauchy} The following identity holds
in the formal power series ring
$\Z[[x_1,x_2,\ldots,y_1,y_2,\ldots]]$:
$$
\Omega_n(x,y) = \sum_{w \in \tS^0} \SSS_{w}(x)\WS_{w}(y).
$$
\end{cor}
\begin{proof}
Put $u=v={\rm id}$ in Corollary~\ref{C:skewCauchy}.
\end{proof}

We now deduce two combinatorial Pieri rules from
Theorem~\ref{T:main}. Recall from Example \ref{ex:strongc} that
$h_r(x)=\SSS_{c_{l,r}}(x)$.

\begin{thm}[Strong Pieri rule]
\label{T:strongPieri} Let $u \in \tS$ and $1 \leq r \leq n-1$.  Then
$$
h_r(x) \, \SSS_u(x) = \sum_{\ws{u}{}{z}} \SSS_z(x),
$$
where the summation is over weak strips $\ws{u}{}{z}$ of size $r$.
\end{thm}
Note that by Proposition~\ref{P:strongGrass}, $\SSS_u(x) = 0$ unless
$u$ is $l$-Grassmannian so all the permutations in the theorem can
be taken to be $l$-Grassmannian.

\begin{proof}
In Theorem~\ref{T:main}, set $v = \id$ and restrict the bijection to
triples $(T,U,m)$ such that $m$ has non-zero entries only in the
first row and such that the entries in the first row of $m$
sum to $r$, where $1 \leq r \leq n-1$.  Since $v = \id$, we must
have $\ins(T) = \ins(U) = \id$ as well, so that in effect
$\Phi_{u,v}$ restricts to a bijection from pairs $(T,m')$ to pairs
$(P,Q)$ where $T$ is a strong tableau of shape $u$, the infinite
vector $m'$ given by $m'_j = m_{1j}$ has non-negative integer
entries summing to $r$, and $Q$ is a weak strip of size $r$ with
$\ins(Q) = u$, and finally $P$ is a strong tableau with shape $z =
\out(P) = \out(Q)$.  Now note that the weight generating function of
the vectors $m'$ is $h_r(x)$.  In view of \eqref{E:cols}, taking the
strong tableau generating functions for the input and output of the
bijection gives our theorem.
\end{proof}

One may define {\it dual weak strips} $\dws{w}{}{v}$ by replacing
cyclically decreasing permutations with cyclically increasing
permutations (defined in the obvious way).

\begin{thm}[Dual strong Pieri rule]
\label{T:dualstrongPieri} Let $u \in \tS$ and $1 \leq r \leq n-1$.
Then
$$
e_r(x) \, \SSS_u(x) = \sum_{\dws{u}{}{z}} \SSS_z(x),
$$
where the sum runs over dual weak strips $\dws{u}{}{z}$ of size $r$.
\end{thm}
\begin{proof}
In Theorem~\ref{T:kSchurStrong} we show that $\SSS_u(x)$ is a
$k$-Schur function when $u \in \tS^l$ and by
Proposition~\ref{P:strongGrass}, $\SSS_u(x) = 0$ otherwise.  It is
known from~\cite{LMproofs} that $k$-Schur functions behave well
under the involution $\omega$ of $\La_{(n)}$ and in our notation we
have $\omega(\SSS_u(x)) = \SSS_{u^*}(x)$ where $u^*$ is defined in
\eqref{E:Dynkinflip}. The involution $u \mapsto u^*$ also
interchanges cyclically increasing and cyclically decreasing
permutations.  Thus the dual strong Pieri rule follows from applying
$\omega$ to Theorem~\ref{T:strongPieri}.
\end{proof}

Recall from Example \ref{ex:weakc} that in $\Lambda^{(n)}$,
$h_r(x)=\WS_{c_{l,r}}(x)$.

\begin{thm}[Weak Pieri rule]
\label{T:weakPieri} Let $w \in \tS$ and $1 \leq r$.  Then the
following identity holds in $\Lambda^{(n)}$:
$$
h_r(x) \, \WS_w(x) = \sum_{S} \WS_{\out(S)}(x),
$$
where the sum runs over strong strips $S$ of size $r$ such that
$\ins(S) = w$.
\end{thm}
\begin{proof}
The proof is similar to the proof of Theorem~\ref{T:strongPieri}. We
take $u = \id$, and restrict to matrices $m$ with non-zero entries
only in the first column.
\end{proof}

For the weak Pieri rule, we can use Theorem~\ref{T:weakconjugate} to
give a dual rule.
\begin{thm}[Dual weak Pieri rule]
\label{T:dualweakPieri} Let $w \in \tS$ and $1 \leq r$.  Then the
following identity holds in $\Lambda^{(n)}$:
$$
e_r(x) \, \WS_w(x) = \sum_{S} \WS_{\out(S)^{-1}}(x),
$$
where the sum runs over strong strips $S$ of size $r$ such that
$\ins(S) = w^{-1}$.
\end{thm}
One can replace $w^{-1}$ and $\out(S)^{-1}$ above by $w^*$ and
$\out(S)^*$, where $w^*$ is defined in \eqref{E:Dynkinflip}.

\begin{proof}
Clearly $\omega^+(h_r) = e_r$ as elements of $\La^{(n)}$.  The
theorem follows from applying $\omega^+$ to
Theorem~\ref{T:weakPieri} and using Theorem~\ref{T:weakconjugate}.
\end{proof}

\section{Pieri rules for the affine Grassmannian} For this section,
we will assume that $l = 0$ and write ``Grassmannian'' instead of
$0$-Grassmannian.  We will use Theorem~\ref{T:Lam} to deduce two
Pieri rules for the affine Grassmannian $\Gr$.

\begin{thm}[Monomial expansion of a $k$-Schur] \label{T:kSchurStrong}
Let $u\in\tS^0$ be Grassmannian.  Then $\SSS_{u}(x)$ coincides with
the $k$-Schur function $s_u^{(k)}(x)$ for $k=n-1$.
\end{thm}

We conjecture in Conjecture~\ref{conj:kSchurt} that the $k$-Schur
functions depending on a parameter $t$ can also be expressed as
generating functions of strong tableaux, using an additional
statistic called spin.

\begin{proof}
The $k$-Schur functions $s_u^{(k)}(x)$ are defined in~\cite{LMproofs} as
the symmetric functions satisfying a certain Pieri rule.
Equivalently, one may define $\{s_u^{(k)}(x) \mid u \in \tS^0\}$ as
the basis of $\Lambda_{(n)}$ dual to the affine Schur basis
$\{\tF_u(x) \mid u \in \tS^0\}$ of $\Lambda^{(n)}$
see~\cite{LMdual,LamAS}.  By definition we have $\WS_{w}(x) = \tF_w(x)$.
By Corollary~\ref{C:cauchy} it suffices to show that $\SSS_{w}(x)
\in \Lambda_{(n)}$, for then the duality (and the fact that $\{
\SSS_w(x)\mid w\in \tS^0\}$ forms a basis of $\La_{(n)}$) will
follow from an argument similar to~\cite[(4.6)]{Mac}.

To show that $\SSS_w(x) \in \Lambda$ we let $\sigma_i$ be the
ring-involution of the ring of formal power series in
$x_1,x_2,\ldots,$ which interchanges $x_i$ and $x_{i+1}$.  Then by
Corollary~\ref{C:cauchy},
\begin{align*}
&\sum_{w \in \tS^0} \left(\sigma_i \cdot
\SSS_{w}(x)\right)\WS_{w}(y)
\\ &= \sigma_i \cdot \prod_{i} \left(1 + y_ih_1(x) + y_i^2 h_2(x) +
\cdots + y_i^{n-1}h_{n-1}(x)\right)\\
 &= \prod_{i} \left(1 + y_ih_1(x) + y_i^2 h_2(x) + \cdots +
y_i^{n-1}h_{n-1}(x)\right) \\
&= \sum_{w \in \tS^0} \SSS_{w}(x)\WS_{w}(y).
\end{align*}
By Theorem~\ref{T:affineSchurbasis}, the $\WS_{w}(x)$ are linearly
independent elements of $\Lambda^{(n)}$. Taking the coefficient of
$\WS_{w}(y)$ in the above equation we obtain $$\SSS_w(x) = \sigma_i
\cdot \SSS_{w}(x).$$ Since this holds for all $i$, we have
$\SSS_w(x) \in \Lambda$. Finally, $\{\WS_{w}(x) \mid w \in \tS^0\}$
is independent and no terms $h_r(x)$ for $r > n-1$ occurs in the
affine Cauchy kernel; so $\SSS_w(x) \in \Lambda_{(n)}$.
\end{proof}
Thus by Theorem~\ref{T:Lam}, the generating functions $\SSS_w(x)$
are explicit combinatorial representatives of $\xi_w$.

Recall from Examples~\ref{ex:strongc} and~\ref{ex:weakc} that
$c_{l,m}:= s_{l+m-1}\cdots s_{l+1} s_l$.

\begin{thm}[Pieri rule for $H_*(\Gr)$]
\label{T:homPieri} Let $w \in \tS^0$ and $1 \leq m \leq n-1$. Then
$$
\xi_{c_{0,m}}\, \xi_w = \sum_{\ws{w}{}{z}} \xi_z,
$$
where the sum runs over weak strips $\ws{w}{}{z}$ of size $m$.
\end{thm}
\begin{proof}
Example~\ref{ex:strongc} says that $\SSS_{c_{0,m}}(x) = h_m(x)$. The
theorem then follows immediately from Theorems~\ref{T:strongPieri},
\ref{T:kSchurStrong} and~\ref{T:Lam}.
\end{proof}

The Pieri rule for $k$-Schur functions was first stated
in~\cite{LMproofs} and in the notation here in~\cite{LMW}. Combining
this with the geometric identification in~\cite{Lam}, one can obtain
Theorem~\ref{T:homPieri} directly as a corollary.

\begin{thm}[Pieri rule for $H^*(\Gr)$]
\label{T:cohomPieri} Let $w \in \tS^0$ and $1 \leq m$.  Then
$$
\xi^{c_{0,m}}\, \xi^w = \sum_{S} \xi^{\out(S)},
$$
where the sum runs over strong strips $S$ of size $m$ such that
$\ins(S) = w$.
\end{thm}
\begin{proof}
Example~\ref{ex:weakc} says that $\WS_{c_{0,m}}(x) = h_m(x)$ in
$\Lambda^{(n)}$.  The theorem then follows immediately from
Theorems~\ref{T:weakPieri} and~\ref{T:Lam}.
\end{proof}
It is not difficult to see that both the weak and strong Pieri rules
reduce to the classical Pieri rule for the finite Grassmannian when
$\ell(w) + m < n$. Using the involution $w\mapsto w^*$ of
\eqref{E:Dynkinflip}, Theorems~\ref{T:dualstrongPieri}
and~\ref{T:dualweakPieri} also gives us a rule for multiplication by
$\xi^{c_{0,m}^*}$ in $H_*(\Gr)$ and $H^*(\Gr)$ (note that if $w \in
\tS^0$ then $w^* \in \tS^0$).

\begin{prop}
Theorems~\ref{T:homPieri} and~\ref{T:cohomPieri} determine the
multiplicative structures of $H_*(\Gr)$ and $H^*(\Gr)$.
\end{prop}
\begin{proof}
It suffices to show that the Pieri rules can be inverted, so that
$\xi^w$ (or $\xi_w$) for any $w \in \tS^0$ can be written in terms
of $\{\xi^{c_{0,m}}\}$ (or $\{\xi_{c_{0,m}}\}$).  The theorem will
then follow from the fact that $\{\xi^w \mid w \in \tS^0\}$ (or
$\{\xi_w \mid w \in \tS^0\}$) forms a basis of $H^*(\Gr)$ (or
$H_*(\Gr)$).

The transition matrix between $\{\WS_w(x) \mid w \in \tS^0\}$ and
$\{h_\lambda(x)\}$, given by Theorem~\ref{T:weakPieri} is the same
as the transition matrix between the monomial symmetric functions
$\{m_\lambda(x)\}$ and the strong Schur functions $\{\SSS_w(x) \mid
w \in \tS^0\}$. Note that $\{h_\lambda(x)\}$ does not form a basis
for $\Lambda^{(n)}$, so the matrix is ``rectangular'' (and
infinite). Since $\{\SSS_w(x) \mid w \in \tS^0\}$ is linearly
independent, the matrix has full rank (when restricted to
submatrices of each degree) and so the Pieri rule can be inverted to
write $\WS_w(x)$ in terms of $h_\lambda(x)$. Applying
Theorem~\ref{T:Lam}, we may write $\xi^w$ in terms of
$\xi^{c_{0,m}}$.  Similarly, the fact that $\{\WS_w(x) \mid w \in
\tS^0\}$ is linearly independent allows one to write $\{\SSS_w(x)
\mid w \in \tS^0\}$ in terms of $\{h_\lambda \mid \lambda_1 < n\}$.
\end{proof}

\section{Conjectured Pieri rule for the affine flag variety}
Theorems~\ref{T:weakPieri} and~\ref{T:cohomPieri} suggest that we
make the following conjecture.  In the following conjecture, we let
$l$ be arbitrary again.

\begin{conj}[Conjectured Pieri rule for $H^*(\G/\Bo)$]
\label{conj:flags} Let $w \in \tS$ and $1 \leq m$.  Then in
$H^*(\G/\Bo)$ we have
$$
\xi^{c_{l,m}}_B\, \xi^w_B = \sum_{S} \xi^{\out(S)}_B,
$$
where the sum runs over strong strips $S$ of size $m$ such that
$\ins(S) = w$.
\end{conj}

In~\cite{Lam} it is observed that $\WS_w(x) = \tF_w(x)$ is the
pullback of $\xi_B^w$ from $H^*(\G/\Bo)$ to $H^*(\Gr)$ under the
map $\Gr \to \G/\Bo$ induced by the map $\Omega SU(n) \to LSU(n)
\to LSU(n)/T$ where $\Omega SU(n)$ denotes the based loop-space,
$LSU(n)$ denotes the loop space, and $T$ denotes the maximal
torus. The pullback of Conjecture~\ref{conj:flags} is consistent
with Theorem~\ref{T:weakPieri}.

\begin{remark}
Recurrences for the structure constants of $H^*(\G/\Bo)$ are given
in~\cite{KK}.  It may be possible to derive
Conjecture~\ref{conj:flags} directly from these recurrences.
\end{remark}

\begin{remark}
Conjecture~\ref{conj:flags} is consistent with the Pieri rule for
the classical finite-dimensional flag manifold.  Indeed if $w,z
\in S_n \subset \tS$ and $l \in [1,n]$, then the existence of a
marked strong cover $\stc{w}{i,j}{z}$ is exactly the combinatorial
condition appearing in Monk's rule, while strong strips agree with
the ``path formulation'' of the Pieri rule in~\cite{Sot}.  Note
that in~\cite{Sot}, using the language that we have introduced, a
strong cover $\stc{w}{i,j}{z}$ would be ``marked'' at $w(i)$
rather than $w(j)$.  However the Pieri rules obtained from the two
different markings agree.
\end{remark}

\section{Geometric interpretation of strong Schur functions}
\label{sec:strongProperty} In this section we list some conjectural
properties of strong Schur functions, assuming $l = 0$ for
simplicity.

\begin{conj}
\label{conj:strongSym}  Let $u, v \in \tS$ be two affine
permutations. We have the following successively stronger
properties.
\begin{enumerate}
\item We have $\SSS_{u/v}(x) \in \Lambda$. \item We have
$\SSS_{u/v}(x) \in \Lambda_{(n)}$. \item We have $\SSS_{u/v}(x) =
\sum_{w \in \tS^0} c_{vw}^u \,    \SSS_{w}(x)$ where $c_{vw}^u$ is
defined by \eqref{E:affineflagconstants}.
\end{enumerate}
\end{conj}

The corresponding properties of weak Schur functions are known.
Symmetry (Theorem~\ref{T:weakSym}) was proven combinatorially
in~\cite{LamAS} (see also~\cite{LMdual}) while positivity was shown
in~\cite{Lam} using geometric work of Peterson~\cite{Pet}.

\begin{prop}
Conjecture~\ref{conj:strongSym} follows from
Conjecture~\ref{conj:flags}.
\end{prop}
\begin{proof}
By making the identification $\Lambda^{(n)} \cong H^*(\Gr)$, we
may consider the affine Cauchy kernel $\Omega_n$ as an element of
the completion $\Lambda_{(n)} \hat\otimes H^*(\Gr)$.  If
Conjecture~\ref{conj:flags} holds, then it completely determines
the action of $H^*(\Gr)$ on $H^*(\G/\Bo)$, obtained from the
inclusion $H^*(\Gr) \subset H^*(\G/\Bo)$.

Let $\inner{.}{.}_{\G/\Bo}$ denote the inner product on
$H^*(\G/\Bo)$ defined by $\inner{\xi^w_B}{\xi^v_B} = \delta_{wv}$.
By the definition of $\Omega_n$ and Conjecture~\ref{conj:flags},
we have
$$
\SSS_{u/v}(x) = \inner{\Omega_n \cdot \xi^v_B}{\xi^u_B}_{\G/\Bo}
$$
where $\Omega_n \cdot \xi^v \in \Lambda_{(n)} \hat\otimes
H^*(\G/\Bo)$.  By Corollary~\ref{C:cauchy} and
Theorem~\ref{T:Lam}, we may also write
$$
\Omega_n = \sum_{w \in \tS^0} \SSS_w(x) \otimes \xi^w
$$
so that
$$
\SSS_{u/v}(x) = \sum_{w \in \tS^0} \SSS_w(x) \inner{\xi_B^w\,
\xi^v_B}{\xi^u_B}_{\G/\Bo}.
$$
Thus $\SSS_{u/v}(x) \in \Lambda_{(n)}$.  By general positivity
results of~\cite{Gra,Kum}, we have $c_{vw}^u \in \Z_{\geq 0}$.
\end{proof}
%

\chapter{The Local Rule $\phi_{u,v}$} \label{ch:forward} In this
section we shall define a local rule as in Proposition
\ref{P:local}, as a sequence of operations called \textit{internal}
and \textit{external insertion} steps.

\section{Internal insertion at a marked strong cover}
\label{SS:intins}

Let $C=(\stc{w}{i,j}{u})$ be a marked strong cover. \textit{Internal
insertion at $C$} is a map that takes as input, a final pair of the
form $(W,S'_1)=(\wcf{w}{A}{v}{\cdot}{S'_1})$ and produces an output
final pair\footnote{This is an abuse of language: in the generality
in which we define internal insertion, $S'$ is only guaranteed to be
a strong tuple, not necessarily a strong strip. However, whenever we
apply internal insertion in the definition of the local rule
$\phi_{u,v}$, $S'$ will be a strong strip.} of the form
$(W',S')=(\wcf{u}{A'}{x}{\cdot}{S'})$ such that
$\ins(W')=u=\out(C)$, $\size(W')=\size(W)$, $\ins(S')=\ins(S'_1)$,
and $\size(S')=\size(S'_1)+1$. This given, we define
\begin{align*}
   C' = \lc(S') = (\stc{v'}{a,b}{x}).
\end{align*}
Internal insertion has three cases named A, B, and C. In Cases A and
B the output takes a particularly pleasant form: $S'$ is obtained by
appending $C'$ to $S'_1$ and $v'=v$. In Case C, $S'$ is obtained
from $S'_1$ by placing a strong cover just \textit{before} the last
cover of $S'_1$.


If $\size(S'_1)>0$ we write
\begin{align} \label{E:ydef}
S'_1 = (\xymatrix{%
{} \cdots {} \ar[r]^{\aspm,\bspm} & {y}
\ar[r]^{\asp,\bsp} & {v}}).
\end{align}
We need only specify $A'$, $(a,b)$, and the rule for obtaining $S'$
for then we set
\begin{align*}
  x &= c_{A'}(u) =\out(S') \\
  v' &= x\,t_{a,b}.
\end{align*}

\subsection{Commuting case} \

\medskip\noindent\textbf{Case A (Commuting case)} Suppose $(W,C)$
commutes. Set $A' = A$, $(a,b)=(i,j)$, and $S' = S'_1 \cup C'$.

\begin{example}
Let $n = 4, l = 0$ and $C = (\stc{[3,5,-2,4]}{-2,1}{[1,7,-2,4]})$.
 Consider internal insertion at $C$ of the final pair $$(W,S_1') =
 (\wcf{[3,5,-2,4]}{\{3\}}{[4,5,-2,3]}{[4,5,-2,3]}{}).$$  Since
 $(W,C)$ commute, the output
 final pair is
 $$
 (W',S') = (\wcf{[1,7,-2,4]}{\{3\}}{[1,8,-2,3]}{[4,5,-2,3]}{S'}).
 $$
 where $S' = (\stc{[4,5,-2,3]}{-2,1}{[1,8,-2,3]})$.

\end{example}

\subsection{Noncommuting cases} Otherwise,  we assume that $(W,C)$
does not commute.  Since $w(j)=u(i)$ is not $A$-nice by
Lemma~\ref{L:Fcyclicobstruction}, we can let
\begin{align}
\label{E:p0def} p_0 &= u(i)-1 \\
\label{E:Ahdef}
  \Ah &= A - \ba{p_0}.
\end{align}

\noindent \textbf{Case B (Normal bumping case)} Suppose that $(W,C)$
does not commute and either $\size(S'_1)=0$, or $\size(S'_1)>0$ and
$i\ne\asp$. Let $q<p$ be the unique pair of consecutive $\Ah$-nice
integers such that $q< u(j)$ and $u^{-1}(q)\le l$ and $q$ is
maximal. We set $A' = \Ah\cup\{\ba{p-1}\}$, $(a,b) =
(u^{-1}(q),u^{-1}(p))$, and $S'= S'_1 \cup C'$.

\begin{example}
Let $n = 6, l = 0$ and $$C =
(\stc{[5,0,1,9,-2,8]}{-2,1}{[3,0,1,11,-2,8]}).$$  Consider the
final pair $(W,S_1')$ given by $$W =
(\ws{[5,0,1,9,-2,8]}{\{3,4,5\}}{[4,-1,1,12,-3,8]})$$ and $S_1' =
(\stc{[4,-1,1,12,-3,8]}{}{[4,-1,1,12,-3,8]})$.  The pair $(W,C)$
does not commute and we have $p_0 = 4$ and $\Ah = \{3,5\}$.  Since
$\size(S_1') = 0$, we are in Case B.  One calculates that $q = 2,
p = 3$, so that $A' = \{2,3,5\}$ and $(a,b) = (0,1)$.  Thus the
output final pair is given by $W' =
(\ws{[3,0,1,11,-2,8]}{\{2,3,5\}}{[2,-1,1,12,-3,10]})$ and $S' =
(\stc{[4,-1,1,12,-3,8]}{0,1}{[2,-1,1,12,-3,10]})$.

\end{example}

\medskip \noindent \textbf{Case C (Replacement Bump)} Suppose that $(W,C)$
does not commute, that $\size(S'_1)>0$, and $i=\asp$. Let $q<p$ be
the unique pair of consecutive $\Ah$-bad integers such that $q <
p_0$ and $y^{-1}(q) \le l$ and $q$ is maximal. Set $A'=\Ah\cup
\{\ba{q}\}$, $(\am,\bm)=(y^{-1}(q),y^{-1}(p))$ and let $S'$ be
obtained by inserting $(\am,\bm)$ just \textit{before} the last pair
of indices $(\asp,\bsp)=(i,\bsp)$ of $S'_1$.

\begin{example}
Let $n = 4, l = 0$ and $C = (\stc{w = [1,7,-2,4]}{-2,4}{u =
[1,8,-2,3]})$. Consider internal insertion at $C$ of the final
pair
$$(W,S_1') =
(\wcf{[1,7,-2,4]}{\{3\}}{[1,8,-2,3]}{[4,5,-2,3]}{S_1'}),$$ where
$S'_1 = (\stc{y = [4,5,-2,3]}{-2,1}{v = [1,8,-2,3]})$. Since the
pair $(W,C)$ does not commute, $\size(S_1') = 1 > 0$ and $i = a_1
= -2$ we are in Case C.  We have $p_0 = 3$ and $\Ah = \emptyset$
so all integers are $\Ah$-bad.  We find that $q = 1, p = 2, A' =
\{1\}$ and obtain
$$
S' =
([4,5,-2,3]\overset{-2,7}{\longrightarrow}[4,6,-3,3]\overset{-2,1}{\longrightarrow}[2,8,-3,3]).
$$
Note that $W' = (\ws{[1,8,-2,3]}{\{1\}}{[2,8,-3,3]})$ is indeed a
weak strip.
\end{example}

\subsection{External Insertion}\

\medskip \noindent
\textbf{Case X (External Insertion)}
\label{sec:externalInsertion} Let $W=(\ws{w}{A}{v})$ be a weak
strip such that $|A|<n-1$. Let $q<p$ be the unique pair of
consecutive $A$-bad integers such that $v^{-1}(q) \le l$ and $q$
is maximal. We set $A' = A \cup \{\ba{q}\}$ and $(a,b) =
(v^{-1}(q),v^{-1}(p))$. By Lemma \ref{L:Faddresidue2} we have a
noncommuting final pair $(W',C')=(\wcf{w}{A'}{x}{v}{a,b})$.


\textit{External insertion} is the map $\phi_1$ that takes as input
a final pair $(W,S'_1)=(\wcf{w}{A}{v}{\cdot}{S'_1})$ with
$\size(W)<n-1$, and produces the final pair $(W',S')$, where
$W\mapsto (W',C')$ as above and $S'=S'_1\cup C'$.

\begin{example}
Let $n = 5, l = 0$ and $W =
(\ws{[2,-4,5,8,4]}{\{3,5\}}{[2,-5,6,9,3]})$.  We find $q = 4$ and
$p = 6$.  Then external insertion of $(W,S'_1)$, where $S'_1 =
(\stc{[2,-5,6,9,3]}{}{[2,-5,6,9,3]})$ produces the final pair
$(W',S')$ where
$$
W' = (\ws{[2,-4,5,8,4]}{\{3,4,5\}}{[2,-5,4,11,3]}),$$ and $S' =
(\stc{[2,-5,6,9,3]}{-1,3}{[2,-5,4,11,3]})$.
\end{example}

\section{Definition of $\phi_{u,v}$} Fix $u,v\in\tS$. We define the
value of $\phi_{u,v}$ on $(W,S,e)\in\IIL_{u,v}$ as the result of a
sequence of steps. Each step, which is either an internal or
external insertion, takes a final pair and produces another.

We start with the final pair $(W^{(0)},S^{(0)})=(W,\stc{v}{}{v})$
where $S^{(0)}$ is the empty strong strip from $v=\out(W)$ to
itself. Iteratively, for $1\le k \le m=\size(S)$, perform the
internal insertion on the final pair $(W^{(k-1)},S^{(k-1)})$ at
$C_k$, and let $(W^{(k)},S^{(k)})$ be the resulting final pair. The
result of this sequence of internal insertions is the final pair
$(W^{(m)},S^{(m)})$. We now perform $e$ external insertions. For
$m<k<m+e$ define
\begin{align*}
  (W^{(k)},S^{(k)}) = \phi_1(W^{(k-1)},S^{(k-1)}).
\end{align*}
We define $\phi_{u,v}(W,S,e)=(W',S'):=(W^{(m+e)},S^{(m+e)})$ to be
the final pair produced by this process.

\section{Proofs for the local rule} We now establish the
well-definedness of the local rule $\phi_{u,v}$ and some of its
properties.

\subsection{Case X} By construction we have a commutative diagram
\begin{align}
\xymatrix{%
{w} \ar@{=}[r] \ar[d]^W_A & {u} \ar@{-->}[d]^{A'}_{W'} \\
{v} \ar@{-->}[r]_{a,b}^{C'} & {x}}
\end{align}
where $(a,b)=(v^{-1}(q),v^{-1}(p))$. By Lemma \ref{L:extins} applied
to $S'_1$, $S'$ is a strong strip, finishing this case.

\begin{lem} \label{L:extins} In Case X let $W\mapsto (W',C')$.
Then for any strong strip $S$ such that $\out(S)=v$, $S'=S\cup C'$
is a strong strip.
\end{lem}
\begin{proof}
If $\size(S)=0$ then $S'$ is automatically a strong strip. Otherwise
let $\lc(S)=(\stc{\cdot}{i,j}{v})$.  By the maximality
of $q$, we have $v(i) \leq q$ since otherwise there would exist a pair of consecutive
$A$-bad integers $r$ and $s$  such that $r<v(i)<s$ with $v^{-1}(s) > l\geq i$, which would
contradict Lemma~\ref{L:weakcriterion} for the weak strip $(\ws{w}{A}{v})$.  Therefore
$m(C')=p>q\geq v(i)=m(\lc(S))$, so $S'$
is a strong strip.
\end{proof}

By Lemma \ref{L:extins}, each external insertion sends a final pair
to a final pair, preserves the inside permutations of both the weak
and strong strip, and adds one to the sizes of the weak and strong
strips. Thus to check that $\phi_{u,v}$ is well-defined, we may
reduce to the case $e=0$ where no external insertions are required.

\subsection{For internal insertion cases} We want to compute
$\phi_{u,v}(W,S,0)$ with $\size(S)=m$. By induction we may assume
that all of the internal insertions have been performed except the
last, which computes the internal insertion on
$(W^{(m-1)},S^{(m-1)})$ at $C_m$, resulting in $(W^{(m)},S^{(m)})$.
To avoid the proliferation of subscripts and superscripts we change
notation, forgetting the global meaning of $u,v,W,S$. We denote this
last internal insertion step as the internal insertion on $(W,S'_1)$
at $C$, resulting in $(W',S')$. We write
\begin{equation} \label{E:Fnotation}
\begin{split}
&(W,C)=(\wci{w}{A}{v}{i,j}{u}) \\
&\xymatrix{%
{S:} & {\dotsm} \ar[r]^{i^{--},j^{--}}_{C^{--}} & {w^-}
\ar[r]^{\im,\jm}_{\Cm} & {w} \ar[r]^{i,j}_C & {u} } \\
&\xymatrix{{S'_1:}&{\dotsm} \ar[r]^{\aspm,\bspm} & {y}
\ar[r]^{\asp,\bsp} & {v} } \\
&\xymatrix{%
{S':} & {\dotsm} \ar[r]^{a^{--},b^{--}}_{{C'}^{--}} & {y}
\ar[r]^{\am,\bm}_{\Cpm} & {v'} \ar[r]^{a,b}_{C'} & {x} } \\
&(W',C') = (\wcf{u}{A'}{x}{v'}{a,b})
\end{split}
\end{equation}

We use the following induction hypothesis.

\begin{property} \label{pro:markmove} \
\begin{enumerate}
\item[(i)] (a) $m(C') < m(C)$ in Case B and (b) $m(C') = c_{A'}(m(C)))$
in Cases A and C.
\item[(ii)] $x(b)<m(C)$.
\item[(iii)] Case C cannot be preceded by Case B, and if Case C
holds then $\im=i$.
\item[(iv)] The final pair $(W',C')$ commutes in Cases A and C and
does not commute in Case B.
\end{enumerate}

\end{property}

\subsection{Case A}
\begin{align} \label{E:FAdiag}
\xymatrix{%
{\dotsm}\ar[r]^{\im,\jm} & {w} \ar[r]^{i,j} \ar[d]_A & {u} \ar@{-->}[d]^A \\
{\dotsm}\ar[r]_{\asp,\bsp} & {v} \ar@{-->}[r]_{i,j} & {x} }
\end{align}
By Lemma \ref{L:initfincommute}, $(\wcf{u}{A}{x}{v}{i,j})$ is a
commuting final pair.

\begin{lem} In Case A, Property \ref{pro:markmove} is satisfied.
\end{lem}
\begin{proof} We have $m(C')=x(i) = c_A u(i)=c_A(m(C))$, which proves (i).
(ii) is equivalent to $c_A w(i)<w(j)$, but this follows by Lemma
\ref{L:Fcyclicobstruction}. (iv) was proved above.
\end{proof}

\begin{lem}
\label{L:Astrong} In Case A, $S'$ is a strong strip.
\end{lem}
\begin{proof} We use the notation in \eqref{E:FAdiag} where the top
row gives $S$ and the bottom row gives $S'$. By induction $S'_1$ is
a strong strip, so we need only check that $m(\Cpm)<m(C')$, that is,
$v(\asp) < v(j)$.

Since $S$ is a strong strip $w(\im)<w(j)$. By induction,
Property~\ref{pro:markmove}(i) asserts that either (a)
$v(\asp)<w(\im)$ or (b) $v(\asp)=c_A(w(\im))$. Suppose (a) holds. By
\eqref{E:cvalue}, $v(j)=c_A(w(j))) \ge w(j)-1 \ge w(\im)>v(\asp)$,
as desired.

Suppose (b) holds. Then $\im=\asp$, since $v(\im)=c_A w(\im) =
v(\asp)$. In particular $w(\asp)<w(j)$.

We shall assume that $v(\asp) \ge v(j)$ and derive a contradiction.
If $v(\asp)=v(j)$ then we have the contradiction $l \ge \asp=j> l$.
So we assume $v(\asp)>v(j)$. In other words $w(\asp)<w(j)$ is
inverted by $c_A$. By Lemma \ref{L:cycdecrinv}, $w(\asp)$ is
$A$-nice, $c_A w(\asp) \ge w(j)$, and $c_A w(j)=w(j)-1$. By
induction, Property~\ref{pro:markmove}(ii) gives
$v(\bsp)<w(\im)=w(\asp)\le w(j)-1 = v(j)$. So
\begin{equation} \label{E:trap}
v(\bsp) < v(j) < v(\asp).
\end{equation}
Next we have
\begin{equation}
\label{eq:casec1} \asp \le l <  \bsp < j.
\end{equation}
This follows from $l<j$ and $j\not\in (\asp,\bsp)$, which follows
from Lemma \ref{L:strongcover} with the strong cover $\Cpm =
(\stc{\cdot}{\asp,\bsp}{v})$ and \eqref{E:trap}. We have
\begin{equation} \label{E:neqsx}
  v(\bsp) < v(i) < v(j) < v(\asp).
\end{equation}
This is obtained from \eqref{E:trap} using Lemma \ref{L:strongcover}
for the strong cover $C'=(\stc{v}{i,j}{\cdot})$ and $i\le l < \bsp <
j$. We have
\begin{equation}\label{E:neqi}
i<\asp\le l<\bsp<j
\end{equation}
since $i\not\in(\asp,\bsp)$ by Lemma \ref{L:strongcover} for the
strong cover $\Cpm$. We have
\begin{equation} \label{E:neqsu}
w(\asp)<w(i)<w(j).
\end{equation}
This follows from Lemma \ref{L:strongcover} for the cover $C$,
$i<\asp\le l < j$, and $w(\asp)<w(j)$.

Now $c_A$ inverts $w(\asp)<w(j)$. By Lemma \ref{L:cycdecrinv} and
\eqref{E:neqsu}, $c_A$ inverts $w(\asp)<w(i)$ as well. With
\eqref{E:neqi}, Lemma \ref{L:weaknoinv} gives a contradiction.
\end{proof}

\subsection{Case B} We first sketch the proof of Case B and fill in
the proofs afterwards. Since $(W,C)$ does not commute, by Lemma
\ref{L:Fremoveresidue}, $\ws{u}{\Ah}{v}$ is a weak strip such that
the diagram commutes:
\begin{align} \label{E:FBstartdiag}
\xymatrix{%
{w} \ar[r]^{i,j} \ar[d]_A & {u} \ar@{-->}[dl]^{\Ah}  \\ %
{v}  & {}
}%
\end{align}
Recall that in Case B, $(a,b)= (u^{-1}(q),u^{-1}(p))$, $A' =
\Ah\cup\{\ba{p-1}\}$, and $x = c_{A'} u$. Using Lemma
\ref{L:casebproperties} we may apply Lemma \ref{L:Faddresidue},
which says that $(\wcf{u}{A'}{x}{v}{a,b})$ is a noncommuting final
pair such that the diagram commutes.
\begin{align*}
\xymatrix{%
{w} \ar[r]^{i,j} \ar[d]_A & {u} \ar@{-->}[dl]^{\Ah} \ar@{-->}[d]^{A'} \\ %
{v} \ar@{-->}[r]_{a,b} & {x}
}%
\end{align*}
Lemma \ref{L:FBstrong} shows that $S'$ is a strong strip, completing
the proof sketch for Case B.

\begin{lem} \label{L:FpreB} In Case B, let
\begin{align} \label{E:q'def}
  q' &= u(i)-n, \\
\label{E:p'def}
  p' &= u(j).
\end{align}
Then
\begin{align}
\label{E:Bprstraddle}
  u^{-1}(q') \le l < u^{-1}(p'), \\
\label{E:pqnprime}
   0 < p'-q' < n,
\end{align}
$q'$ is $\Ah$-nice, $p'$ is $\Ah$-nice and $A$-nice, and
$|\Ah|=|A|-1$.
\end{lem}
\begin{proof} $q'$ is $\Ah$-nice by definition.
$p'$ is $A$-nice and $\Ah$-nice by Lemma \ref{L:Fcyclicobstruction}.

We have $p'-q'=u(j)-u(i)+n=w(i)-w(j)+n$. The second inequality in
\eqref{E:pqnprime} holds since $C$ is a strong cover. Lemma
\ref{L:Fcyclicobstruction} says that $w(i)<w(j)$ are inverted by
$c_A$, and in particular, $w(j)=u(i)$ is not $A$-nice and
$w(j)-w(i)<n$, which implies that $|\Ah|=|A|-1$ and the first
inequality in \eqref{E:pqnprime}.

By straddling, $u^{-1}(q')=i-n<i\le l<j =u^{-1}(p')$, so
\eqref{E:Bprstraddle} holds.

\end{proof}

\begin{lem}
\label{L:casebproperties} In Case B, $q$ is well-defined and
\begin{align}
\label{eq:Bbounds} q' = u(i)-n &\le q < p \le u(j) = p' \\
\label{eq:Bstraddle} u^{-1}(q) &\le l < u^{-1}(p)
\end{align}
Moreover $|A'|=|A|$, $q$ is $A'$-nice and $p$ is $A$-nice.
\end{lem}
\begin{proof} $q$ is well-defined, is $\Ah$-nice, and
satisfies the above inequalities by Lemma \ref{L:FpreB}, which
assures that $q$ is the maximum of a set that contains $q'$. $p$ is
$A$-nice and satisfies the above inequalities, also thanks to Lemma
\ref{L:FpreB}. $q$ is $A'$-nice since it is $\Ah$-nice and too close
to $p$ and distinct from $p$ to be congruent to $p$, by
\eqref{eq:Bbounds} and \eqref{E:pqnprime}. Since $p$ is $\Ah$-nice
it follows that $|A'|=|\Ah|+1=|A|$.
\end{proof}

\begin{lem}
\label{L:CaseBProp} In Case B Property \ref{pro:markmove} holds.
\end{lem}
\begin{proof} Property \ref{pro:markmove}(iv) was already proved
above. Note that if $m(C')<m(C)$ then Property
\ref{pro:markmove}(ii) follows immediately: we have $x(b) < x(a) =
m(C')<m(C)$ since $C'$ is a strong cover.

To prove Property \ref{pro:markmove}(i), we have
\begin{equation} \label{E:casebmr}
m(C') = x(u^{-1}(q)) = c_{A'}(q) = c_{\Ah}(p)=\begin{cases}
u(i)-1 & \text{if $p=u(j)$} \\
c_A(p) &\text{otherwise.}
\end{cases}
\end{equation}
The first two equalities hold by definition. The third follows from
\eqref{E:cvalue2}.  For the last equality, $u(j)$ is $A$-nice by
Lemma \ref{L:Fcyclicobstruction}. If $p=u(j)$, then by Lemma
\ref{L:Fcyclicobstruction}, $u(i)$ is the next $\Ah$-nice integer
after $u(j)$, and $c_{\Ah}(p)=u(i)-1$ by \eqref{E:cvalue2}.
Otherwise by \eqref{eq:Bbounds} $p<u(j)$ and the cyclic component of
$\ba{p}$ is the same in $A$ and $\Ah$. Either way the last equality
holds.

For $p=u(j)$, \eqref{E:casebmr} implies that $m(C')<u(i)=m(C)$. For
$p<u(j)$, since $u(j)$ is $A$-nice, by \eqref{E:casebmr} and the
fact that $C$ is a strong cover,  \eqref{E:cvalue2} gives $m(C')=c_A(p)<u(j)<u(i)=m(C)$ as
desired.
\end{proof}

\begin{lem} \label{L:FBstrong}
In Case B, $S'$ is a strong strip.
\end{lem}
\begin{proof} We use the notation in the diagram below, where the
top row is $S$ and the bottom row is $S'$.
\begin{align*}
\xymatrix{%
{\dotsm} \ar[r] & {\cdot} \ar[r]^{\im,\jm}_{\Cm} & {w} \ar[r]^{i,j}_C \ar[d]_A & {u}  \ar[d]^{A'} \\ %
{\dotsm} \ar[r] & {y} \ar[r]_{\asp,\bsp}^{\Cpm} & {v}
\ar[r]_{a,b}^{C'} & {x}
}%
\end{align*}
By induction it suffices to show that  $m(\Cpm)<m(C')$. Since
$v=c_{\Ah}u$ (see \eqref{E:FBstartdiag}) this is equivalent to
\begin{align} \label{E:FBdownstr}
 c_{\Ah}(u(\asp)) < c_{\Ah}(p).
\end{align}
Since $p$ is $\Ah$-nice it suffices to show that
\begin{align} \label{E:FBineq}
  u(\asp) < p.
\end{align}
Since $S$ is a strong strip we have
\begin{align} \label{E:FBupstr}
w(\im)<w(j).
\end{align}
The cases of Property~\ref{pro:markmove}(i) lead in (a)
to $w(i^-)> c_A w(a_1) \geq w(a_1)-1$ from  \eqref{E:cvalue2}, and in (b)
to $c_A w(a_1)=c_A w(i^-)$.  We thus find that
\begin{align} \label{E:FBwa1im}
  w(\asp) \le w(\im).
\end{align}
Suppose that $\ba{\asp}\not\in \{\ba{i},\ba{j}\}$. Then $u(\asp)\le
w(\im)$. We claim that
\begin{equation} \label{E:casebless}
u(\asp) < u(j).
\end{equation}
If $\im=i$ then since Case B has $i\ne \asp$ this implies
$u(\asp)<w(i)=u(j)$ and \eqref{E:casebless} holds. So we assume
$\im\ne i$. Suppose $w(\im)>w(i)$. Since $C$ is a strong cover and
$\im<l$, by Lemma \ref{L:strongcover} and \eqref{E:FBupstr} we have
$\im < i \le l < j$. By Lemma \ref{L:Fcyclicobstruction}(iv) and  \eqref{E:cvalue2},
$c_A$ inverts $w(i)<w(\im)$ since $w(i)<w(\im)<w(j)$, thus contradicting
Lemma \ref{L:weaknoinv}.
Therefore $w(\im)<w(i)$, which, together with \eqref{E:FBwa1im}
yields \eqref{E:casebless}.

Let $q_1$ be the maximum $\Ah$-nice integer such that $q_1\le
u(\asp)$. By Lemma \ref{L:weakcriterion}, $u^{-1}(q_1)\le \asp\le
l$. By the definition of $q$ and~\eqref{E:casebless}, we have $q\ge
q_1$. But $q<p$ and $p$ is $\Ah$-nice so \eqref{E:FBineq} holds.

Next suppose $\ba{\asp}=\ba{j}$. Let $\asp=j+kn$ for some $k\in\Z$.
Since $\asp\le l<j$ we have $k<0$ and we have $u(\asp)=u(j)+kn\le
u(j)-n<u(i)-n<p$ by \eqref{eq:Bbounds}, proving \eqref{E:FBineq}.

Finally suppose $\ba{\asp}=\ba{i}$ but $\asp\ne i$. Let $\asp=i+kn$
for some $k\ne0$. By \eqref{E:FBupstr} and \eqref{E:FBwa1im} we have
$w(\asp)<w(j)$, that is, $w(i)+kn<w(j)$. Now $0 < w(j)-w(i)<n$ by
Lemma \ref{L:Fcyclicobstruction}, so $k<0$. We have
\begin{align*}
  u(\asp) = u(i)+kn \le u(i)-n < p
\end{align*}
by \eqref{eq:Bbounds}, proving \eqref{E:FBineq}.
\end{proof}

\subsection{Case C} We first sketch the overall strategy for Case C
using diagrams and fill in the proofs afterwards. By Lemma
\ref{L:Fremoveresidue}, $\ws{u}{\Ah}{v}$ is a weak strip such that
the diagram commutes.
\begin{align} \label{E:FCstartdiag}
\xymatrix{%
{} & {w} \ar[r]^{i,j} \ar[d]_A & {u} \ar@{-->}[dl]^{\Ah} \\ %
{y} \ar[r]_{i,\bsp} & {v} &
}%
\end{align}
By Lemma \ref{L:FCcom1}, $(\wcf{u}{\Ah}{v}{y}{i,\bsp})$ is a
commuting final pair. Let $w' = u\, t_{i,\bsp}=c_{\Ah}^{-1}y$. By
Lemma \ref{L:initfincommute}, $(\wci{w'}{\Ah}{y}{i,\bsp}{u})$ is a
commuting initial pair.
\begin{align} \label{E:FCw'diag}
\xymatrix{%
{} & {w'} \ar@{-->}[r]^{i,\bsp} \ar@{-->}[dl]_{\Ah}  & {u}
\ar[dl]^{\Ah}
 \\ %
{y} \ar[r]_{i,\bsp} & {v} & {}
}%
\end{align}
Let $(\am,\bm)=(y^{-1}(q),y^{-1}(p))$ and $v'= t_{qp}\,y =
y\,t_{\am,\bm}$. Using Lemma \ref{L:casecproperties} we may apply
Lemma \ref{L:Faddresidue2}, which shows that
$(\wcf{w'}{A'}{v'}{y}{(\am,\bm)})$ is a noncommuting final pair such
that the diagram commutes.
\begin{align*}
\xymatrix{%
{} & {w'} \ar[r]^{i,\bsp} \ar[dl]_{\Ah} \ar@{-->}[d]_{A'} & {u}
 \\ %
{y} \ar@{-->}[r]_{\am,\bm} & {v'}  & {}
}%
\end{align*}
Let $x = c_{A'} u = v'\,t_{i,\bsp}$. By Lemma \ref{L:FCcom2} the
initial pair $(\wci{w'}{A'}{v'}{i,\bsp}{u})$ commutes, so that by
Lemma \ref{L:initfincommute} we have a commuting final pair
$(\wcf{u}{A'}{x}{v'}{i,\bsp})$.

\begin{align} \label{E:FCenddiag}
\xymatrix{%
{} & {w'} \ar[r]^{i,\bsp} \ar[dl]_{\Ah} \ar[d]_{A'} & {u}
\ar@{-->}[d]^{A'}
 \\ %
{y} \ar[r]_{\am,\bm} & {v'} \ar@{-->}[r]_{i,\bsp} & {x}.
}%
\end{align}
This shows that the definition of $S'$ produces a strong tuple.
Lemma \ref{L:FCstrong} shows that $S'$ is a strong strip, completing
the argument for Case C.

\begin{lem} \label{L:cpropearly} In Case C, Property~\ref{pro:markmove}(iii) holds.
\end{lem}
\begin{proof}
Case B could not have occurred in the previous step for otherwise
we would have $c_A w(i) = v(i) = m(\Cpm)< m(\Cm) < m(C) = w(j)$,
contradicting Lemma~\ref{L:Fcyclicobstruction}.

So the previous step must be Case A or Case C. The input reflection
for that step is $(\im,\jm)$. In Case A the output reflection is
also $(\im,\jm)$, and in Case C it has the form $(\im,\cdot)$. But
the output reflection for this step is $(i,\bsp)$ so $\im=i$.
\end{proof}

\begin{lem} \label{L:p0} In Case C,
\begin{align*}
 (i) & {\rm ~If~} \bar j \neq \bar b_1 {\rm ~then~} p_0 =y(j).  \\
(ii) & {\rm ~Otherwise~} j=b_1+kn {\rm ~with~} k>0,  p_0=y(i)+kn  {\rm ~and~} b_1-i <n.
\end{align*}
Furthermore, letting
\begin{align} \label{E:q0def}
q_0 = y(\bsp)-n \, ,
\end{align}
we have
\begin{align}
\label{E:pq0} 0 &< p_0 - q_0 < n \\
\label{E:p0pos} l &< y^{-1}(p_0) \, .
\end{align}
\end{lem}
\begin{proof}

Since $\lc(S'_1)$ is a marked strong cover, we have
\begin{equation} \label{E:ystrong}
y(i)<y(\bsp) \qquad \text{and}\qquad i \le l < \bsp.
\end{equation}
To prove \eqref{E:pq0}, we have
\begin{equation} \label{E:yz=cujr}
y(\bsp)=y\,t_{i,\bsp}(i)=v(i)=c_A w(i)=c_A u(j)
\end{equation}
so that
\begin{align} \label{E:q0val}
  q_0 = c_A u(j) - n.
\end{align}
Then $p_0-q_0=(u(i)-1)-(c_Au(j)-n)=n-1-(c_A u(j)-u(i))$ and
\eqref{E:pq0} is equivalent to $0 \le c_A u(j)-u(i)<n-1$. But this
follows from
\begin{equation} \label{E:noncomm}
c_A u(j)-n < u(j) < u(i) \le c_A u(j),
\end{equation}
which holds by \eqref{E:cvalue2}, the fact that $C$ is a marked
strong cover, and Lemma \ref{L:Fcyclicobstruction}. So \eqref{E:pq0}
holds.

We have
\begin{equation} \label{E:p0comp}
\begin{split}
p_0&=u(i)-1=c_A(u(i))=c_A(w(j))=v(j) \\
&= y \,t_{i,\bsp}(j)=
\begin{cases} %
y(i)+kn&\text{if $j=\bsp+kn$ for some $k\in\Z$.} \\ %
y(j) &\text{otherwise.} %
\end{cases}
\end{split}
\end{equation}
This gives that in case $(i)$ $p_0=y(j)$, which yields
\eqref{E:p0pos}.

Otherwise, we are in  case $(ii)$,  with $j=\bsp+kn$ for
some $k\in\Z$ and $p_0=y(i)+kn$.  We thus need to show that $k>0$ and $b_1-i<n$.
By \eqref{E:noncomm} and \eqref{E:yz=cujr} we have
$y(\bsp)-n<u(i)-1=p_0$. Substituting for $p_0$ using
\eqref{E:p0comp} we have $y(\bsp)-n<y(i)+kn$, which gives that
 $k\geq0$ since $y(i)<y(b_1)$.  Now, suppose $k=0$ so that
$j=\bsp$ and $p_0=y(i)$. By Property
\ref{pro:markmove}(iii) and Lemma~\ref{L:cpropearly}, Case A must
have occurred in some step and Case C must have occurred in every
step thereafter. Let $(C_1,C_2,\dotsc,C_r = C)$ be the corresponding
sequence of marked strong covers in $S$ (starting with the above
instance of Case A) with $C_k=(\stc{u_{k-1}}{i,j_k}{u_k})$, where
the first index of the reflection is always $i$. The Case A step had
the same output reflection as its input reflection $(i,j_1)$. In the
output strong strip, the applications of Case C inserted reflections
just before this reflection, which stayed at the end, so that
$(i,j_1)=(i,j)$. Now, since $C_1$ is a strong cover we have $u_1(i)
> u_1(j)$.  And since each $C_k$ for $k > 1$ is a strong cover we
have $u_k(i) > u_{k-1}(i) > u_{k-1}(j) \geq u_k(j)$, using $\bar j
\neq \bar i$. This contradicts $C_r=(\stc{u_{r-1}}{i,j}{u_r})$ being
a strong cover and proves that $k>0$ in case $(ii)$.

Since $k>0$, we have from Lemma~\ref{L:Rcyclicobstruction}
$$
v(b_1) \leq v(j)-n < v(j) < v(i) < v(j)+n \, .
$$
We thus have $|v(i)-v(b_1)|>n$, which yields $b_1 < i+n$ from Lemma~\ref{L:strongcover}.  Therefore, $(ii)$ holds.

Finally, given that $p_0=y(i)+kn$ with $k>0$ implies $y^{-1}(p_0)=i+kn > b_1 > l$, we have that
\eqref{E:p0pos} follows.
\end{proof}

\begin{lem} \label{L:FCcom1} In Case C,
\begin{align} \label{E:FCmidcomm}
  u(i)>u(\bsp).
\end{align}
\end{lem}
\begin{proof}
By Lemma~\ref{L:cpropearly} the previous step had to be Case A or
Case C, and the final pair $(\wcf{w}{A}{v}{y}{i,\bsp})$ it produced
was commuting, by Property \ref{pro:markmove}(iv). By Lemma
\ref{L:Rcyclicobstruction}, $w(i)>w(\bsp)$. By Lemma~\ref{L:p0}, we have either
$(i)$ $\bar b_1 \neq \bar j$ or $(ii)$ $j=b_1+kn$,
with both implying $u(j)>u(\bsp)$. Since $C$ is a strong cover,
$u(i)>u(j)$, so that $u(i)>u(\bsp)$ as desired.
\end{proof}

\begin{lem} \label{L:casecproperties} In Case C, $q$ and $p$ are well-defined
and satisfy
\begin{align}
\label{E:pqbound} q_0 &\le q < p \le p_0 \\
\label{E:pqstraddle} y^{-1}(q) &\le l < y^{-1}(p).
\end{align}
Moreover $q$ is $A$-bad and $p$ is $A'$-bad.
\end{lem}
\begin{proof}
Given the other assertions, it is easy to see that $q$ is $A$-bad from  \eqref{E:pq0}, and
that $p$ is $A'$-bad.

\begin{sublemma} \label{sl:q}
There is an $\Ah$-bad integer $q$ such that
\begin{align}
\label{E:qbound} q_0 &\le q < p_0 \\
\label{E:qpos} y^{-1}(q) &\le l.
\end{align}
\end{sublemma}

This suffices: if $q$ is maximal and $p$ is the next $\Ah$-bad integer
after $q$, then $p$ satisfies \eqref{E:pqstraddle} by the choice of
$q$ and \eqref{E:pqbound} holds by \eqref{E:p0pos}.%
\end{proof}

\begin{proof}[Proof of Sublemma \ref{sl:q}]
Note that
\begin{equation} \label{E:qp0}
y(i) = v(\bsp) < m(\Cm) = w(\im) = w(i)
\end{equation}
by Property \ref{pro:markmove}(ii) applied to the previous step, and
Lemma~\ref{L:cpropearly}.

By Lemma \ref{L:strongcover} applied to the strong cover $y \lessdot
y \,t_{i,\bsp}$ we have either (a) $\bsp-i<n$ or (b)
$y(\bsp)-y(i)<n$. Suppose (a) holds. We shall show that $q=q_0$
satisfies Sublemma \ref{sl:q}. We have $q_0=v(i)-n= c_A w(i)-n$. By
Lemma \ref{L:Fcyclicobstruction}, $w(i)$ is $A$-nice. By Lemma
\ref{L:nicebad}, $q_0$ is $A$-bad and $\Ah$-bad. Clearly $q=q_0$
satisfies \eqref{E:qbound}. We have $y^{-1}(q_0)=\bsp-n<i\le l$
using (a) and \eqref{E:ystrong}, so that $q=q_0$ satisfies
\eqref{E:qpos}.

Now suppose (b) holds. Let $q'$ be the minimum $\Ah$-bad integer
such that $q' \ge y(i)$. We shall show that $q=q'$ satisfies
Sublemma \ref{sl:q}. We have
\begin{equation} \label{E:q0q}
q_0 = y(\bsp)-n<y(i) \le q' \le w(i)-1 <  p_0<y(\bsp)
\end{equation}
by assumption (b), the definition of $q'$, the fact that $w(i) - 1$
is $\Ah$-bad and greater or equal to $y(i)$ by \eqref{E:qp0}, and
\eqref{E:pq0}. It suffices to show
\begin{equation} \label{E:q'pos}
y^{-1}(q')\le l.
\end{equation}
If $y(i)$ is $\Ah$-bad, then $q'=y(i)$ and \eqref{E:q'pos} follows.
So we may assume that $y(i)$ is not $\Ah$-bad and $y(i)<q'$. Since
$\ws{w'}{\Ah}{y}$ is a weak strip (see \eqref{E:FCw'diag}), applying
Lemma \ref{L:weakcriterion} to $y$, the $\Ah$-bad integer $q'$ and
the next smaller one, we obtain $y^{-1}(q') < i \le l$ and
\eqref{E:q'pos} follows.
\end{proof}

\begin{lem} \label{L:FCcom2} In Case C, $v'(i)<v'(\bsp)$.
\end{lem}
\begin{proof} It is equivalent to show that
\begin{align} \label{E:FCcommneed}
  t_{q,p} \, y(i) < t_{q,p} \,y(\bsp).
\end{align}
We have
\begin{align} \label{E:ysand}
  y(\bsp)-n = q_0 \le q < p \le p_0 < y(\bsp)
\end{align}
by \eqref{E:pqbound} and \eqref{E:pq0}. Using this one may deduce
that
\begin{align} \label{E:tqpyb1}
  y(\bsp) \le  t_{q,p}(y(\bsp)).
\end{align}
If $\ba{y(i)}\not\in\{\ba{q},\ba{p}\}$, then $t_{q,p}\,y(i) =
y(i)<y(\bsp)$ using \eqref{E:ystrong}. If
$\ba{y(i)}\in\{\ba{q},\ba{p}\}$, then $t_{q,p}\,y(i) \le p<y(\bsp)$
by \eqref{E:ysand}. Either way \eqref{E:FCcommneed} holds by
\eqref{E:tqpyb1}.
\end{proof}

\begin{lem} \label{L:cprop} In Case C, Property~\ref{pro:markmove} holds.
\end{lem}
\begin{proof}
Property~\ref{pro:markmove}(iii) is Lemma~\ref{L:cpropearly}.

We have $m(C') = x(i) = c_{A'}(u(i)) = c_{A'}(m(C))$, proving
Property~\ref{pro:markmove}(i). For Property~\ref{pro:markmove}(ii),
we must show that $t_{q,p}\,y(i) \leq w(j) - 1 = p_0$. By the proof
of Lemma \ref{L:FCcom2}, either $t_{q,p}y(i) = y(i)$ or
$t_{q,p}y(i)\le p$.  Using either \eqref{E:qp0} or \eqref{E:pqbound}
we have $t_{q,p}y(i) \le p_0$ as desired.


\end{proof}

\begin{lem} \label{L:FCstrong}
In Case C, $S'$ is a strong strip.
\end{lem}
\begin{proof} We use the following notation:
\begin{align*}
&\xymatrix{%
{S:\,\,}&{\dotsm} \ar[r]_{C^{--}} & {\cdot} \ar[r]^{i,\jm}_{\Cm} &
{w} \ar[r]^{i,j}_C & {u}
} \\ %
&\xymatrix{%
{S'_1:}&{\dotsm} \ar[r]^{\aspm,\bspm} & {y} \ar[r]^{i,\bsp} & {v}
} \\ %
&\xymatrix{%
{S':\,}&{\dotsm} \ar[r]^{\aspm,\bspm}_{{C'}^{--}} & {y}
\ar[r]^{\am,\bm}_{\Cpm} & {v'} \ar[r]^{i,\bsp}_{C'} & {x}}
\end{align*}
where $(\am,\bm)=(y^{-1}(q),y^{-1}(p))$. We have $\im=i$ by Lemma
\ref{L:cprop}.

By induction and the definition of $S'$, we need only check its last
two pairs of consecutive marks. We have $m(C')=x(i)= c_{A'} \, u(i)
> u(i)-1= p_0 \ge p = m(\Cpm)$, using the fact that $u(i)$ is
$A'$-nice (since $\bar p_0 \neq \bar q$ from \eqref{E:pq0} and Lemma~\ref{L:casecproperties})
and \eqref{E:pqbound}. It remains to verify that $m(\Cpm)
> m({C'}^{--})$, that is,
\begin{align} \label{E:FCpya1}
p > y(\aspm).
\end{align}
Suppose that $y(\aspm) \le p_0$. Since $y^{-1}(p_0)>l$ and $\aspm\le
l$ we have $y(\aspm)<p_0$. Let $r<r'$ be consecutive $\Ah$-bad
integers such that $r \le y(\aspm) < r'$. Since $p_0$ is $\Ah$-bad
we have $r'\le p_0$.  If $r=y(a_1^-)$, then $y^{-1}(r)=a_1^- \leq l$, and thus
$p>q\geq r=y(a_1^-)$ by the maximality of $q$.  Otherwise, by Lemma \ref{L:weakcriterion} for the weak
strip $\ws{w'}{\Ah}{y}$ we have $y^{-1}(r') < \aspm \le l$. As
before $r'<p_0$. Therefore by the definition of $p$ and $q$, $p>q\ge
r'>y(\aspm)$ as desired.

Otherwise $p_0 < y(\aspm)$. Since $S'_1$ is a strong strip, its last
two marks satisfy $y(\aspm) < y(\bsp)$. By \eqref{E:qp0},
\begin{align} \label{E:yijab}
  y(i)<w(i)\le p_0 <y(\aspm)<y(\bsp).
\end{align}
By Lemma \ref{L:strongcover} applied to the strong cover $y\lessdot
y\,t_{i,\bsp}$, we have
\begin{equation}\label{E:ar2}
\aspm < i
\end{equation}
since $\aspm \leq l < \bsp$. We claim that
\begin{align} \label{E:vja1mi}
  v(j) < v(\aspm) < v(i).
\end{align}
According to Lemma~\ref{L:p0}, we have two cases to
consider.  First, suppose that Lemma~\ref{L:p0}$(i)$ holds.
We have $y(\bsp)=vt_{i,\bsp}(\bsp)=v(i)$.
Since $t_{ij}$ is a
reflection and $\bar b_1 \neq \bar j$, $y(j)=v(j)$. Therefore
\eqref{E:yijab} gives $v(j)<y(\aspm)<v(i)$. We claim that
$\ba{\aspm}\not\in\{\ba{i},\ba{\bsp}\}$. Suppose
$\ba{\aspm}=\ba{i}$. By \eqref{E:yijab} $y(i)<y(\aspm)$, which
implies that $i<\aspm$, contradicting \eqref{E:ar2}. Suppose
$\ba{\aspm}=\ba{\bsp}$. Write $\aspm=\bsp+kn$. By \eqref{E:yijab}
$y(\aspm)<y(\bsp)$ so $k<0$. Then since $\bar b_1 \neq \bar j$, we have
\begin{align*}
  v(j)=y(j)<y(\aspm)=y(\bsp)+kn \le y(\bsp)-n=v(i)-n.
\end{align*}
But $\wci{w}{A}{v}{i,j}{u}$ is a noncommuting pair so by Lemma
\ref{L:Fcyclicobstruction}
$w(j)-w(i)<n$, which leads to the contradiction $v(i)-v(j)<n$.
Therefore our claim holds. Consequently $y(\aspm)=v(\aspm)$ and
\eqref{E:vja1mi} holds in this case.

Now suppose that Lemma~\ref{L:p0}$(ii)$ holds.  We have $y(b_1)=v(i)$ and $p_0=v(j)$ from \eqref{E:p0comp}.
Since we still have $v(i)-v(j)<n$, we obtain $y(b_1)-p_0<n$, and so $\ba{\aspm}\not\in\{\ba{i},\ba{\bsp}\}$
holds.  This leads again to \eqref{E:vja1mi}.

Finally, by Lemma
\ref{L:Fcyclicobstruction}, we know that there are no $A$-nice integer
in the interval $(w(i),w(j)]$ and thus by Lemma~\ref{L:nicebad}, that there are no
$A$-bad integer in the interval $[v(j),v(i))$.  Letting $r$ be the
$A$-bad integer immediately before $v(i)$, we have from \eqref{E:vja1mi} that
$r<v(a_1^-)<v(i)$.  Hence, from
Lemma \ref{L:weakcriterion} applied
to the weak strip $\ws{w}{A}{v}$ and the consecutive $A$-bad
integers $r$ and $v(i)$, that $i<\aspm$, contradicting
\eqref{E:ar2}.
\end{proof}

\chapter{Reverse Local Rule} \label{ch:inverse} We describe an
algorithm to compute the inverse
$\psi=\psi_{u,v}:\OOL_{u,v}\rightarrow\IIL_{u,v}$ of the local rule
$\phi_{u,v}$ defined in the previous chapter.

\section{Reverse insertion at a cover} \label{SS:Rcover}

Let $C'=(\stc{v}{a,b}{x})$ be a marked strong cover. \textit{Reverse
insertion at $C'$} is a map that takes as input an initial pair
$(W',S_1)=(\wci{u}{A'}{x}{S_1}{\cdot})$ such that
$\out(W')=\out(C')$, and produces an initial pair of the form
$(W,S)=(\wci{w}{A}{v}{S}{\cdot})$, such that $\out(S)=\out(S_1)$ and
$\out(W)=\ins(C')$.

There are four cases, RA, RB, RC, and RX, which denote the inverses
of cases A, B, C, and X in the forward insertion algorithm.


We only need to specify $A$, $(i,j)$, and $S$, because $w=c_A^{-1}v$
and $C=\fc(S)=(\stc{w}{i,j}{\cdot})$. If $\size(S_1)>0$ let
\begin{align} \label{E:fcS1}
\xymatrix{%
{S_1:} & {u} \ar[r]^{\is,\js} & {z} \ar[r]^{\is^+,\js^+} & {\dotsm}
 }
\end{align}

\subsection{Commuting case} \

\medskip \noindent \textbf{Case RA (Commuting Case)}
If $(W',C')$ commutes, then we set $A = A'$, $(i,j) = (a,b)$, and
$S = C \cup S_1$.

\begin{example}
Let $n = 4, l = 0$ and $C' =
(\stc{[4,6,-3,3]}{-2,1}{[2,8,-3,3]})$. Consider the initial pair
$(W',S_1)$ where $W' = (\ws{[1,8,-2,3]}{\{1\}}{[2,8,-3,3]})$  and
$S_1 = (\stc{[1,8,-2,3]}{}{[1,8,-2,3]})$ has size 0.  Since the
pair $(W',C')$ commutes, we obtain the output initial pair $(W,S)$
where $$W = (\ws{[4,5,-2,3]}{\{1\}}{[4,6,-3,3]})$$ and $S =
(\stc{[4,5,-2,3]}{-2,1}{[1,8,-2,3]})$.
\end{example}


\subsection{Noncommuting cases} In the rest of the cases we assume
that the pair $(W',C')$ does not commute. By
Lemma~\ref{L:Rcyclicobstruction}, $\ba{u(b)-1} \in A'$ and we set
\begin{align} \label{E:RAhat}
 \Ah = A' - \{\ba{u(b)-1}\}.
\end{align}

We say that condition B (not to be confused with case B) holds if
there exists an $\Ah$-nice integer $q$ such that $u(b) < q$ and
$u^{-1}(q) \leq l$. If condition B holds let $q_B$ be the minimal
such $q$.

Say that condition C holds if $\size(S_1)>0$ (so that $(\is,\js)$ is
defined) and $u(\js)$ is $\Ah$-nice.
If condition C holds let $q_C = u(\js)$; in this case it follows
that $u(b)<u(\js)$ (see Lemma \ref{L:b1bound}).

Note that if conditions B and C both hold then $q_B \ne q_C$ since
$\js>l$ by the straddling condition.

\medskip \noindent \textbf{Case RX}
Suppose condition B does not hold and $\size(S_1)=0$ (and so in
particular condition C does not hold). We set $w=u$, $A=\Ah$, and
$S=(\stc{w}{}{w})$.

\begin{example}
Let $n = 5, l = 0$ and $C' =
(\stc{[2,-5,6,9,3]}{-1,3}{[2,-5,4,11,3]})$.  Consider the initial
pair $(W',S_1)$ where $$W' =
(\ws{[2,-4,5,8,4]}{\{3,4,5\}}{[2,-5,4,11,3]})$$ and $S =
(\stc{[2,-4,5,8,4]}{}{[2,-4,5,8,4]})$ is empty.  Then $(W',C')$
does not commute and we have $\Ah = \{3,5\}$. Neither condition B
nor condition C holds and we have the output initial pair is
$(W,S)$ where $W = (\ws{[2,-4,5,8,4]}{\{3,5\}}{[2,-5,6,9,3]})$ and
$S = S_1$.
\end{example}

\medskip \noindent \textbf{Case RB}
Suppose condition B holds, and, in addition, either condition C
does not hold or $q_B<q_C$. Set $q=q_B$ and let $p$ be the maximum
$\Ah$-nice integer such that $p<q$. Let $ A = \Ah \cup
\{\overline{q-1}\}$, $(i,j)= (u^{-1}(q),u^{-1}(p))$, and $S= C
\cup S_1$.

\begin{example}
Let $n = 6$, $l = 0$ and $$C' =
(\stc{[4,-1,1,12,-3,8]}{\{0,1\}}{[2,-1,1,12,-3,10]}).$$  Consider
the initial pair $(W',S_1)$ given by $$W' =
(\ws{[3,0,1,11,-2,8]}{\{2,3,5\}}{[2,-1,1,12,-3,10]})$$  and $S_1 =
(\stc{[3,0,1,11,-2,8]}{}{[3,0,1,11,-2,8]})$.  The pair $(W',C')$
does not commute and we have $\Ah = \{3,5\}$.  Condition B holds
with $q = q_B = 5$ and $p = 3$. Thus $A = \{3,4,5\}$, $(i,j) =
(-2,1)$ and the output initial pair $(W,S)$ is given by $W =
(\ws{[5,0,1,9,-2,8]}{\{3,4,5\}}{[4,-1,1,12,-3,8]})$ and $S =
(\stc{[5,0,1,9,-2,8]}{-2,1}{[3,0,1,11,-2,8]})$.
\end{example}

\medskip \noindent \textbf{Case RC}
Suppose condition C holds, and, in addition, either condition B
does not hold or $q_C<q_B$. Set $q=q_C$ and let $p$ be the maximum
$\Ah$-nice integer such that $p<q$. We set $A = \Ah \cup
\{\overline{q-1}\}$, $(i,j) = (\is,\js)$, and $S$ is obtained by
inserting $(\is,z^{-1}(p))$ into $S_1$ \textit{after} the first
reflection.

\begin{example}
Let $n = 4, l = 0$ and $C' = (\stc{v =
[4,5,-2,3]}{-2,7}{[4,6,-3,3]})$.  Consider the initial pair
$(W',S_1)$ given by $W' = (\ws{u =
[4,5,-2,3]}{\{1\}}{[4,6,-3,3]})$ and $S_1
=(\stc{[4,5,-2,3]}{-2,1}{z = [1,8,-2,3]})$.  The pair $(W',C')$
does not commute.  We have $u(b) = 2$ and $\Ah = \emptyset$.
Condition B does not hold but  Condition C holds with $q_C = u(1)
= 4$.  We have $p = 3$, $A =\{3\}$ and $z^{-1}(3) = 4$ giving us
the output initial pair $(W,S)$ where $W =
(\ws{[3,5,-2,4]}{\{3\}}{[4,5,-2,3]})$ and
$$S = ([3,5,-2,4] \overset{-2,1}{\longrightarrow} [1,7,-2,4]
\overset{-2,4}{\longrightarrow}[1,8,-2,3]).$$

\end{example}

\section{The reverse local rule} \label{SS:Rlocalrule} The reverse
algorithm applied to a final pair $(W',S')\in\OOL_{u,v}$ consists of
$m':=\size(S')$ steps, one for each cover in $S'$. The each step of
the reverse algorithm is called a reverse insertion. Each reverse
insertion takes as its input an \textit{initial} pair and produces
another as output. Write $S'=(C'_1,C'_2,\dotsc,C'_{m'})$. We
initialize $(W_{(m')},S_{(m')})=(W',\stc{u}{}{u})$, where
$u=\ins(W')$. For $k$ going from $m'$ down to $1$, we compute the
reverse insertion at $C'_k$ on the initial pair $(W_{(k)},S_{(k)})$
(which has the property that $(W_{(k)},C'_k)$ is a final pair),
which produces an initial pair $(W_{(k-1)},S_{(k-1)})$ such that
$\out(S_{(k-1)})=u$. Let $\psi_{u,v}(W',S')=(W,S,e)$ where
$(W,S)=(W_{(0)},S_{(0)})$ and $e=\size(S')-\size(S)$ is the number
of times Case RX occurred.

\section{Proofs for the reverse insertion} We want to compute
$\psi_{u,v}(W',S')$ with $\size(S')=m'$. By induction we may assume
that all of the reverse insertions have been performed except the
last step, which computes the reverse insertion on
$(W_{(m'-1)},S_{(m'-1)})$ at $C'_{m'}$, resulting in
$(W_{(m')},S_{(m')})$. Again we change notation, forgetting the
global meaning of $u,v,W',S'$. We denote this last reverse insertion
step as the reverse insertion on $(W',S_1)$ at $C'$, resulting in
$(W,S)$. We write
\begin{equation} \label{E:Rnotation}
\begin{split}
&(W',C')=(\wcf{u}{A'}{x}{v}{a,b}) \\
&\xymatrix{%
{S':\,\,} & {v} \ar[r]^{a,b}_{C'} & {x\,}
\ar[r]^{\ap,\bp}_{\Cpp} & {x^+} \ar[r]^{a^{++},b^{++}}_{{C'}^{++}} & {\dotsm} } \\
&\xymatrix{%
{S_1:} & {\hphantom{x}\,} & {u\,} \ar[r]^{\is,\js} & {z}
\ar[r]^{\is^+,\js^+} &
{\dotsm}} \\ %
&\xymatrix{%
{S:\,\,} & {w} \ar[r]^{i,j}_C & {u'} \ar[r]^{\ip,\jp}_{\Cp} & {z}
\ar[r]^{i^{++},j^{++}}_{C^{++}} & {\dotsm}} \\ %
&(W,C) = (\wci{w}{A}{v}{i,j}{u'})
\end{split}
\end{equation}
In Cases RA, RB, and RC we have $\size(S)>0$ so that $C=\fc(S)$ is
well-defined. In Case RX we make the convention that $u'=w$ and
write $(W,\vn)$ instead of $(W,C)$ where $\vn=(\stc{w}{}{w})$ is the
empty strong strip going from $w$ to itself.

The following inductive hypothesis will be useful. In each case it
must be re-established.

\begin{property} \label{pro:markmove2}\
\begin{enumerate}
\item[(i)] In Case RA, $m(C) = c_{A'}^{-1}(m(C'))$. In Case RB,
$m(C)>m(C')$. In Case RC, $m(C) \ge c_{\Ah}^{-1}(m(C'))$.
\item[(ii)] $m(C)>x(b)$.
\item[(iii)] Case RC cannot be preceded by Case RB.
\item[(iv)] The initial pair $(W,C)$ commutes in Cases RA and RC
and does not in Case RB.
\end{enumerate}
\end{property}

\subsection{Case RA} By Lemma \ref{L:initfincommute}
$(\wci{w}{a,b}{u}{A'}{v})$ is a commuting initial pair such that the
diagram commutes.
\begin{align*}
\xymatrix{%
{w}\ar@{-->}[r]^{a,b} \ar@{-->}[d]_{A'} & {u} \ar[d]^{A'}  \\
{v}\ar[r]_{a,b} & {x}}
\end{align*}
By Lemma \ref{L:RAstrong} $S$ is a strong strip.

For the proofs below, for Case RA we specialize the general notation
of \eqref{E:Rnotation} as follows.
\begin{align} \label{E:RAnotation}
&\xymatrix{%
{S:\,\,} & {w} \ar[r]^{a,b}_C & {u} \ar[r]^{\is,\js}_{\Cp} & {z} \ar[r] & {\dotsm}} \\ %
&\xymatrix{%
{S':\,\,} & {v} \ar[r]^{a,b}_{C'} & {x\,} \ar[r]^{\ap,\bp}_{\Cpp} &
{\xp} \ar[r] &  {\dotsm} }
\end{align}

\begin{lem}
In Case RA, Property~\ref{pro:markmove2} holds.
\end{lem}
\begin{proof} (iv) was proved above.
For (i) we have $m(C)=u(a)=c_{A'}^{-1}(x(a))=c_{A'}^{-1}(m(C'))$.
For (ii) we need to show $x(b)<u(a)$. Since $\stc{v}{a,b}{x}$ and
$\stc{w}{a,b}{u}$ are
strong covers, we have  $x(b)<x(a)$ and $u(b)<u(a)$.
If $u(b)$ is not $A'$-nice or is $A'$-bad then by \eqref{E:cvalue},
$x(b)=c_{A'}(u(b))\le u(b) <u(a)$ as desired. Otherwise by Lemma
\ref{L:cycdecrinv}, $x(b)<u(a)$ as desired.
\end{proof}

\begin{lem}
\label{L:RAstrong} In Case RA, $S$ is a strong strip.
\end{lem}
\begin{proof} We use the notation \eqref{E:RAnotation}.
Since $S_1$ is a strong strip by induction and $S=C\cup S_1$ we need
only show that
\begin{align} \label{E:RAstrong}
u(a) = m(C) < m(\Cp) = u(\js).
\end{align}
It suffices to show that %
\begin{align} \label{E:RAstrongsuff}
  u(\js)>x(a),
\end{align}
since $x(a) =c_{A'}(u(a)) \ge u(a)-1$ and $a\le l<\js$. Since $S'$
is a strong strip,
\begin{align} \label{E:RAS'strong}
x(a) = m(C') < m(\Cpp) = x(\bp).
\end{align}
We apply Property~\ref{pro:markmove2}(i) to the previous step.
Suppose the previous step was Case RA. Then $\js=\bp$ and
$x(a)<x(\js)$. Since $a\le l < \js$, \eqref{E:RAstrong} holds by
Lemma \ref{L:weaknoinv} applied to the weak strip $\ws{u}{}{x}$.
Suppose the previous step was Case RB. Then $m(\Cp)>m(\Cpp)>m(C')$
gives \eqref{E:RAstrongsuff}. Suppose the previous step was Case RC.
We have
\begin{align} \label{E:xapxbp}
  x(\bp)-x(\ap)<n
\end{align}
by Lemma \ref{L:Rcyclicobstruction} applied to the noncommutative
final pair $(\wcf{z}{}{\xp}{x}{\ap,\bp})$ at the beginning of the
previous step.

Suppose first that $\ba{a}=\ba{\ap}$. Let $a=\ap+kn$. By
\eqref{E:RAS'strong} and \eqref{E:xapxbp} it follows that $k \le0$.
Then $x(a)=x(\ap)+kn=\xp(\bp)+kn \le \xp(\bp)<m(\Cp)$, the last step
holding by Property \ref{pro:markmove2}(ii) applied to the previous
step. This gives \eqref{E:RAstrongsuff}.

Suppose next that $\ba{a}=\ba{\bp}$. Since $a\le l<\bp$ it follows
that $a=\bp+kn$ for some $k<0$. Then, since \eqref{E:xapxbp} gives $x^+(a^+)-x^+(b^+) < n$, we have
$$
x(a) =x(b^+)+kn=x^+(a^+)+kn \leq x^+(b^+)+(k+1)n \leq x^+(b^+)<m(C^+) \, ,
$$
by
Property \ref{pro:markmove2}(ii).

Finally suppose that $\ba{a}\not\in\{\ba{\ap},\ba{\bp}\}$. Then
$m(C')=x(a)=\xp(a)$. Suppose $\xp(a)\ge \xp(\bp)$. Equality cannot
hold since $a\le l<\bp$ so $\xp(a)>\xp(\bp)$. By
\eqref{E:RAS'strong} $\xp(\ap)>\xp(a)$.
By Lemma
\ref{L:weakcriterion} for the weak strip $\ws{z}{}{\xp}$ and the
noncommutativity of the final pair $(\wcf{z}{}{\xp}{x}{\ap,\bp})$,
we have $\ap<a$. But this contradicts Lemma \ref{L:strongcover} for
the strong cover $\stc{x}{\ap,\bp}{\xp}$.
\end{proof}

\subsection{Reverse noncommuting cases} By Lemma
\ref{L:Rremoveresidue}, $\ws{u}{\Ah}{v}$ is a weak strip such that
the diagram commutes.
\begin{align*}
\xymatrix{%
{} & {u} \ar@{-->}[dl]_{\Ah} \ar[d]^{A'} \\
{v} \ar[r]_{a,b} & {x}
}%
\end{align*}

\subsection{Case RX} In this case we have $w=u$ and $A=\Ah$ so that
$W=(\ws{w}{A}{v})=(\ws{u}{\Ah}{v})$ is a weak strip.

\subsection{Case RB} Thanks to Lemma
\ref{L:inverseExternalInsertion}, we may apply Lemma
\ref{L:Raddresidue}, which says that $(\wci{w}{A}{v}{i,j}{u})$ is a
noncommuting initial pair such that the diagram commutes.
\begin{align*}
\xymatrix{%
{w} \ar[r]^{i,j} \ar[d]_A & {u} \ar@{-->}[dl]_{\Ah} \ar[d]^{A'} \\
{v} \ar[r]_{a,b} & {x}
}%
\end{align*}
This case is finished by Lemma \ref{L:RBstrong} which shows that $S$
is a strong strip.

\begin{lem}
\label{L:inverseExternalInsertion} Suppose condition B holds. Let
$q=q_B$ and $p$ be the maximum $\Ah$-nice integer such that $p<q$. Then
\begin{align}\label{E:pqinv}
p_0 := u(b) &\le p < q \le u(a) + n \\
\label{E:iBppos}
  u^{-1}(q) &\le l < u^{-1}(p).
\end{align}
\end{lem}
\begin{proof} $p_0$ is $\Ah$-nice by construction.
Also $u^{-1}(p_0)=b>l$. Since $q>p_0$ it follows that $p \ge p_0$
and that \eqref{E:iBppos} holds. For the upper bound on $q$, suppose
$q>u(a)+n$. Then $q':=q-n>u(a)$ is $\Ah$-nice. Since $(W',C')$ does
not commute, $u(a)<u(b)$ and by Lemma \ref{L:Rcyclicobstruction}
these are consecutive $\Ah$-nice integers. Therefore $q' \ge u(b)$.
Also $u^{-1}(q')=u^{-1}(q)-n\le l$, contradicting the minimality of
$q$. This proves the upper bound in \eqref{E:pqinv}.
\end{proof}

For the proofs below, for Case RB we specialize the general notation
of \eqref{E:Rnotation} as follows.
\begin{align*}
&\xymatrix{%
{S:\,\,} & {w} \ar[r]^{i,j}_C & {u} \ar[r]^{\is,\js}_{\Cp} &  {\dotsm}} \\ %
&\xymatrix{%
{S':\,\,} & {v} \ar[r]^{a,b}_{C'} & {x\,} \ar[r]^{\ap,\bp}_{\Cpp} &
{\dotsm} }
\end{align*}
where $(i,j)=(u^{-1}(q),u^{-1}(p))$.

\begin{lem}
In Case RB, Property~\ref{pro:markmove2} holds.
\end{lem}
\begin{proof} (iv) was already shown. (ii) follows from (i) since
$m(C)>m(C')=x(a)>x(b)$, the last inequality holding because
$\stc{v}{a,b}{x}$ is a strong cover. For (i) it is equivalent to
show that $q>c_{A'}(u(a))$. By \eqref{E:pqinv} the only difference
between $\Ah$-niceness and $A'$-niceness for integers in the
interval $[u(a),u(a)+n]$, is that $u(b)$ is $\Ah$-nice but not
$A'$-nice. Since $(W',C')$ does not commute, by Lemma
\ref{L:Rcyclicobstruction}, $u(a)$ is $A'$-nice and the minimum
$A'$-nice integer $r$ such that $r>u(a)$, satisfies $r>u(b)$. By
\eqref{E:pqinv} we see that $q$ is $A'$-nice and $q>u(b)$. Therefore
$r\le q$ and from \eqref{E:cvalue2}, $c_{A'}(u(a))=r-1<q$ as desired.
\end{proof}

\begin{lem} \label{L:RBstrong} In Case RB, $S$ is a strong strip.
\end{lem}
\begin{proof} By induction and the construction of $S$, it suffices
to show that $m(C)<m(\Cp)$.

We have $m(C) = q < u(\js) = m(\Cp)$, which follows from Lemma
\ref{L:inversebounds}, as $u(\is) < u(\js)$ holds since
$\Cp=(\stc{u}{\is,\js}{\cdot})$ is a strong cover.
\end{proof}

\begin{lem}
\label{L:b1bound} Suppose $\size(S_1)>0$ so that $(\is,\js)$ is
defined. Then
\begin{align} \label{E:RCubj1}
u(b) \le x(a)<u(\bp)\le u(\js).
\end{align}
In particular $\js\ne b$.
\end{lem}
\begin{proof}
By Lemma \ref{L:Rcyclicobstruction}, $u(a)$ is $A'$-nice and the
first inequality in \eqref{E:RCubj1} holds. Let $S'$ be as in
\eqref{E:Rnotation}. Since $S'$ is a strong strip we have
$$c_{A'} u(a)=x(a) = m(C') < m(\Cpp) = x(\bp)=c_{A'} u(\bp).$$ Since
$u(a)$ is $A'$-nice and $c_{A'}$ is cyclically decreasing, it
follows that the second inequality in \eqref{E:RCubj1} holds. By
Property \ref{pro:markmove2} (ii) for the previous step we have
$u(\js)>x(\bp) \ge u(\bp)-1$, proving the last inequality in
\eqref{E:RCubj1}.
\end{proof}

\begin{lem}\label{L:inversebounds}
Suppose that $S_1$ is nonempty so that $(\is,\js)$ is defined. Let
$q=q_B$ if Case RB holds and $q=q_C$ if Case RC holds. Let $p$ be
the maximum $\Ah$-nice integer with $p<q$. Then
\begin{enumerate}
\item If $u(b) < u(\is)$ then Case RB holds and
\begin{equation}\label{E:pqinv2}
u(b) \leq p < q \le u(\is).
\end{equation}
\item Otherwise we have $u(\is) < u(b)<u(\js)$ and
\begin{equation}\label{E:pqinv3}
u(b) \leq p < q \leq u(\js).
\end{equation}
\end{enumerate}
Moreover, whenever Case RB holds we have $q \le u(a)+n$.
\end{lem}


\begin{proof}[Proof of Lemma~\ref{L:inversebounds}] Since
$\fc(S_1)=(\stc{u}{\is,\js}{\cdot})$ is a strong cover,
$u(\is)<u(\js)$. Lemma~\ref{L:b1bound} implies that either
$u(b)<u(\is)$ or $u(\is)<u(b)<u(\js)$.

For either \eqref{E:pqinv2} or \eqref{E:pqinv3} when Case RB holds,
the argument that $q \leq u(a) + n$ is the same as in the proof of
Lemma~\ref{L:inverseExternalInsertion}.

Suppose $u(b)<u(\is)$. Suppose there is no $\Ah$-nice integer $q'$
such that $u(b)<q'\le u(\is)$. Applying Lemma \ref{L:weakcriterion}
to the weak strip $\ws{u}{\Ah}{v}$ and the $\Ah$-nice integer $u(b)$
and the next larger one (which is greater than $u(\is)$), we have $b
< \is$ which contradicts the straddling inequality $\is\le l< b$. So
let $q'$ be the maximum $\Ah$-nice integer such that $u(b)<q'\le
u(\is)$. If $u(\is)$ is $\Ah$-nice then $q'=u(\is)$ and
$u^{-1}(q')=\is\le l$. Otherwise $u(\is)$ is not $\Ah$-nice and
$q'<u(\is)$. Let $p'$ be the maximum $\Ah$-nice integer with
$p'<q'$. Applying Lemma \ref{L:weakcriterion} to the weak strip
$\ws{u}{\Ah}{v}$ and consecutive $\Ah$-nice integers $p'<q'$, we
have $u^{-1}(q')\le \is\le l$. Therefore Case B holds.

Suppose $u(\is) < u(b) < u(\js)$. Let $q'$ be the maximum $\Ah$-nice
integer with $q'\le u(\js)$; it exists and satisfies $u(\is) < q'$
since $u(b)$ is $\Ah$-nice.

If $q' = u(\js)$ then either Case RB with $q < q'$ will hold or Case
RC will hold for $q = q'$. \eqref{E:pqinv3} follows.

Now suppose $q' < u(\js)$.  By Lemma~\ref{L:weakcriterion} applied
to $W'$, we have $u^{-1}(q') < \js$. But by
Lemma~\ref{L:strongcover} applied to $\fc(S_1)$ we must have
$u^{-1}(q') < \is \le l$. Therefore $q'$ is a witness that Condition
B holds. (2) follows.
\end{proof}

\subsection{Case RC} Again we sketch the proof and then fill in the
Lemmata which prove the details.

By Lemma \ref{L:Rremoveresidue},
$\ws{u}{\Ah}{v}$ is a weak strip such that the diagram commutes:
\begin{align*}
\xymatrix{%
{} & {u} \ar@{-->}[dl]_{\Ah} \ar[r]^{\is,\js}\ar[d]^{A'} & {z} \\
{v} \ar[r]_{a,b} &{x}&{} }%
\end{align*}
By Lemma \ref{L:RCcom1},
$(\wci{u}{\Ah}{v}{\is,\js}{z})$ is a commuting initial pair.
Let $x' = v \,t_{\is,\js} = c_{\Ah}z$. By
Lemma \ref{L:initfincommute} $(\wcf{z}{\Ah}{x'}{v}{\is,\js})$ is a
commuting final pair such that the diagram commutes.
\begin{align*}
\xymatrix{%
{} & {u}\ar[dl]_{\Ah}\ar[r]^{\is,\js} & {z} \ar@{-->}[dl]^{\Ah} \\
{v} \ar@{-->}[r]_{\is,\js} &{x'}&{} }%
\end{align*}
Recall that $q_C=u(\js)=z(\is)$. Let $(\ip,\jp)= (\is,z^{-1}(p))$
and $u' = z \,t_{\ip,\jp}=t_{p,q_C} z$. By Lemma \ref{L:RCzinvp} we
may apply Lemma \ref{L:Raddresidue} to the weak strip
$\ws{z}{\Ah}{v}$ and pair $p<q_C$ of consecutive $\Ah$-nice
integers, so that $(\wci{u'}{A}{x'}{\ip,\jp}{z})$ is a noncommuting
initial pair such that the diagram commutes.
\begin{align*}
\xymatrix{%
{} & {u'} \ar@{-->}[r]^{\ip,\jp}\ar@{-->}[d]_{A} & {z} \ar[dl]^{\Ah} \\
{v} \ar[r]_{\is,\js} & {x'} & {} }
\end{align*}
By Lemma \ref{L:RCcom2} the final pair
$(\wcf{u'}{A}{x}{v}{\is,\js})$ commutes. Define $w=
u'\,t_{\is,\js}=c_{A}^{-1}(v)$. By Lemma \ref{L:initfincommute},
$(\wci{w}{A}{v}{\is,\js}{u'})$ is a commuting initial pair such that
the diagram commutes.
\begin{align*}
\xymatrix{%
{w} \ar@{-->}[r]^{\is,\js} \ar@{-->}[d]_{A} & {u'} \ar[r]^{\ip,\jp}\ar[d]_{A} & {z} \ar[dl]^{\Ah} \\
{v} \ar[r]_{\is,\js} & {x'} & {} }
\end{align*}
In particular, as defined, $S$ is a strong tuple and $W$ is a weak
strip. Lemma \ref{L:RCstrong} shows that $S$ is a strong strip.

\begin{lem} \label{L:RCcom1} In Case RC, $v(\is)<v(\js)$.
\end{lem}
\begin{proof}
Since $\stc{u}{\is,\js}{z}$ is a strong cover, $u(\is)<u(\js)$. Now
$v = c_{\Ah} u$ and $u(\js)$ is $\Ah$-nice since Condition C holds.
It follows that $v(\is)<v(\js)$.
\end{proof}

\begin{lem} \label{L:RCzinvp} In Case RC, $z^{-1}(p)>l$.
\end{lem}
\begin{proof} Suppose not, that is, $z^{-1}(p)\le l$. By Lemma
\ref{L:inversebounds}, $u(b)<u(\js)=q_C$. Now $p<q_C$ are
consecutive $\Ah$-nice integers. Since $u(b)$ is also $\Ah$-nice it
follows that $u(b)\le p$. Suppose $u(b)=p$.  Since $z^{-1}(p) = t_{i_1 j_1} u^{-1}(p)=u^{-1}(p)$
given that $u(j_1)-u(i_1)<n$ and $u(i_1) <p< u(j_1)$ from Lemma~\ref{L:Fcyclicobstruction}
on the non-commuting initial pair $(\wci{u}{A'}{x}{i_1,j_1}{z})$, we have
the contradiction $l<b=u^{-1}(u(b))=u^{-1}(p)=z^{-1}(p) \le l$. So $u(b)<p$.
But then condition B is satisfied by the integer $p$ which is less
than $q_C$, meaning that Case RB holds, which is a contradiction.
\end{proof}

\begin{lem} \label{L:RCcom2}
In Case RC, $u'(\is)>u'(\js)$.
\end{lem}
\begin{proof} Since $u'(i_1)=t_{p , q_C}z(i_1)=t_{p, q_C}(q_C)=p$, and $u'(j_1)=t_{p ,q_C}z(j_1)=t_{p, q_C}u(i_1)$,
this is equivalent to
\begin{align} \label{E:ptui1}
p> t_{p,q_C} u(\is).
\end{align}
Since we are in Case RC, by Lemma \ref{L:inversebounds} we have
$u(\is)<u(b)\le p < q_C=u(\js)$. If
$\ba{u(\is)}\not\in\{\ba{p},\ba{q_C}\}$ then \eqref{E:ptui1} clearly
holds. Otherwise $u(\is)\in\{p+kn,q_C+kn\}$ for some integer $k<0$.
Then $t_{p,q_C}u(\is)\in \{p+kn,q_C+kn\}$. Then \eqref{E:ptui1}
follows from $p+kn<q_C+kn\le q_C-n<p$, the last inequality holding
since $p$ and $q_C$ are consecutive $\Ah$-nice integers.
\end{proof}

For the rest of the proofs we use the following notation for Case
RC.
\begin{equation} \label{E:RCnotation}
\begin{split}
&\xymatrix{%
{S':\,}&{v} \ar[r]^{a,b}_{C'} & {x\,} \ar[r]^{\ap,\bp}_{\Cpp} & {x^+} \ar[r]^{a^{++},b^{++}}_{{C'}^{++}} & {\dotsm}} \\
&\xymatrix{%
{S_1:}& {\hphantom{x}} & {u\,} \ar[r]^{\is,\js} & {z}
\ar[r]^{\is^+,\js^+} & {}}
 \\ %
&\xymatrix{%
{S:\,}&{w} \ar[r]^{\is,\js}_C & {u'} \ar[r]^{\is,\jp}_{\Cp} & {z} \ar[r]^{\is^+,\js^+}_{C^{++}} & {}}  %
\end{split}
\end{equation}
where $\jp=z^{-1}(p)$.

\begin{lem}
Property~\ref{pro:markmove2} holds in Case RC.
\end{lem}
\begin{proof} We use the notation \eqref{E:RCnotation}.
Since $m(C)=u'(i_1)=p$, as seen above equation \eqref{E:ptui1}, (i) is equivalent to $p \ge u(b)$
which holds by Lemma \ref{L:inversebounds}.
(ii) then follows since the non-commutativity of the final pair
$(\wci{u}{A'}{v}{a,b}{x})$ implies from Lemma~\ref{L:Rcyclicobstruction}
that $u(b)$ is not $A'$-nice, and thus that $x(b)=c_{A'}u(b)=u(b)-1<p$ from \eqref{E:cvalue2}.
(iv) was shown above. For
(iii), suppose Case RC was preceded by Case RB. We use the following
diagram.
\begin{align*}
\xymatrix{%
  {} & {u} \ar[r]^{\is,\js} \ar[d]^{A'} \ar[dl]_{\Ah} & {z}\ar[d]^{A^+}  & {\cdots} \\
  {v}\ar[r]_{a,b} & {x}\ar[r]_{\ap,\bp} & {\xp} & {\cdots}
}%
\end{align*}
By assumption the final pairs $(W',C')=(\wcf{u}{A'}{x}{v}{a,b})$ and
$(W^+,\Cpp)=(\wcf{z}{A^+}{\xp}{x}{\ap,\bp})$ both do not commute.

We will see that it suffices to show that
\begin{align} \label{E:RCpropsuff}
\text{There is an integer $i'\le l$ such that $u(b)<u(i')<u(\js)$.}
\end{align}
Let $i'$ be as in \eqref{E:RCpropsuff}. Then $u(i')$ is not
$\Ah$-nice; otherwise Condition B would hold for $u(i')$. Let $r<r'$
be the pair of consecutive $\Ah$-nice integers such that $r <
u(i')<r'$. Since $u(b)$ and $u(\js)$ are $\Ah$-nice (the latter by
the Case RC assumption) we have $u(b)\le r < u(i')<r'\le u(\js)$. By
Lemma \ref{L:weakcriterion} for the weak strip $\ws{u}{\Ah}{v}$ we
have $u^{-1}(r)<i'\le l$. This gives the contradiction that
Condition B holds for $r$ if $u(b)<r$. But if $u(b)=r$ then
$l<b=u^{-1}(u(b))=u^{-1}(r)\le l$, a contradiction.

We now prove \eqref{E:RCpropsuff}. Since $S'$ is a strong strip,
\begin{align} \label{E:RCproppresand}
  x(b)<x(a)<x(\bp).
\end{align}
Suppose $a\ne \ap$. We claim that \eqref{E:RCpropsuff} holds for
$i'=\ap$. By Property \ref{pro:markmove2}(ii) for the previous step,
we have $u(\js)>\xp(\bp)=x(\ap) \ge u(\ap)-1$, that is, $u(\js) \ge
u(\ap)$. But $\ap\le l<\js$ so $u(\js)>u(\ap)$, giving the right
hand inequality in \eqref{E:RCpropsuff} for $i'=\ap$. We are done if
$u(b)<u(\ap)$. Suppose not. Since $\ap\le l<b$ we must have
$u(b)>u(\ap)$. We claim that
\begin{align} \label{E:RCpropsand}
  x(\ap)<x(b)<x(a)<x(\bp).
\end{align}
Since $(W',C')$ is a noncommuting final pair, by Lemma
\ref{L:Rcyclicobstruction}, $u(a)$ is $A'$-nice, $u(b)$ is not
$A'$-nice, and $x(a)\ge u(b)>u(b)-1=x(b)$. If $u(a)<u(\ap)$ then
$u(\ap)$ is not $A'$-nice and $x(\ap)=u(\ap)-1<u(b)-1=x(b)$ so that
\eqref{E:RCpropsand} holds. If $u(\ap)<u(a)$ then since $u(a)$ is
$A'$-nice, $x(\ap)<u(a)\le x(b)$ and again \eqref{E:RCpropsand}
holds. By Lemma \ref{L:strongcover} applied to the strong cover
$\Cpp$ we have $a<\ap$. Now $0<x(\bp)-x(\ap)<n$ by Lemma
\ref{L:Rcyclicobstruction} for the noncommuting final pair
$(W^+,\Cpp)$. Therefore \eqref{E:RCpropsand} gives
$\xp(\bp)<\xp(b)<\xp(a)<\xp(\ap)$. By Lemma \ref{L:weakcriterion}
for the weak strip $W^+$ we have $\ap<a$, a contradiction.

Suppose $a=\ap$. We claim that \eqref{E:RCpropsuff} holds for
$i'=\is$. Since $\stc{u}{\is,\js}{z}$ is a strong cover we have
$u(\is)<u(\js)$, the right hand inequality in \eqref{E:RCpropsuff}
for $i'=\is$. It suffices to show that
\begin{align*}
  u(\is) = z(\js) \ge z(\bp) > z(\bp)-1 = \xp(\bp) = x(\ap)=x(a) \ge
  u(b).
\end{align*}
The first inequality holds by \eqref{E:pqinv2} for the previous
(Case RB) step. The second equality holds by Lemma
\ref{L:Rcyclicobstruction} for the final pair $(W^+,\Cpp)$. The last
inequality holds by Lemma \ref{L:Rcyclicobstruction} for the
noncommutative final pair $(W',C')$.
\end{proof}

\begin{lem} \label{L:RCstrong}
In Case RC, $S$ is a strong strip.
\end{lem}
\begin{proof} By the construction of $S$ and the fact that
$S_1$ is a strong strip by induction, we need only show that
\begin{align*}
  m(C) < m(\Cp) < m(C^{++}).
\end{align*}
The second inequality holds because
$m(\Cp)=z(i^+)=z(\is)=q_C=m(\fc(S_1))<m(C^{++})$, the inequality holding
by the strong strip condition of $S_1$. Given that $m(C)=p$, the first inequality holds
since  $p<q_C$ follows from Lemma~\ref{L:inversebounds}.
\end{proof}

\chapter{Bijectivity} \label{ch:bijectivity}


Our main theorem is the following.
\begin{thm}
\label{thm:ins} The maps $\psi$ and $\phi$ are inverses to each
other.  Thus the map $\phi$ is a bijection.
\end{thm}

Our approach to proving bijectivity is to reduce to the case of at
most two steps, exploiting the fact that the maps $\phi$ and $\psi$
can be ``factorized" into ``smaller" instances of $\phi$ and $\psi$,
to which induction may be applied.

Consider the sequence of steps involved in computing $\phi(W,S,e)$.
With respect to this sequence, we call a subsequence of consecutive
steps \textit{irreducible} if it consists of:
\begin{enumerate}
\item a Case X step.
\item a Case A step followed by some maximum number $m$ (possibly zero)
of consecutive Case C steps.
\item a Case B step.
\end{enumerate}
Mnemonically we denote such irreducible sequences by $\IX$,
$\IA(m)$, and $\IB$ respectively.

Dually, consider the sequence of steps involved in computing
$\psi(W',S')$. With respect to this sequence, we call a subsequence
of consecutive steps \textit{irreducible} if it consists of:
\begin{enumerate}
\item a Case RX step.
\item a Case RA step followed by some maximum
number $m$ (possibly zero) of Case RC steps.
\item a Case RB step.
\end{enumerate}
Denote these mnemonically by $\IX^{-1}$, $\IA(m)^{-1}$ and
$\IB^{-1}$ respectively. As a warning, note that, for example, when
$m=2$ the inverse of Case A followed by two Case Cs, is Case RA
followed by two Case RCs.

It is clear that every sequence has a unique factorization into
irreducible subsequences.

We shall show in Section \ref{sec:bijext} that if the forward
algorithm ends with $\IX$ then $\psi\circ\phi=\id$ and if the
reverse algorithm begins with $\IX^{-1}$ then $\phi\circ\psi=\id$.

In Section \ref{sec:bija} we show that if the forward algorithm
starts with $\IA(0)$ or the reverse algorithm ends with
$\IA^{-1}(0)$ then we have bijectivity.

In Section \ref{sec:bijb} we show bijectivity when the forward
algorithm begins with $\IB$ and when the reverse algorithm ends with
$\IB^{-1}$.

In Section \ref{sec:bijc} we show bijectivity when the forward
algorithm begins with $\IA(m)$ and when the reverse algorithm ends
with $\IA(m)^{-1}$ for $m>0$. This is accomplished by reducing to
the case that $m=1$.

This covers all cases, so Theorem \ref{thm:ins} follows.

In each of the following sections, for the proof of
$\psi\circ\phi=\id$, we suppose that $(W,S,e)\in \IIL_{u,v}$ and
denote $(W',S')=\phi(W,S,e)$ and $\IS$ for the sequence of steps in
this computation. For the proof of $\phi\circ\psi=\id$ we suppose
$(W',S')\in\OOL_{u,v}$, write $(W,S,e)=\psi(W',S')$ and by abuse of
notation write  $\IS^{-1}$ to denote the sequence of reverse
insertion steps in the computation of $(W,S,e)$.

\section{External insertion} \label{sec:bijext} For
$\psi\circ\phi$, suppose $\IS$ has the form $\IS=\IS_1\IX$, that is,
it ends with $\IX$. The last (Case X) step depends only on the weak
strip entering that step, and $\IS_1$ is itself a ``smaller"
instance of $\phi$. By induction we may reduce to the case that
$\IS=\IX$ consists of a single Case X step.

For $\phi\circ\psi$, suppose $\IS^{-1}$ has the form
$\IS^{-1}=\IX^{-1}\IS_1^{-1}$. The first step, which is Case RX,
depends only on the last cover in $S'$. The second step will not be
Case RC, so $\IS_1^{-1}$ is a smaller instance of $\psi$. By
induction we may assume that $\IS^{-1}=\IX^{-1}$, so that $S'$
consists of a single cover $C'$.

So for $\psi\circ\phi$, let $(W,S)$ be a final pair such that
$W=(\ws{w}{A}{v})$ with $|A|<n-1$ and $S=(\stc{v}{}{v})=\vn$ the
empty strong strip. Let $W \mapsto (W',C')$ be the result of
external insertion with $(W',C')=(\wcf{u}{A'}{x}{v}{a,b})$,
$A'=A\cup\{\ba{q}\}$ and $(a,b)=(v^{-1}(q),v^{-1}(p))$ where $q<p$
are the pair of consecutive $A$-bad integers such that $v^{-1}(q)\le
l$ and $q$ is maximal. By construction the final pair $(W',C')$ is
noncommuting.
\begin{align*}
\xymatrix{%
{} & {u} \ar[dl]_A \ar@{-->}[d]^{A'} \\
{v}\ar@{-->}[r]_{a,b} & {x}
}%
\end{align*}
Applying the first step of $\psi$ to $(W',C')$, we cannot be in
Cases RA or RC. By Lemma~\ref{L:Fcyclicobstruction}, $u(b)$ is not $A'$-nice and thus
$u(b)-1 = c_{A'}(u(b))=x(b)=v(a)=q$.
Therefore $\Ah$ as defined in \eqref{E:RAhat}, equals our $A$. By
Lemma \ref{L:nicebad}, $u(q')$ is $A$-nice if and only if
$v(q')=c_A(u(q'))$ is $A$-bad. So condition B fails due to the
maximality of $q$ in the external insertion. Therefore Case RX
occurs and it produces the original weak strip $W$ as desired.

For $\phi\circ\psi$, let $S'=C'$ be a single strong cover and
$(W',C')=(\wcf{u}{A'}{x}{v}{a,b})$ a final pair such that Case RX
occurs. In particular $(W',C')$ does not commute. We define $A=\Ah$
as in \eqref{E:RAhat} and have the commutative diagram with weak
strip $W=(\ws{u}{A}{v})$.
\begin{align*}
\xymatrix{%
{} & {u} \ar@{-->}[dl]_A \ar[d]^{A'} \\
{v}\ar[r]_{a,b} & {x}
}%
\end{align*}
By Lemma \ref{L:Rcyclicobstruction} $u(a)<u(b)$ are consecutive
$\Ah$-nice integers and $u(b)$ is not $A'$-nice. We now apply
external insertion. Let $q<p$ be consecutive $\Ah$-bad integers such
that $v^{-1}(q)\le l$ with $q$ maximum. We have $u(b)-1=x(b)=v(a)$
so it suffices to show that $q=v(a)$ and $p=v(b)$. By Lemma
\ref{L:nicebad} $v(a)<v(b)$ are consecutive $A=\Ah$-bad integers and
$v^{-1}(v(a))=a\le l$. If there was an $\Ah$-bad integer $q'>v(a)$
with $v^{-1}(q')\le l$, then $q'>v(b)$ and $c_{\Ah}^{-1}(q')$ is an
$\Ah$-nice integer greater than $u(b)$ with
$u^{-1}(c_{\Ah}^{-1}(q'))=v^{-1}(q')\le l$, so that condition B
holds, contradicting the assumption that Case RX occurs. Therefore
$q=v(a)$ and $p=v(b)$ as desired.

\section{Case A (commuting case)} \label{sec:bija}

In the case that $\IS$ (resp. $\IS^{-1}$) consists of a single Case
A (resp. RA) step, bijectivity holds by Lemma
\ref{L:initfincommute}.

In general, for $\psi\circ\phi$, suppose $\IS$ starts with a Case A
step which is not followed by a Case C step. Write
$\IS=\IA(0)\IS_1$. Since the first step of $\IS_1$ is not Case C by
assumption, $\IS_1$ is an instance of $\phi$ involving fewer steps.
By induction we may assume that $\IS$ is a single Case A step, since
a single Case RA step is unaffected by the strong covers that were
produced previously.

For $\phi\circ\psi$ we may similarly reduce to the single RA step
case.

\section{Case B (bumping case):} \label{sec:bijb}

For $\psi\circ\phi$ we first suppose that the entire sequence $\IS$
is a single Case B step. Let $(W,C)=(\wci{w}{A}{v}{i,j}{u})$ be a
noncommuting initial pair. As in Case B we define
$\Ah=A-\{\ba{u(i)-1}\}$ and let $q<p$ be the pair of consecutive
$\Ah$-nice integers such that $q<u(j)$ and $u^{-1}(q)\le l$, with
$q$ maximum. Defining $A'=\Ah\cup\{p-1\}$,
$(a,b)=(u^{-1}(q),u^{-1}(p))$ and $x=v\,t_{a,b}=c_{A'} u$, we have
that $\ws{u}{\Ah}{v}$ is a weak strip and
$(W',C')=(\wcf{u}{A'}{x}{v}{a,b})$ is a noncommuting final pair such
that the diagram commutes:
\begin{equation*}
\xymatrix{%
{w} \ar[r]^{i,j} \ar[d]_A & {u} \ar@{-->}[dl]_{\Ah}
\ar@{-->}[d]^{A'} \\
{v} \ar@{-->}[r]_{a,b} & {x}}
\end{equation*}

We now apply the reverse algorithm to $(W',C')$. In the single cover
context Case RC does not occur, and $(W',C')$ does not commute. By
Lemma \ref{L:Rcyclicobstruction} $u(b)=p$ is not $A'$-nice. Then
$A'-\{\ba{u(b)-1}\}=A'-\{\ba{p-1}\}=\Ah$ so that the above
definition of $\Ah$ agrees with the one in the noncommutative case
of the reverse algorithm. Now $q=u(a)$ satisfies $q<u(j)$ by its
definition, so $u(a)<u(j)$. Since $(W',C')$ is a noncommuting final
pair, by Lemma \ref{L:Rcyclicobstruction} $u(a)<u(b)$ are
consecutive $\Ah$-nice integers. Since $(W,C)$ is a noncommuting
initial pair, by Lemma \ref{L:Fcyclicobstruction} $u(j)<u(i)$ are
consecutive $\Ah$-nice integers. Therefore we have $u(a)<u(b)\le
u(j)< u(i)$. In particular condition B holds with the $\Ah$-nice
integer $u(i)$. So Case RB holds, and say it selects the consecutive
$\Ah$-nice integers $u(b)\le p_B<q_B\le u(i)$. If $q_B<u(i)$ then
since $u(j)<u(i)$ are consecutive $\Ah$-nice integers and $j>l$, it
follows that $u(b)<q_B<u(j)$, contradicting the maximality of $q$ in
Case B. Therefore $q_B=u(i)$ and $p_B=u(j)$. It is now clear that
$\psi\circ\phi=\id$.

We now return to the general case of $\psi\circ\phi$. Let
$S=(\xymatrix{{w}\ar[r]^{i,j} & {u} \ar[r]^{i,j_1} \ar[r] &
{\dotsm}})$. We know that the output of the first (Case B) step is a
noncommuting final pair. By Property \ref{pro:markmove}(iii) the
second step of $\phi$ is not Case C. By induction we suppose that
performing the rest of $\phi$ and then all of $\psi$ except for the
last step, is the identity. We now consider the last step of the
reverse algorithm $\psi$; it is the noncommutative case. It suffices
to show that Case RB occurs in this last step. By the single cover
case we know that condition B holds with $q_B=u(i)$, so that Case RX
cannot occur. Since $S$ is a strong strip we have $u(i)<u(\js)$ so
that Case RC cannot hold. Therefore Case RB holds and we have
reduced to the single cover case.

The reasoning for $\phi\circ\psi$ in the single cover case, is
entirely similar to that for $\psi\circ\phi$ in the single cover
case. Suppose now that the last step of $\psi$ on $(W',S')\in\OOL$
is Case RB. After applying $\psi$ to $(W',S')$ we apply $\phi$. Its
first step is Case B and undoes the last step of $\psi$ by the
single cover case. The second step of $\phi$ cannot be Case C, and
by induction the rest of $\phi$ is the inverse of the rest of
$\psi$, and we are done.

\section{Case C (replacement bump)} \label{sec:bijc}

\subsection{Reduction to $m=1$}

For $\psi\circ\phi$ suppose $\IS=\IA(m)\IS_1$ with $m>0$ maximal.
Then $\IS_1$ does not start with Case RC and the output strong cover
of the last (Case C) step of $\IA(m)$ is involved in a commuting
final pair by Property \ref{pro:markmove}(iv). It follows that the
usual reduction works and we may assume that $\IS=\IA(m)$.

For $\phi\circ\psi$ suppose $\IS^{-1}=\IS_1^{-1}\IA(m)^{-1}$ for
$m>0$, that is, $\psi$ produces a sequence of steps that ends with
Case RA followed by some positive number of Case RC steps. Since
$\IS_1^{-1}$ is followed by a Case RA step, it is a smaller instance
of $\psi$. In particular its inverse cannot begin with Case C. Since
Case RA steps are unaffected by previously produced strong covers,
again by induction we may assume that $\IS=\IA(m)$.

Let us consider $\psi\circ\phi$ on $(W,S,0)$ such that $\phi$ on
$(W,S,0)$ consists of a single Case A step followed by $m$ Case C
steps for some $m>0$. By Section \ref{sec:bijext} the output
$(W',S')$ is such that $\psi$ on $(W',S')$ does not start with RX.
We shall reduce to the case $m=1$ by showing that $\phi$ can be
achieved by splicing together two operations: the first two steps of
$\phi$, which consists of Case A followed by a Case C, and another
application of $\phi$ which consists of a Case A step (the ``second
half" of the first Case C step) followed by $m-1$ Case C steps.

Suppose $m>1$. Let
\begin{align*}
  S &=
  (\xymatrix{{w}\ar[r]^{i,j} & {u_0} \ar[r]^{i,j_1} &{u_1}
  \ar[r]^{i,j_2} & \dotsm  \ar[r]^{i,j_m} &
  {u_m=u}}).
\end{align*}
Since $S$ is a strong strip we have
\begin{align} \label{E:FACstrong2}
  u_0(i) < u_0(j_1) = u_1(i) < u_1(j_2).
\end{align}
Let $W=(\ws{w}{A}{v})$ and $\phi(W,S,0)=(W',S')$; since this is a
Case A step followed by $m$ Case C steps, $S'$ has the form
\begin{align*}
  S' = [v;((a_1,b_1),(a_2,b_2),\dotsc,(a_m,b_m),(i,j));x].
\end{align*}
The computation of $\phi(W,S,0)$ starts with a two step computation
$$\phi(W,[w;((i,j),(i,j_1));u_1],0)=(W_1,[v;((a_1,b_1),(i,j));x_1]),$$
where $W_1=(\ws{u_1}{A_1}{x_1})$. By Property \ref{pro:markmove}(iv)
applied to the Case C step, the final pair
$(W_1,\stc{x_1\,t_{ij}}{i,j}{x_1})$ commutes. By Lemma
\ref{L:initfincommute},
$(\wci{u_1\,t_{ij}}{A^*}{x_1\,t_{ij}}{i,j}{u_1})$ is a commuting
initial pair such that the diagram commutes.
\begin{align*}
\xymatrix{%
{u_1\,t_{ij}} \ar[r]^{i,j} \ar[d]_{A^*}^{W^*} & {u_1}\ar[d]^{A_1} \\
{x_1\,t_{ij}} \ar[r]_{i,j} & {x_1}
}%
\end{align*}
Let
\begin{align*}
W^* &=(\ws{u_1\,t_{ij}}{A^*}{x_1\,t_{ij}}) \\
S^*&=[u_1\,t_{i,j};((i,j),(i,j_2),\dotsc,(i,j_m));u].
\end{align*}
Due to \eqref{E:FACstrong2} we see that $(W^*,S^*)$ is an initial
pair. It is clear from the definitions that $$\phi(W^*,S^*) = (W',
[x_1\,t_{ij};((a_2,b_2),\dotsc,(a_m,b_m),(i,j));x])$$ which is a
Case A step followed by $m-1$ Case C steps. By induction we have
$\psi \circ \phi(W^*,S^*) = (W^*,S^*)$. Thus we have reduced to the
case that $m=1$.

In the case that $\psi$ on $(W',S')\in\OOL_{u,v}$ consists of an RA
step followed by $m>1$ RC steps, we may apply a similar reduction to
the $m=1$ case.

\subsection{$m=1$}

So we assume $m=1$. For $\psi\circ\phi$, we start with the initial
pair $(W,S)$ where
\begin{align*}
W&=(\ws{w^-}{A}{y}) \\
S&=(\xymatrix{{w^-}\ar[r]^{i,j} &{w}\ar[r]^{i,j'} & {u}}).
\end{align*}
We give diagrams before and after the Case C step.
\begin{align*}
\xymatrix{%
{w^-} \ar[r]^{i,j} \ar[d]_A & {w} \ar[r]^{i,j'} \ar[d]_A & {u}\ar[dl]^{\Ah} \\
{y} \ar[r]_{i,j} & {v} & {}}\qquad
\xymatrix{%
{}  & {w'} \ar[r]^{i,j} \ar[d]^{A'} \ar[dl]_{\Ah} & {u} \ar[d]^{A'} \\
{y} \ar[r]_{a,b} & {v'} \ar[r]_{i,j} & {x}}
\end{align*}
Here $\Ah$ is defined by \eqref{E:Ahdef} and $A'=\Ah\cup\{\ba{q}\}$
where $(a,b)=(y^{-1}(q),y^{-1}(p))$ and $q<p$ is the consecutive
pair of $\Ah$-bad integers such that $q<u(i)-1$ and $u^{-1}(q)\le l$
with $q$ maximum. We have $S'=(\xymatrix{{y}\ar[r]^{a,b} & {v'}
\ar[r]^{i,j} \ar[r] &{x}})$ and $W'=(\ws{u}{A'}{x})$. By Property
\ref{pro:markmove}(iv), the final pair $(\wcf{u}{A'}{x}{v'}{i,j})$
commutes. We now apply $\psi$. The first step is the commuting step
RA, which produces the commuting initial pair
$(\wci{w'}{A'}{v'}{i,j}{u})$. We now apply the next reverse
insertion step. By the Case C construction, the final pair
$(\wcf{w'}{A'}{v'}{y}{a,b})$ is noncommuting. In particular by Lemma
\ref{L:Rcyclicobstruction}, $w'(b)$ is not $A'$-nice and
\begin{align} \label{E:w'bq}
 w'(b)-1= c_{A'} w'(b)=v'(b)=y(a)=q \, .
\end{align}
In this situation
the reverse insertion algorithm defines
$$\Ah=A'-\{\overline{w'(b)-1}\}=A'-\{\ba{q}\}$$ which agrees with $\Ah$ as
defined above.

Condition C holds since $w'(j)=u(i)$ is $\Ah$-nice by definition. If
Case C holds, then it is easy to verify bijectivity: we have
$q_C=w'(j)=u(i)$, the previous $\Ah$-nice integer is indeed $p_C =
u(j')$ by Lemma \ref{L:Fcyclicobstruction} for the initial pair
$(\wci{w}{A}{v}{i,j'}{u})$ for the forward direction, and the
definition of $A$ for the current RC step agrees with the above
definition of $A$, as both are obtained from $\Ah$ by adding the
element $q_C-1=u(i)-1$.

So it suffices to show that Case C holds, that is, there is no
$\Ah$-nice integer $q_B$ such that $w'(b)<q_B<w'(j)$ and $g :=
w'^{-1}(q_B)\le l$. Suppose such an integer exists. We must derive a
contradiction.

We have
\begin{align} \label{E:yj'ui}
  p_0 = u(i) - 1 = w'(j) - 1 \, .
\end{align}
We claim that $q'=y(g)$ contradicts the
definition of $q$ in Case C. By definition $y^{-1}(q')=g\le l$. We
have that $y(g)=c_{\Ah}(w'(g))=c_{\Ah}(q_B)$. Since $q < q_B-1 <
u(i)-1=p_0$ are all $\Ah$-bad integers, it follows that
$q_B-1<y(g)$ are consecutive $\Ah$-bad integers and that $y(g)\le
p_0$. Furthermore $y(g)<p_0=u(i)-1$ since $g\le l < y^{-1}(p_0)$ by \eqref{E:p0pos}.  It only
remains to show that $y(g)>q$. We have $\ba{y(g)}\not\in
\{\ba{q},\ba{p}\}$, by Lemma~\ref{L:casecproperties} and the fact
that $g \le l < y^{-1}(p)$. Therefore
$$y(g)  = t_{qp} y(g) = t_{qp} c_{\Ah}(q_B) = c_{A'} (q_B) \ge q+1 >
q$$ and we arrive at the desired contradiction.

For $\phi\circ\psi$, the argument is nearly the same. By the
definition of Case RC, it is clear that in the subsequent
calculation of $\phi$, Case C will be invoked at the second step.
The only thing which needs to be checked is the maximality of $q$ in
Case C, which follows from the fact that Case B does not hold in the
second step of $\psi$.


\chapter{Grassmannian Elements, Cores, and Bounded Partitions}

\label{ch:coresbijection}

Let $k=n-1$ from now on. The $k$-Schur functions, denoted
$s_u^{(k)}(x)$ for $u\in\tS^0$ in Theorem \ref{T:kSchurStrong}, are
traditionally written $s_\la^{(k)}(x)$ where $\la$ is a partition
such that $\la_1 \le k$. Weak tableaux for Grassmannian elements
were first introduced~\cite{LMcore} as $k$-tableaux, which are
defined in terms of $(k+1)$-cores. In this chapter we recall
bijections between Grassmannian elements, offset sequences, cores,
and bounded partitions.

Let $l=0$ in this chapter.  The set $\tS^0 = \tS^l$ of Grassmannian
elements is defined in Section \ref{SS:maxparabolic}.

\section{Translation elements} \label{SS:trans}

Let $Q$ be the coroot lattice of $\mathfrak{sl}_n$, realized as the
set of $n$-tuples of integers with sum zero. $S_n$ acts on $Q$ by
permuting coordinates. Given $\beta=(\beta_1,\dotsc,\beta_n)\in Q$,
the translation element $\tau_\beta\in \tS$ is uniquely defined by
$\tau_\beta(i)=i+n\beta_i$ for $1\le i\le n$. We have $\tau_\beta
\tau_\gamma = \tau_{\beta+\gamma}$, so that
$T(Q)=\{\tau_\beta\mid\beta\in Q\}$ forms an abelian subgroup of
$\tS$ isomorphic to $Q$. By \eqref{E:Shi} we have
\begin{align} \label{E:translen}
  \ell(\tau_\beta) = \sum_{1\le i<j\le n} | \beta_j-\beta_i|
\end{align}
from which it follows that
\begin{align}
\label{E:transperm} %
\ell(\tau_\beta) &=\ell(\tau_{u\beta})
&\qquad&\text{for $\beta\in Q$ and $u\in S_n$} \\
\label{E:transantidom} %
\ell(\tau_\beta) &= 2 \sum_{i=1}^n i\beta_i= 2
\inner{\rho'}{\beta} &\qquad&\text{for $\beta$ antidominant,}
\end{align}
where $\rho'=(1,2,\dotsc,n)$ and an antidominant element is one that
is weakly increasing.  The inner product $\inner{.}{.}$ is the
standard one on ${\mathbb R}^n$.

$T(Q)$ acts on $Q$ by translations:
$\tau_\beta(\gamma)=\beta+\gamma$ for all $\gamma\in Q$. $S_n$ acts
on $T(Q)$ by conjugation: $s_i \tau_\beta s_i = \tau_{s_i\beta}$ for
$1\le i\le n-1$ and $\beta\in Q$. There is a well-known isomorphism
(\cite[Prop. 6.5]{Kac}) $\tS \cong S_n \ltimes T(Q)$ under which
$s_0\mapsto s_\theta \tau_{-\theta}$ where
$\theta=(1,0,\dotsc,0,-1)\in Q$ is the highest coroot and $s_\theta$
is the associated reflection, which satisfies
$s_\theta=t_{1,n}=s_1s_2\dotsc s_{n-2}s_{n-1}s_{n-2}\dotsm s_2 s_1$
and acts on $\Z^n$ by
$s_\theta(a_1,\dotsc,a_n)=(a_n,a_2,\dotsc,a_{n-1},a_1)$.

Again using \eqref{E:Shi}, for $v\in S_n$ and $\gamma\in Q$ we have
\begin{equation}\label{E:factlength}
 \ell(v\tau_\gamma) = \ell(\tau_\gamma) + \sum_{1\le i<j\le n} \chi(v(i)>v(j))
 (-1)^{\chi(\gamma_i<\gamma_j)}
\end{equation}
where $\chi(S)=1$ if $S$ is true and $\chi(S)=0$ if $S$ is false.

\begin{example} \label{X:semidirect}
Let $n=4$ and
$$
w=s_1s_2s_3s_0s_3s_2s_1s_0s_3s_2s_0s_3s_1s_0 = [-7,-1,4,14].
$$
Then $ws_3s_2=[-7,14,-1,4]$. In $\Z^4$ we have the equality
$$(-7,14,-1,4)=(1,2,3,4)+4(-2,3,-1,0),$$ so $ws_3s_2 = t_{-2,3,1,0}$
and $w=t_{-2,3,-1,0}s_2 s_3 = s_2 s_3 t_{-2,-1,0,3}$. By
\eqref{E:factlength} we have $\ell(w)=-2+\ell(t_{-2,-1,0,3})=
-2+2(1(-2)+2(-1)+3(0)+4(3))=14$, matching the length of the above
reduced word.
\end{example}

For $\beta \in Q$ let $S_n^\beta \subset S_n$ denote its
stabilizer.
\begin{prop} \label{P:replength}
Let $u \in S_n$.
\begin{enumerate}
\item \label{I:one}
If $\beta\in Q$ is antidominant and $u$ is of minimum length in its
coset $u S_n^\beta$ then
\begin{equation} \label{E:replength}
\ell(u\tau_\beta)=\ell(\tau_\beta)-\ell(u).
\end{equation}
\item \label{I:two} Let $u$ and $\beta$ be as in (\ref{I:one}).  If
$\gamma\in Q$ and $v\in S_n$ are such that $v\gamma=u\beta$ then
$\ell(v\tau_\gamma) \ge \ell(u\tau_\beta)$ with equality if and
only if $v=u$ and $\gamma=\beta$.
\end{enumerate}
\end{prop}
\begin{proof} Equation \eqref{E:factlength} says that
$\ell(v\tau_\gamma)-\ell(\tau_\gamma)$ is the number of inversions
$(i,j)$ of $v$ such that $\gamma_i\ge\gamma_j$, minus the number
of inversions of $v$ such that $\gamma_i<\gamma_j$. For $u$ and
$\beta$ in the hypotheses, we have $\beta_i\le\beta_j$ and
$\beta_j=\beta_i$ implies $u(i)<u(j)$.  This is precisely the
condition that the first set of inversions is empty and the second
is the set of all inversions of $u$. This proves \eqref{I:one}.

For \eqref{I:two} we must show that for $(v,\gamma)\in S_n \times
(S_n\cdot\beta)$ such that $v\gamma=u\beta$, $\ell(v\tau_\gamma)$ is
uniquely minimized by the pair $(u,\beta)$.

First let $\gamma\in S_n\cdot\beta$ be fixed. Suppose for some $1\le
p<q\le n$, $\gamma_p=\gamma_q$ and $v(p)<v(q)$. Let $v'=v t_{p,q}$.
We shall show that $\ell(v\tau_\gamma) < \ell(v'\tau_\gamma)$, which
implies that the desired minimum element $(v,\gamma)$ must have the
property that $v$ is of minimum length in the coset $v S_n^\gamma$.
For $1\le i<j\le n$ let
\begin{equation}
f(i,j) =
(-1)^{\chi(\gamma_i<\gamma_j)}(\chi(v'(i)>v'(j))-\chi(v(i)>v(j)))
\end{equation}
so that using \eqref{E:factlength}
\begin{equation}
\ell(v'\tau_\gamma)-\ell(v\tau_\gamma)=\sum_{1\le i<j\le n} f(i,j).
\end{equation}
For $(i,j)$ such that $\{i,j\}\cap\{p,q\}=\varnothing$, $f(i,j)=0$.
For $(i,j)$ such that $i=p$ and $j\not=q$ we have $p<j$ and
\begin{equation} \label{E:ip}
\begin{split}
  f(p,j) &= (-1)^{\chi(\gamma_p<\gamma_j)}(\chi(v(q)>v(j))-\chi(v(p)>v(j)))
  \\
  &= (-1)^{\chi(\gamma_p<\gamma_j)} \chi(v(p)<v(j)<v(q))
\end{split}
\end{equation}
For $(i,j)$ such that $i=q$ we have $q<j$ and
\begin{equation} \label{E:iq}
\begin{split}
  f(q,j) &= (-1)^{\chi(\gamma_q<\gamma_j)} (\chi(v(p)>v(j))-\chi(v(q)>v(j)))
  \\
  &= -(-1)^{\chi(\gamma_p<\gamma_j)} \chi(v(p)<v(j)<v(q))
\end{split}
\end{equation}
using $\gamma_q=\gamma_p$. We have
\begin{equation} \label{E:1}
\begin{split}
&\qquad\!\!\sum_{\substack{j>p\\ j\not=q}} f(p,j) + \sum_{j>q}
f(q,j) =
\sum_{p<j<q} f(p,j) \\
&= \sum_{p<j<q} (-1)^{\chi(\gamma_p<\gamma_j)} \chi(v(p)<v(j)<v(q)).
\end{split}
\end{equation}
Similarly
\begin{equation} \label{E:2}
\sum_{i<p} f(i,p)+\sum_{\substack{i<q\\i\not=p}} f(i,q)=
\sum_{p<i<q} (-1)^{\chi(\gamma_i<\gamma_p)} \chi(v(p)<v(i)<v(q)).
\end{equation}
The remaining term is $(i,j)=(p,q)$:
\begin{equation}\label{E:3}
f(p,q) =
(-1)^{\chi(\gamma_p<\gamma_q)}(\chi(v(q)>v(p))-\chi(v(p)>v(q))) \\
=1.
\end{equation}
Combining \eqref{E:1}, \eqref{E:2}, and \eqref{E:3} we have
\begin{equation}
\ell(v'\tau_\gamma)-\ell(v\tau_\gamma)= 1+2 \sum_{p<i<q}
\chi(\gamma_p=\gamma_i) \chi(v(p)<v(i)<v(q)) \ge 1.
\end{equation}
We may therefore assume that $v$ is of minimum length in its coset
$v S_n^\gamma$. Suppose next that $\gamma\not=\beta$, so that there
is an index $r$ such that $\gamma_r>\gamma_{r+1}$. Let
$\gamma'=s_r\gamma$ and $v'=vs_r$ so that $v'\gamma'=v\gamma$. It
suffices to show that $\ell(v'\tau_{\gamma'})<\ell(v\tau_\gamma)$.
Let
\begin{equation}
g(i,j)= (-1)^{\chi(\gamma'_i<\gamma'_j)}\chi(v'(i)>v'(j))
-(-1)^{\chi(\gamma_i<\gamma_j)}\chi(v(i)>v(j))
\end{equation}
Due to \eqref{E:transperm} and \eqref{E:factlength} we have
\begin{equation}
\ell(v'\tau_{\gamma'})-\ell(v\tau_\gamma)=\sum_{1\le i<j\le n}
g(i,j).
\end{equation}
If $\{i,j\}\cap\{r,r+1\}=\varnothing$ then $g(i,j)=0$. For $i=r$ and
$j>r+1$ we have
\begin{equation}
  g(r,j) = (-1)^{\chi(\gamma_{r+1}<\gamma_j)}\chi(v(r+1)>v(j))
-(-1)^{\chi(\gamma_r<\gamma_j)}\chi(v(r)>v(j)))
\end{equation}
For $i=r+1$ and $j>r+1$ we have
\begin{equation}
g(r+1,j)=(-1)^{\chi(\gamma_r<\gamma_j)}\chi(v(r)>v(j))
-(-1)^{\chi(\gamma_{r+1}<\gamma_j)}\chi(v(r+1)>v(j)))
\end{equation}
These cancel: for $j>r+1$ we have $g(r,j)+g(r+1,j)=0$. Similarly for
$i<r$ we have $g(i,r)+g(i,r+1)=0$. For $(i,j)=(r,r+1)$ we have
\begin{equation}
\begin{split}
g(r,r+1) &= (-1)^{\chi(\gamma_{r+1}<\gamma_r)}\chi(v(r+1)>v(r)) \\
&-(-1)^{\chi(\gamma_r<\gamma_{r+1})}\chi(v(r)>v(r+1)) \\
&=-\chi(v(r+1)>v(r))-\chi(v(r)>v(r+1)) \\
&= -1
\end{split}
\end{equation}
so that $\ell(v'\tau_{\gamma'})=\ell(v\tau_\gamma)-1$, which
suffices.
\end{proof}

\section{The action of $\tS$ on partitions} \label{sec:action}
Consider the positive quadrant $\Z_{>0}^2$ in the plane, where an
element $(i,j)$ is depicted as a \textit{cell} (square) in the
plane, indexed using standard Cartesian coordinates. The
\textit{diagram} of the partition $\la=(\la_1,\la_2,\dotsc,\la_p)$
is the set of cells $\{(i,j)\mid 1\le i\le p, 1\le j\le \la_i\}$
with $\la_i$ left-justified cells in row $i$ for all $i$. Define the
\textit{diagonal index} of a cell by $\diag(i,j)=j-i$ and the
\textit{residue} of a cell by $\res(i,j)=\ba{j-i} \in\Z/n\Z$.

Given a partition $\la$, one may associate a bi-infinite binary word
$p(\la)=p=\dotsm p_{-1} p_0 p_1 \dotsm $ called its \textit{edge
sequence}. The edge sequence $p(\la)$ traces the border of the
diagram of $\la$, going from northwest to southeast, such that every
letter $0$ (resp. $1$) represents a south (resp. east) step, such
that some cell in the $i$-th diagonal is touched by the steps
$p_{i-1}$ and $p_i$.

\begin{example} \label{X:edgeseq}
The diagram of $\la=(10,7,4,3,2,1,1,1)$ is pictured below.
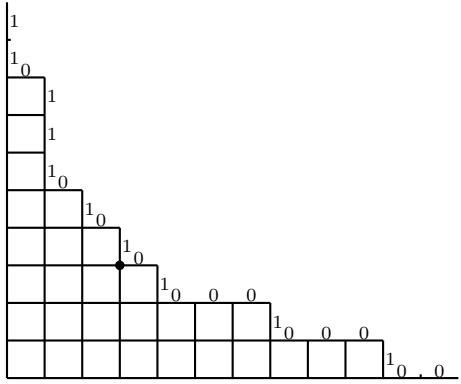
\begin{align*}
\psset{xunit=5mm,yunit=5mm}%
\pspicture(0,0)(12.5,10.5)%
\psdots(3,3)%
\psline(0,0)(0,10)%
\psline(1,0)(1,8)%
\psline(2,0)(2,5)%
\psline(3,0)(3,4)%
\psline(4,0)(4,3)%
\psline(5,0)(5,2)%
\psline(6,0)(6,2)%
\psline(7,0)(7,2)%
\psline(8,0)(8,1)%
\psline(9,0)(9,1)%
\psline(10,0)(10,1)%
\psline(0,0)(12,0)%
\psline(0,1)(10,1)%
\psline(0,2)(7,2)%
\psline(0,3)(4,3)%
\psline(0,4)(3,4)%
\psline(0,5)(2,5)%
\psline(0,6)(1,6)%
\psline(0,7)(1,7)%
\psline(0,8)(1,8)%
\psline(0,9)(.1,9)
\psline(11,0)(11,.1)
\rput(.2,9.5){\SB}%
\rput(.2,8.5){\SB}%
\rput(1.2,7.5){\SB}%
\rput(1.2,6.5){\SB}%
\rput(1.2,5.5){\SB}%
\rput(2.2,4.5){\SB}%
\rput(3.2,3.5){\SB}%
\rput(4.2,2.5){\SB}%
\rput(7.2,1.5){\SB}%
\rput(10.2,.5){\SB}%
\rput(.5,8.2){\SA}%
\rput(1.5,5.2){\SA}%
\rput(2.5,4.2){\SA}%
\rput(3.5,3.2){\SA}%
\rput(4.5,2.2){\SA}%
\rput(5.5,2.2){\SA}%
\rput(6.5,2.2){\SA}%
\rput(7.5,1.2){\SA}%
\rput(8.5,1.2){\SA}%
\rput(9.5,1.2){\SA}%
\rput(10.5,.2){\SA}%
\rput(11.5,.2){\SA}%
\endpspicture
\end{align*}
The edge sequence is $p(\la)=\dotsm
11|0111|0101\bullet0100|0100|0100|\dotsm$ where $\bullet$ indicates
the $0$ diagonal, which separates the bits $p_{-1}$ and $p_0$.
\end{example}

The affine symmetric group $\tS$ acts on partitions in an obvious
way, if we identify elements of $\tS$ with functions
$\Z\rightarrow\Z$ and partitions with their edge sequences, which
are certain functions $\Z\rightarrow\{0,1\}$. Say that the cell $x$
is \textit{$\la$-addable} (resp. \textit{$\la$-removable}) if adding
(resp. removing) the cell $x$ to (resp. from) the diagram of $\la$
results in the diagram of a partition. Then for $i\in\Z/n\Z$,
$s_i\la$ is obtained by removing from $\la$ every $\la$-removable
cell of residue $i$, and adding to $\la$ every $\la$-addable cell of
residue $i$.

\begin{example} \label{X:affpermonpart}
For $n=4$, applying the reflections of the reduced word of $w$ in
Example \ref{X:semidirect} to the empty partition from right to
left, we obtain the sequence of partitions $()\subset(1)\subset (2)
\subset (2,1) \subset (2,2) \subset (3,2,1) \subset (4,2,2) \subset
(5,2,2) \subset (6,3,2,1) \subset (7,4,2,2) \subset (8,5,2,2)
\subset (9,6,3,2,1) \subset (9,6,3,3,1,1)\subset (9,6,3,3,1,1,1)
\subset (10,7,4,3,2,1,1,1)$.
\end{example}

In fact, the action of $\tS$ on partitions coincides with the action
of the Kashiwara reflection operators on the crystal graph of Fock
space \cite{MM}.

\section{Cores and the coroot lattice} An \textit{$n$-ribbon} is a
skew partition diagram $\la/\mu$ (the difference of the diagrams of
the partitions $\la$ and $\mu$) consisting of $n$ rookwise connected
cells, all with distinct residues. We say that this ribbon is
\textit{$\la$-removable} and \textit{$\mu$-addable}. An
\textit{$n$-core} is a partition that admits no removable
$n$-ribbon. Henceforth when we say ``core" we mean ``$n$-core".
Since the removal of an $n$-ribbon is the same thing as exchanging
bits $p_i=0$ and $p_{i+n}=1$ in the edge sequence for some $i$, it
follows that $\la$ is a core if and only if for every $i$, the
sequence $p^{(i)}(\la) := \dotsm p_{i-2n} p_{i-n} p_i p_{i+n}
p_{i+2n}\dotsm$ consisting of the subsequence of bits indexed by $i$
mod $n$, has the form $\dotsm 1111100000\dotsm$. Thus the core is
specified by the positions where this sequence changes from $1$ to
$0$ for various $i$. For a core $\la$ and $i\in\Z$ let
$d_i(\la)=d_i\in\Z$ be the integer such that $p_{i-1+n (d_i-1)}=1$
and $p_{i-1+n d_i}=0$.\footnote{Using $i-1$ rather than $i$  allows
an elegant statement of Proposition \ref{P:scover}.} The sequence
$(d_i)_{i\in\Z}$ is called the \textit{extended offset sequence} of
$\la$, and $d(\la)=(d_1,d_2,\dotsc,d_n)$ the \textit{offset
sequence} of $\la$. Observe that
\begin{align} \label{E:offsetshift}
  d_{i-n} = d_i + 1\qquad\text{for all $i\in\Z$}
\end{align}
and that $\sum_{i =1}^n d_i = 0$.

\begin{example} \label{X:coreoffset} For $n=4$ and $\la=(10,7,4,3,2,1,1,1)$,
we have $d(\la)=(-2, 3, -1, 0)$. These are the heights ``above sea
level" that the 1s attain if we draw the edge sequence $p(\la)$ in
an $\infty \times n$ array with the $r$-th row given by the binary
word with bits $p_{i+nr}$ for $i=0,1,2,\dotsc,n-1$, and ``sea level"
is the division between rows $-1$ and $0$. Below we depict the bit
sequence $p(\la)$, whose columns, read from bottom to top, are
$p^{(0)}(\la)$ through $p^{(n-1)}(\la)$.
\begin{align*}
\begin{matrix}
\vdots&\vdots&\vdots&\vdots\\
0&0&0&0\\
0&1&0&0\\
0&1&0&0\\
0&1&0&0\\ \hline
0&1&0&1\\
0&1&1&1\\
1&1&1&1 \\
\vdots&\vdots&\vdots&\vdots
\end{matrix}
\end{align*}
\end{example}

The following results are well known; see for example~\cite{L, vL}.

\begin{lem} \label{L:coreres} For any core $\la$ and $i\in \Z/n\Z$,
either there are no $\la$-addable cells of residue $i$ or there are
no $\la$-removable cells of residue $i$.
\end{lem}

\begin{prop} \label{P:grasscore} There is a bijection from the set
$\Core$ of $n$-cores to the root lattice $Q$ of $\mathfrak{sl}_n$
given by $\la\mapsto d(\la)$.
\end{prop}

$Q$ already has an action of $\tS$ defined in Section
\ref{SS:trans}, which is transitive since the subgroup $T(Q)$ of
$\tS$ acts transitively on $Q$. Via the bijection $\la\mapsto
d(\la)$, $\tS$ acts transitively on the set $\Core$ of cores. This
induced action coincides with the action of $\tS$ on partitions
defined above.

\begin{prop} \label{P:affpermaction} The action of $\tS$ on
partitions restricts to an action on $\Core$. Moreover, the
bijection of Proposition \ref{P:grasscore} is an isomorphism of sets
with $\tS$-action:
\begin{align} \label{E:affpermaction}
d(w\cdot \la)=w \cdot d(\la) \qquad\text{for every $w\in\tS$.}
\end{align}
\end{prop}
\begin{proof} It suffices to consider the case
$w=s_i$. This is easily verified, the most interesting case being
$i=0$, where one uses \eqref{E:offsetshift}.
\end{proof}

\begin{prop} \label{P:coreorbit} $\Core=\tS\cdot \vn$.
In particular there is a bijection
\begin{align} \label{E:coremap}
  \core: \tS^0 &\rightarrow \tS/S_n \rightarrow \Core \\
  w&\mapsto w S_n \mapsto w \cdot \vn.
\end{align}
\end{prop}
\begin{proof} We already know that $\tS$ acts transitively on
$\Core$ and $\vn\in\Core$, so that $\Core=\tS\cdot\vn$. Since $S_n$
is the stabilizer of $\vn$ the result follows.
\end{proof}

\section{Grassmannian elements and the coroot lattice}

The bijection $d \circ\core:\tS^0 \rightarrow  Q$ has the following
affine Lie-theoretic interpretation.

\begin{prop} \label{P:semidirect} Write $w\in
\tS^0$ as $w=u\tau_\beta$ where $u\in S_n$ and $\beta\in Q$.
\begin{enumerate}
\item $d(\core(w))= u \beta\in Q$. \item $\beta$ is antidominant
and $u$ is of minimum length in $u S_n^\beta$ where $S_n^\beta$ is
defined before Proposition~\ref{P:replength}. \item
$\ell(w)=\ell(t_\beta)-\ell(u)$.
\end{enumerate}
\end{prop}
\begin{proof} Let $\la=\core(w)$ (that is, $\la=w\cdot \vn$)
and $d=d(\la)$. We have $d(\core(w))=d=d(w\cdot \vn)=w\cdot
d(\vn)=w\cdot (0^n) = u \tau_\beta (0^n)=u\beta$ using Proposition
\ref{P:affpermaction}. In fact this works for any $w\in\tS$ such
that $w\cdot \vn=\la$ and $w=u\tau_\beta$.

Suppose $w'=v \tau_\gamma$ is such that $w'\cdot \vn =\la$. Then
$v\gamma=d$ as well. We have
\begin{equation*}
\begin{split}
  w' = v \tau_\gamma &= v \tau_{v^{-1}u\beta} \\
  &=   v v^{-1} u \tau_\beta u^{-1} v \\
  &= u\tau_\beta (u^{-1}v) = w (u^{-1} v).
\end{split}
\end{equation*}
The result follows from Proposition \ref{P:replength}.
\end{proof}

\begin{example} \label{X:dc} For $n=4$ and $w=[-7,-1,4,14]$, by
Example \ref{X:semidirect} we have $u=s_2s_3$ and
$\beta=(-2,-1,0,3)$. So $d=u\beta=(-2,3,-1,0)$, which agrees with
Example \ref{X:coreoffset}.
\end{example}

\section{Bijection from cores to bounded partitions} Say that a
partition $\la$ is $n$-bounded\footnote{Our language differs
slightly from that of \cite{LMcore}.} if $\la_1<n$. Denote by
$\Bounded$ the set of $n$-bounded partitions.

The \textit{hook length} of a cell $(i,j)$ in a skew shape $\la/\mu$
is the number of cells in $\la/\mu$ in the $i$-th row to the right
of $(i,j)$, plus the number of cells in $\la/\mu$ in the $j$-th
column above $(i,j)$, plus one (for the cell $(i,j)$ itself).

\begin{prop} \label{P:corebounded} \cite{LMcore}
Let $\la\in \Core$. Let $\mu\subset\la$ be the smallest partition
(in the order by containment) such that the hook length of every
cell in $\la/\mu$ is less than $n$. Then $\la\mapsto
\bounded(\la):=(\la_1-\mu_1,\la_2-\mu_2,\dotsc)$ defines a bijection
$\Core\rightarrow\Bounded$.
\end{prop}

\begin{example} \label{X:corebounded} For
$n=4$ and $\la$ as in Example \ref{X:coreoffset}, we have
$\mu=(7,4,2,1,1)$ so that $\la\mapsto (3,3,2,2,1,1,1,1)$, which is
obtained by reading the sizes of the rows of the skew diagram
$\la/\mu$, which is pictured below.
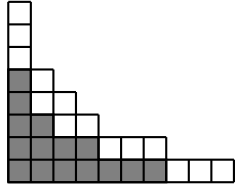
\begin{align*}
\psset{xunit=3mm,yunit=3mm}%
\pspicture(-.5,-.5)(10.5,8.5)%
\pspolygon[fillstyle=solid,fillcolor=gray](0,0)(7,0)(7,1)(4,1)(4,2)(2,2)(2,3)(1,3)(1,5)(0,5)(0,0)%
\psline(0,0)(10,0)%
\psline(0,1)(10,1)%
\psline(0,2)(7,2)%
\psline(0,3)(4,3)%
\psline(0,4)(3,4)%
\psline(0,5)(2,5)%
\psline(0,6)(1,6)%
\psline(0,7)(1,7)%
\psline(0,8)(1,8)%
\psline(0,0)(0,8)%
\psline(1,0)(1,8)%
\psline(2,0)(2,5)%
\psline(3,0)(3,4)%
\psline(4,0)(4,3)%
\psline(5,0)(5,2)%
\psline(6,0)(6,2)%
\psline(7,0)(7,2)%
\psline(8,0)(8,1)%
\psline(9,0)(9,1)%
\psline(10,0)(10,1)%
\endpspicture
\end{align*}
\end{example}

\section{$k$-conjugate} Transposition defines an involution
$\omega$ on the set of partitions, which is easily seen to restrict
to an involution on the set $\Core$ of $n$-cores. Denote by
$\omega^{(n)} = \bounded \circ \omega\circ \bounded^{-1}$ the
induced involution on $\Bounded$. For $k=n-1$, $\omega^{(n)}$ is the
$k$-conjugate map of \cite{LMcore}.

\begin{example} \label{X:ktranspose} Reading the column sizes of the
skew shape in Example \ref{X:corebounded} we have
$\omega^{(4)}(3,3,2,2,1,1,1,1)=(3,2,2,1,1,1,1,1,1,1)$.
\end{example}

\section{From Grassmannian elements to bounded partitions}

Recall from Section \ref{SS:coxeter} the definition of an inversion
of $w\in\tS$. Let the \textit{code} of $w$ be the sequence
$\code(w)=(c_1,c_2,\dotsc,c_n)$ where $c_i$ is the number of
inversions of $w$ of the form $(i,j)$ for some $j$; see
also~\cite{BB} where the code is called the inversion table. For a
Grassmannian permutation it is easy to see that the code is a weakly
increasing sequence of nonnegative integers that starts with $0$;
reversing this sequence yields a partition with fewer than $n$
parts. Applying the transpose $\omega$, we obtain an $(n-1)$-bounded
partition.

\begin{prop} \label{P:shapetrans} The composite bijection $\bounded \circ
\core:\tS^0\rightarrow \Bounded$ is given by $w\mapsto
\omega^{(n)}(\omega(c_n,\dotsc,c_2,c_1))$ where
$\code(w)=(c_1,c_2,\dotsc,c_n)$.  This bijection satisfies
$\ell(w) = |\bounded \circ \core(w)|$.
\end{prop}

\begin{proof}
Bj\"{o}rner and Brenti~\cite[Theorem 4.4]{BB} recursively
construct a bijection from bounded partitions to Grassmannian
affine permutations $w$, from which one obtains the code
$\code(w)$. This recursive construction is compatible with the
construction of a core from a Grassmannian permutation as in
Section~\ref{sec:action}.
\end{proof}

\begin{example} Let $n=4$ and $w=[-7,-1,4,14]$. Then $w$ has
inversions $(2,5)$, $(3,5)$, $(3,6)$, $(3,9)$, $(4,5)$, $(4,6)$,
$(4,7)$, $(4,9)$, $(4,10)$, $(4,11)$, $(4,13)$, $(4,14)$, $(4,17)$,
$(4,21)$ so that $\code(w)=(0,1,3,10)$. Now
$\omega(\code(w))=(3,2,2,1,1,1,1,1,1,1)$ and by Example
\ref{X:ktranspose} we have
$\omega^{(4)}(\omega(\code(w))=(3,3,2,2,1,1,1,1)$, which agrees with
Example \ref{X:corebounded}.
\end{example}

\chapter{Strong and Weak Tableaux Using Cores}
\label{ch:corestableaux}

We now specialize all constructions to the special case of
Grassmannian elements. By the previous chapter we may work instead
with cores.

\section{Weak tableaux on cores are $k$-tableaux} Consider a weak
tableau $U=(W_1,W_2,\dotsc)$ where $\ins(W_i)=u_{i-1}$ and
$\out(W_i)=u_i$ for all $i$, and $u_i\in \tS^0$ for all $i\ge0$. Let
$\mu=\core(\ins(U))$ and $\la=\core(\out(U))$. The weak tableau $U$
can be depicted as a usual tableau of shape $\la/\mu$ in which all
cells of the skew shape $\core(u_i)/\core(u_{i-1})$ have been filled
with the letter $i$ for all $i$. Such a tableau is called a {\it
$k$-tableau} \cite{LMcore}.  The following characterization of
$k$-tableaux is given by Lapointe, Morse and Wachs in~\cite{LMW}.

\begin{lem} \label{L:weakktab} $k$-tableaux are semistandard
(increasing in rows and strictly increasing in columns). Conversely,
if a chain of cores in (left) weak order defines a semistandard
tableau then it comes from a weak tableau.
\end{lem}


\begin{proof}
If $s_\sigma$ is a cyclically decreasing permutation, then $s_{j}$
is never applied after $s_{j+1}$ in  $s_{\sigma}(\gamma)$. This
means that two cells are never added on top of each other in the
same column and thus $s_{\sigma}(\gamma)/\gamma$ is a horizontal
strip.  This shows the first statement of the lemma.

Now suppose that $\delta = \core(w)$ and $\gamma = \core(v)$ where
$w \wle v\in \tS^0$ are such that $\delta/\gamma$ is a horizontal
strip of size $m$ (meaning that the difference $\ell(w)-\ell(v)$ is
$m$). Theorem 56 of \cite{LMcore} then says that $\delta = s_{i_1}
\cdots s_{i_m} (\gamma)$, for some $i_1,\dots,i_m$ that are all
distinct. Furthermore, in the proof of Theorem 56, it is shown that
$i_m$ can be chosen as the residue of the southeasternmost cell
$c_{SE}$ in $\delta/\gamma$. By induction on the size of
$\delta/\gamma$, it suffices to show that $s_{i_m+1}$ is never
applied at some point after $s_{i_m}$ so that $s_{i_1} \cdots
s_{i_m}$ is cyclically decreasing. Suppose this is the case. Then,
since $\delta/\gamma$ is a horizontal strip, this means that there
is a cell of residue $i_m+1$ to the northwest of $c_{SE}$ without a
cell above it, and thus that the extremal cell of residue $i_m$ to
its left is not at the end of its row. This is a contradiction to
\cite[Proposition 15(1)]{LMcore} which says that all extremal cells
of residue $i_m$ to the northwest of $c_{SE}$ are at the end of
their row.
\end{proof}

\begin{example} \label{X:ktab} Let $n=4$ and $l=0$ and consider the weak
tableau defined by the length-additive factorization of the
Grassmannian element
$$w=(s_1)(s_2)(s_3)(s_0)(s_3s_2s_1)(s_0s_3s_2)(s_0s_3)(s_1s_0)$$ of Example
\ref{X:semidirect} into cyclically decreasing elements. The
corresponding sequence of cores is $()\subset (2) \subset (2,2)
\subset (5,2,2) \subset (8,5,2,2) \subset (9,6,3,2,1) \subset
(9,6,3,3,1,1)\subset (9,6,3,3,1,1,1) \subset (10,7,4,3,2,1,1,1)$,
which gives the $k$-tableau \setcellsize{11}
\begin{align*}
\tableau{{8}\\{7}\\{6}\\{5}&{8}\\{4}&{4}&{6}\\{3}&{3}&{5}&{8}\\{2}&{2}&{4}&{4}&{4}&{5}&{8}\\
{1}&{1}&{3}&{3}&{3}&{4}&{4}&{4}&{5}&{8}}
\end{align*}
\end{example}

\section{Strong tableaux on cores} We assume that $l=0$ with $l$ as
in Section \ref{SS:maxparabolic}. Since $\id$ is Grassmannian, by
Proposition \ref{P:strongGrass} any strong tableau with inner shape
$\id$, involves only Grassmannian elements of $\tS$.

The strong order on Grassmannian permutations corresponds to
containment of cores.

\begin{prop}[\cite{L, MM}] \label{P:crystal} For $v,w\in \tS^0$,
$v\le w$ if and only if $\core(v)\subset\core(w)$. %
\end{prop}

\begin{lem} \label{L:strongoffset}
Let $\mu$ be a core and $t_{r,s}\in \tS$ a reflection with $r<s$.
Then
\begin{enumerate}
\item \label{I:offset1} $t_{r,s} \mu \ge \mu$ if and only if $d_r
\ge d_s$ using the extended offset sequence $d$ of $\mu$. \item
\label{I:offset2} $t_{r,s}\mu \scov \mu$ if and only if $d_r>d_s$
and for all $r<i<s$, $d_i\not\in [d_s,d_r]$.
\end{enumerate}
\end{lem}
\begin{proof} (2) follows from (1). To prove (1),
for cores $\mu$ and $\la$, by Proposition \ref{P:crystal}, $\mu\le
\la$ if and only if for all $i\in\Z$, $\sum_{j\le i}
(p_j(\mu)-p_j(\la)) \ge 0$. Therefore $t_{r,s} \mu \ge \mu$ if and
only if $p_{r+qn}(\mu)\ge p_{s+qn}(\mu)$ for all $q\in\Z$. $d_r\ge
d_s$ implies that this holds for $q=0$ and \eqref{E:offsetshift}
implies it for all other $q$.
\end{proof}

The \textit{head} of a ribbon is its southeastmost cell.

\begin{prop} \label{P:scover} Let $\mu\scovby \la$ be cores with
$t_{r,s} \mu=\la$ and $0<r<s$. Then
\begin{enumerate}
\item \label{I:lessn} $s-r<n$.
\item Each connected component of $\la/\mu$ is a ribbon with $s-r$
cells in diagonals of residue $\ba{r},\ba{r+1},\dotsc,\ba{s-1}$.
\item The components are translates of each other and their heads
lie on ``consecutive" diagonals of residue $s-1$.
\item \label{I:ncomp} The skew shape $\la/\mu$ has
$d_r-d_s$ components.
\end{enumerate}
\end{prop}

\begin{example} \label{X:offsetcover} Let $n=4$ and $l=0$ and
consider the strong marked cover $C=(\stc{w}{i,j}{u})$ where
$w=[-8,-3,6,15]$, $u=[-8,-6,9,15]$, and $(i,j)=(-1,10)$. We have
$m(C)=w(j)=5$.  The affine Grassmannian permutations $w$ and $u$
have associated cores $\mu=(11,8,5,5,3,3,1,1,1)$ and
$\la=(11,8,7,6,5,4,3,2,1)$ and offsets $d(\mu)=(-1,1,3,-3)$ and
$d(\la)=(2,-2,3,-3)$ respectively. Letting $(r,s)$ be such that
$t_{r,s} \mu = \la$, we have $t_{r,s} w = u = w t_{i,j}$, so that
$t_{r,s} = w t_{i,j} w^{-1} = t_{w(i),w(j)}$. Thus we may take
$r=w(i)=2$ and $s=w(j)=5$. Letting $d=d(\mu)$, the skew shape
$\la/\mu$ has $d_2-d_5=1-(-2)$ components, each of size $5-2$.
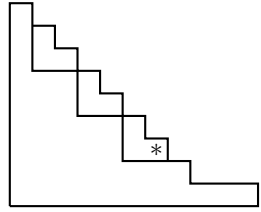
\begin{align*}
\psset{xunit=3mm,yunit=3mm}%
\pspicture(-.5,-.5)(11.5,9.5)%
\psline(0,0)(11,0)(11,1)(8,1)(8,2)(5,2)(5,4)(3,4)(3,6)(1,6)(1,9)(0,9)(0,0)%
\psline(7,2)(7,3)(6,3)(6,4)(5,4)(5,5)(4,5)(4,6)(3,6)(3,7)(2,7)(2,8)(1,8)%
\rput(6.5,2.5){$*$}
\endpspicture
\end{align*}
\end{example}

\begin{proof}[Proof of Proposition \ref{P:scover}]
All of the assertions follow from the first and Lemma
\ref{L:strongoffset}.  The first assertion follows from
Lemma~\ref{L:strongoffset}(\ref{I:offset2}) and
\eqref{E:offsetshift}.
\end{proof}

\begin{prop} \label{P:strongcovercore} Let $C=(\stc{w}{i,j}{u})$ be a
marked strong cover with $w,u\in \tS^0$ and write $\core(w)=\mu$ and
$\core(u)=\la$. Then the head of one of the ribbons forming the
connected components of $\la/\mu$, is on the diagonal $m(C)-1$.
\end{prop}
\begin{proof} We apply Proposition \ref{P:scover} with $r=w(i)$ and
$s=w(j)=m(C)$. Since the $r$-th bit $0$ is being exchanged with the
$s$-th bit $1$ in the strong cover, the head of one of the ribbons
must lie on the previous diagonal $s-1$.
\end{proof}

\begin{example} \label{X:offsetmark} In Example \ref{X:offsetcover}
the head of one of the ribbons in $\la/\mu$ is marked with a $*$; it
is in the diagonal $m(C)-1=5-1$.
\end{example}

In light of Proposition \ref{P:strongcovercore}, a marked strong
cover of cores $\stc{\mu}{i,j}{\la}$ is given by the disjoint union
of ribbons comprising the skew shape $\la/\mu$, together with a
marking on the head of one of the ribbons. For a strong strip, the
sequence of marked cells must have strictly increasing diagonal
indices. Since each strong cover on cores is a skew partition shape,
this is equivalent to saying that the marked heads must proceed
weakly to the south and strictly to the east.

We make the following convention for strong tableaux on cores. Let
$T=(S_1,S_2,\dotsc)$ be a strong tableau going between the elements
$u_0 \le u_1 \le \dotsc$ of $\tS^0$. Let $\la^{(i)}=\core(u_i)$,
$\mu=\core(\ins(T))$, and $\la=\core(\out(T))$. We depict $T$ as a
tableau of shape $\la/\mu$. The strong strip $S_i$ is depicted by
the skew subtableau of shape $\la^{(i)}/\la^{(i-1)}$. We distinguish
the cells of the strong covers in $S_i$ by placing the subscripted
letter $i_j$ in a cell if it occurs in the skew shape corresponding
to the $j$-th strong cover in $S_i$. Finally, a mark ${}^*$ must be
placed on the head of some ribbon in each strong cover, such that
the diagonals of the heads increase in each strong strip.

\begin{example} \label{X:strongtableaucore} Here is a strong tableau
for $n=3$. %
\setcellsize{13}%
$$P=\tableau{{3_1} \\ {2_1^*}&{3_3}&{3_3}
\\{1_1^*}&{3_1^*}&{3_2^*}&{3_3}&{3_3^*}}.$$
It has three nonempty strong strips $S_1,S_2,S_3$. $S_1$ consists of
a single marked strong cover $\stc{\vn}{0,1}{(1)}$. $S_2$ also
consists of a single marked strong cover $\stc{(1)}{-1,1}{(1,1)}$.
$S_3$ is given by $(1,1)\overset{0,4}{\rightarrow} (2,1,1)
\overset{0,2}{\rightarrow} (3,1,1) \overset{0,5}{\rightarrow}
(5,3,1)$. To understand the markings, we convert the $3$-cores in
$S_3$ into elements of $\tS$, giving the sequence
$[-1,3,4],[-2,3,5],[-2,2,6],[-2,0,8]$ in window notation.  Call
these $u_0,u_1,u_2,u_3$. Then the marks of the three covers in $S_3$
are given by $u_0(4)=2$, $u_1(2)=3$, and $u_2(5)=5$. This means that
the heads of the marked ribbons in $S_3$ occur in diagonals
$2-1,3-1,5-1$, which indeed is the case.
\end{example}

\begin{example}
We compute an example of Theorem~\ref{T:weakPieri} using cores.
Let $n = 3$, $r = 2$ and $w = s_2 s_0$.  Then we have the equality
(in $\Lambda^{(3)})$ $$ h_2(x) \WS_{s_2s_0}(x) =
2\WS_{s_2s_1s_2s_0}(x) + \WS_{s_0s_1s_2s_0}(x) +
\WS_{s_0s_2s_1s_0}(x).$$  The right hand side corresponds to the
following four strong strips on cores.

$$\tableau{{1_1^*} \\ {} \\ {}&{1_1} &{1_2^*}} \hspace{20pt} \tableau{{1_1} \\ {}\\ {}&{1_1^*}
&{1_2^*}}\hspace{20pt} \tableau{{1_2}\\{1_1^*} \\ {}&{1_2^*} \\
{}&{1_1} } \hspace{20pt} \tableau{{}&{1_2}
\\ {}&{1_1} &{1_1^*}&{1_2^*}}
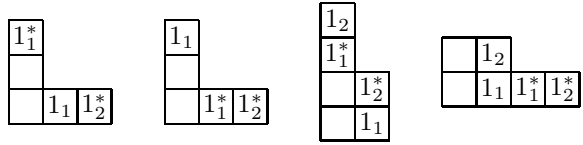
$$
\end{example}

\section{Monomial expansion of $t$-dependent $k$-Schur functions}

The $k$-Schur functions
were discovered by Lapointe, Lascoux and Morse  when studying
the Macdonald polynomials $H_\lambda(x;q,t)$.  To summarize,
great attention has been given to the study of
the {\it $q,t$-Kostka coefficients} in the Schur expansion,
\begin{equation}
H_\mu(x;q,t) = \sum_\lambda K_{\lambda\mu}(q,t) \,s_\lambda
\,.
\label{macdo}
\end{equation}
By way of geometry of Hilbert schemes, Haiman proved  \cite{[Ha]}
that the $K_{\lambda\mu}(q,t)$ lie in $\mathbb N[q,t]$ and are
associated to the characters of irreducible representations in a
bigraded representation of the symmetric group \cite{[GH]}. A
long-standing open problem is to find a tableaux characterization
for these expansion coefficients.  For one thing, it is known that
the number of monomial terms in $q$ and $t$ that occur in a given
$K_{\lambda\mu}(q,t)$ equals the number of standard tableaux of
shape $\lambda$.  Furthermore, in the case that $q=0$ (the
Hall-Littlewood case), there is a beautiful solution to this problem
given by the {\it charge} statistic on tableaux \cite{LaSc}.

The $k$-Schur functions arose in the study of this problem
when it was empirically observed  that there
exists a remarkable refinement of the expansion \eqref{macdo}.
To be more precise, computer evidence suggested the existence of
a family of polynomials \cite{LLM} defined by certain sets of
tableaux $\mathcal A^k_\mu$ as:
\begin{equation} \label{eqatom}
s_{\mu}^{(k)}(x;t) = \sum_{T\in \mathcal A^k_\mu}
t^{charge(T)} \, s_{shape(T)}
\end{equation}
with the property that
any Macdonald polynomial indexed
by a $k$-bounded partition $\lambda$ can be decomposed as:
\begin{equation}
H_{\lambda}(x;q,t) = \sum_{\mu;\mu_1\leq k}
K_{\mu \lambda}^{(k)}(q,t) \, s_{\mu}^{(k)}(x;t) \, , \qquad
K_{\mu \lambda}^{(k)}(q,t) \in \mathbb N[q,t] \, .
\label{(1.7)}
\end{equation}
Moreover, given that in this setting $s_\mu^{(k)}(x;t)=s_\mu$
for large $k$, such a decomposition reduces
to \eqref{macdo} when $k$ is large enough.

The original definition of $s_\mu^{(k)}(x;t)$ given
in \cite{LLM} is a combinatorial characterization of
the sets $\mathcal A_\mu$.  As it relies too heavily
on tableaux combinatorics to warrant full presentation
here, we will instead give a less cumbersome definition
of $\tilde s_{\mu}^{(k)}(x;t)$ from \cite{LMfil}
that is conjectured to equal $s_{\lambda}^{(k)}(x;t)$.

Formally, for a $k$-bounded partition
$\la=(\la_1,\ldots ,\la_\ell)$,
 the definition is simply
\begin{equation}
\tilde s_{\lambda}^{(k)}(x;t) \, =
{T}_{\lambda_1}^{(k)} B_{\lambda_1}^t \cdots
{T}_{\lambda_\ell}^{(k)} B_{\lambda_\ell}^t \cdot 1 \, .
\label{(1.9)}
\end{equation}
In this formula, the operators $B_a^t$ are the vertex operators
introduced in \cite{[Ji]} that recursively build
Hall-Littlewood polynomials:
$$
H_{\lambda}^{(k)}(x;0,t) =
\, \,    B_{\lambda_1}^t \cdots
B_{\lambda_\ell}^t \cdot 1 \, ,
$$
and the operators ${T}_{m}^{(k)}$, as we will see, encode the complexity of
$k$-Schur functions.  To define
${T}_{m}^{(k)}$, we require the  construction of an auxiliary basis
$\big\{G_\la(x;t)\big\}_{\la_1\leq k}$ for $\Lambda^k_t$ (the $\mathbb
Q$-linear span of
Hall-Littlewood polynomials indexed by $k$-bounded partitions). This
done, ${T}_{m}^{(k)}$
is defined by setting
\begin{equation}
{T}_m^{(k)} \, G_{\lambda}^{(k)}(x;t) =
\begin{cases}
G_{\lambda}^{(k)}(x;t) &   \text{if}\quad  \lambda_1=m \, ,\\
0 &   \text{otherwise }\,.
\cr
\end{cases}
\label{1.10}
\end{equation}

The basis $G_{\lambda}^{(k)}(x;t)$ is best understood by working
first under the specialization $t=1$.  Associate to any $k$-bounded
partition $\lambda$ a sequence of partitions,
$$\lambda^{\to k}= (\lambda^{(1)},\lambda^{(2)},\ldots,\lambda^{(r)})\,,$$
called the $k$-{\it split} of $\lambda\, $.
The sequence $\lambda^{\to k}$ is obtained by partitioning $\lambda$
(without rearranging the entries) into partitions $\lambda^{(i)}$
with main hook-length equal to   $k$, for all $i<r $.
The last partition in   $\lambda^{\to k}$ may have main
hook-length less than $k$.
For example,
$$(3,2,2,2,1,1)^{\to 3}=\bigl((3),(2,2),(2,1),(1)\bigr)$$
and
$$(3,2,2,2,1,1)^{\to 4}=\bigl((3,2),(2,2,1),(1)\bigr)\,.$$
It is important to note that $\lambda^{\to k}=(\lambda)$
when the main hook-length of $\la$ is $\leq k$. This given,
let $G_{\lambda}^{(k)}(x;1) $ be the ordinary
Schur function product
$$
G_{\lambda}^{(k)}(x;1)\,=\, s_{\la^{(1)}}s_{\la^{(2)}}\cdots s_{\la^{(r)}}\, .
\eqno (1.11)
$$
These functions are called ``{\it $k$-split}'' polynomials.

The $k$-split polynomials when $t\neq 1$ are a subfamily of
a family of
polynomials, $\mathcal H_S(x;t)$, known as generalized
Kostka polynomials (see for instance \cite{SW,SZ}).  These
polynomials, indexed by arbitrary sequences $S$ of partitions,
are $t$-analogs of products of Schur functions.  The
$k$-split polynomials correspond to the
generalized
Kostka polynomials for which
$S=\lambda^{\to k}\,.$
That is,
$$
G_{\lambda}^{(k)}(x;t):=\mathcal H_{\lambda^{\to k}}(x;t) \,.
$$

The primary goal of investigations carried out in \cite{LLM,LMfil}
was to find a fruitful characterization of $k$-Schur functions.
The empirical evidence that these functions
satisfy analogs of Schur function properties such as Pieri and
Littlewood-Richardson rules at $t=1$ motivated an independent study
of the $s_\mu^{(k)}(x) =s^{(k)}_\mu(x;1)$ case.  The intricacy of
the original definitions prompted a third (conjecturally equivalent)
characterization for the special case that $t=1$ in \cite{LMproofs}.
That is the definition in terms of weak $k$-tableaux that we have
used for $s_\mu^{(k)}(x)$ throughout this article.

The developments in this article have led us to prove in
Theorem~\ref{T:kSchurStrong} that these $s_\mu^{(k)}(x)$
can be defined by the weight generating function of
strong $k$-tableaux, providing a striking refinement of the
classical combinatorial definition of Schur functions.
Although our work here concerns only the special case when $t=1$,
the family of strong $k$-tableaux also seems to play a natural
role in the general theory of $k$-Schur functions.

Proposition~\ref{P:scover} allows us to give a conjectural monomial
expansion for $s_\lambda^{(k)}(x;t)$.
Suppose $C = \stc{\mu}{i,j}{\lambda}$ is
a marked strong cover of $k+1$-cores.  By Proposition~\ref{P:scover}
the skew shape $\lambda/\mu$ consists of $m$ identical ribbons each
of height $h$, where the height of a ribbon is the number of rows it
occupies. Define the {\it spin} of $C$ to be ${\rm spin}(C) = m(h-1)
+ (l-1)$ where the marked ribbon in $C$ is the $l$-th one from the
top.  The ${\rm spin}(T)$ of a strong tableau $T$ is the sum of
the spins of its marked strong covers.  For example, the strong
tableau of Example~\ref{X:strongtableaucore} has spin 2 = 0 + 0 + 1
+ 0 + 1. Our use of the name ``spin'' is due to the similarities
between this statistic and Lascoux, Leclerc and Thibon's spin
statistic for ribbon tableaux~\cite{LLT}.

\begin{conj} \label{conj:kSchurt}
Let $\lambda \in \Bounded$ where $n = k+1$.  The $k$-Schur functions
of~\cite{LLM,LMfil} are such that
$$
s_\lambda^{(k)}(x;t)= \tilde s_{\lambda}^{(k)}(x;t)
= \sum_Q x^{\wt(Q)} \, t^{{\rm spin}(Q)} \, ,
$$
where the summation runs over strong tableaux $Q$ of shape given by
the $n$-core $\bounded^{-1}(\lambda)$.
\end{conj}
The conjecture has been tested for all $k$ and all
$k$-bounded partitions up to degree
10. We should note that this conjecture is still open in the case $t=1$, since
the $k$-Schur functions at $t=1$ that we use in this article are those
of \cite{LMproofs}, which are only conjectured to coincide with the special
case
$t=1$ of the $k$-Schur functions of \cite{LLM,LMfil}.


\section{Enumeration of standard strong and weak tableaux}
\label{sec:enu} A strong (or weak) tableau of shape $w/v$ is {\it
standard} if its weight is the composition $(1,1,\ldots,1,0,\ldots)$
where the number of 1's is equal to $\ell(w)-\ell(v)$.  We will now
briefly discuss the enumeration of standard strong and weak tableaux
for small values of $n$.  This will also serve as examples for the
material of this section. Unfortunately, we have not been able to
find elegant formulae, similar to the hook-length formula, for these
numbers in general.

For $w,v \in \tS$, let $f^{w/v}_{\strong}$ (or $f^{w/v}_{\weak}$)
denote the (finite) number of standard strong (or weak) tableaux of
shape $w/v$.  As usual when $v = \id$, we will just write
$f^{w}_{\strong}$ (or $f^{w}_{\weak}$).  The affine insertion
bijection of Theorem~\ref{T:main} (with $l = 0$) has the following
immediate consequence when restricted to permutation matrices: for
each $m \in \Z_{>0}$ we have
\begin{equation}\label{E:fact}
m! = \sum_{\substack{w \in \tS^0 \\ \ell(w) = m}} f^w_\strong \,
f^w_\weak. \end{equation} Note that we have used the fact that
$f^w_\strong = 0$ if $w \notin \tS^0$.

\begin{remark}
The number $f^w_\weak$ for $w \in \tS$ is equal to the number of
reduced words of $w$.  They were studied to some extent
in~\cite{LamAS,LMcore}.  It follows from \cite{LamAS} that for each
$w \in \tS$ we have a formula
$$
f^w_\weak = \sum_{\substack{v \in \tS^0 \\ \ell(v) = \ell(w)}}
a_{wv} f^v_\weak
$$
for some integers $a_{wv} \in \Z$ arising from geometry.  It is
known from \cite{Lam,Pet} that $a_{wv}$ are in fact positive.
\end{remark}

\subsection{Case $n = 2$} Let $n = 2$.  For each $m \in \Z_{>0}$
there is a unique $w = w_m \in \tS^0$, which has a unique reduced
word of the form $w = \cdots s_0s_1s_0$.  The core $c(w_m)$ is the
staircase $(m,m-1,\ldots,1)$ which has bounded partition $(1^m)$.
There is a unique weak tableau with shape $c(w_m)$, of the form
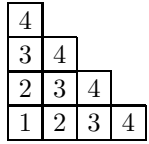
\begin{align*}
\tableau{{4}\\{3}&{4}\\{2}&{3}&{4}\\{1}&{2}&{3}&{4}}
\end{align*}
Thus \eqref{E:fact} says that $m! = f^{w_m}_\strong$.  Indeed, up to
marking, there is only one standard strong tableaux with shape
$w_m$, and each of the skew shapes $c(w_i)/c(w_{i-1})$ is a disjoint
union of $i$ squares.  Thus there are a total of $m!$ choices for
the markings.  Note that there is no condition on the markings to be
increasing, since all strong strips are of size 1 (or in other
words, all strong covers belong to different strong strips).

\subsection{Case $n=3$} Let $n = 3$ and $m \in \Z_{>0}$.  There are
$\lfloor m/2 \rfloor + 1$ Grassmannian elements with length $m$. We
denote by $w_{m,\ell} \in \tS^0$ the element with bounded partition
$(b \circ c)(w_{m,\ell}) = (2^\ell\,1^{m-2\ell})$ where $\ell \in
[0,\lfloor m/2 \rfloor]$.  Thus $c(w_{m,0})$ is ``tall'' and
$c(w_{m,\lfloor m/2\rfloor})$ is ``wide''.

For example when $m = 3$, we have $w_{m,0} = s_1s_2s_0$ and $w_{m,1}
= s_2s_1s_0$.  The weak order is described in the following result,
which has a straightforward proof.
\begin{prop}\label{P:3weakorder}
Let $n = 3$.  If $m$ is even then the weak covers $w \prec v$ with
$\ell(w) = m$ are given completely by $w_{m,\ell} \prec w_{m+1,\ell}$.  If
$m$ is odd then the weak covers $w \prec v$ with $\ell(w) = m$ are
given completely by $w_{m,\ell} \prec w_{m+1,\ell}$ and $w_{m,\ell} \prec
w_{m+1,\ell+1}$.
\end{prop}

\begin{remark} The Hasse diagram of the weak order of $\tilde{S}^0_3$ is in
fact (part of) a honeycomb lattice (see \cite[Figure 1]{LMcore})
which had earlier appeared in related work of Fomin \cite{Fom} where
it was called ``Hex''.
\end{remark}

It follows immediately from Proposition~\ref{P:3weakorder} that
$$
f^{w_{m,\ell}}_\weak = \begin{cases} f^{w_{m-1,\ell}}_\weak &\mbox{if $m$ is odd,} \\
f^{w_{m-1,\ell}}_\weak + f^{w_{m-1,\ell-1}}_\weak &\mbox{if $m$ is even
and $\ell \notin \{0,m/2\}$,} \\
f^{w_{m-1,0}}_\weak &\mbox{if $m$ is even and $\ell = 0$,} \\
f^{w_{m-1,m/2-1}}_\weak &\mbox{if $m$ is even and $\ell = m/2$.}
\end{cases}
$$
Solving these recursions, one obtains
\begin{prop}\label{P:3weak}
The numbers of standard weak tableaux for $n = 3$ are given by
$$
f^{w_{m,\ell}}_\weak = \binom{\lfloor m/2 \rfloor}{\ell}.
$$
In particular the total number of weak tableaux with Grassmannian
shape of length $m$ is given by
$$
\sum_{\substack{w \in \tilde{S}^0_3 \\ \ell(w) = m}} f^w_\weak =
2^{\lfloor m/2 \rfloor}.
$$
\end{prop}

Now we turn to the calculation of $f^{w_{m,\ell}}_\strong$.  Rather
surprisingly, we find that $f^{w_{m,\ell}}_\strong$ depends only on
$m$.  We first describe the strong covers.

\begin{prop}\label{P:3strongcover}
Let $n = 3$.  Suppose $m$ is even.  Then we have $w_{m,\ell} \lessdot
w_{m+1,\ell}$, $w_{m,\ell} \lessdot w_{m+1,\ell+1}$ and $w_{m,\ell} \lessdot
w_{m+1,\ell-1}$ whenever these permutations exist.  The corresponding
difference of cores have $m/2 +  1$, $\ell+1$, and $m/2 - \ell + 1$
components respectively.  Suppose $m$ is odd. Then $w_{m,\ell} \lessdot
w_{m+1,\ell}$ and $w_{m,\ell} \lessdot w_{m+1,\ell+1}$.  The corresponding
difference of cores have $\lfloor m/2 \rfloor -\ell + 1$ and $\ell+1$
components respectively.  Furthermore this describes all the
(marked) strong covers in $\tilde{S}_3^0$.
\end{prop}
Note that in particular if $\ell(w)$ is odd then $w \lessdot v$ if
and only if $w \prec v$ and $\ell(v) = \ell(w) + 1$.
\begin{proof}
To prove the result, we use the language of cores and
Proposition~\ref{P:crystal} which says that strong order corresponds
to inclusion of cores.  The {\it width} of a partition $\lambda$ is
its first part $\lambda_1$ and its {\it height} is equal to the
number of parts.  We first note that the 3-core $c(w_{m,\ell})$ has
height $m - \ell$ and width $m - (\lfloor m/2 \rfloor - \ell)$.  This
observation and Proposition~\ref{P:crystal} is enough to show that
the stated strong covers are the only ones possible.  To show that
the stated pairs are indeed strong covers, we use the ``$k$-skew''
shape of \cite{LMcore}, which is denoted $\lambda/\mu$ in
Proposition~\ref{P:corebounded}.  The $k$-skew shapes of $w_{m,\ell}$
(which have outside shape $c(w_{m,\ell})$) are of the form
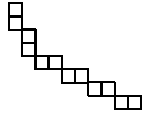
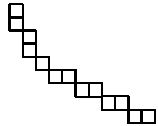
\setcellsize{5}
$$
\tableau{{}\\{}\\&{}\\&{}\\&&{}&{} \\ &&&&{}&{} \\&&&&&&{}&{} \\
&&&&&&&&{}&{}}
$$
for $m$ even, and
$$
\tableau{{}\\{}\\&{}\\&{}\\&&{} \\&&&{}&{} \\ &&&&&{}&{} \\&&&&&&&{}&{} \\
&&&&&&&&&{}&{}}
$$
for $m$ odd.  In both cases, one has $\ell$-horizontal dominoes and
$\lfloor m/2 \rfloor -\ell$ vertical dominoes.  With this description
it is easy to see the pairs stated in the proposition are indeed
strong covers with the stated number of possible markings.
\end{proof}

Using Proposition~\ref{P:3strongcover} one obtains the following
recursions for $f^w_\strong$:
$$
f^{w_{m,\ell}}_\strong = \ell\,f^{w_{m-1,\ell-1}}_\strong +(m/2 -
\ell)\,f^{w_{m-1,\ell}}_\strong$$ for $m$ even, and
$$ f^{w_{m,\ell}}_\strong = \ell\,f^{w_{m-1,\ell-1}}_\strong +\lceil m/2 \rceil\,
f^{w_{m-1,\ell}}_\strong + (\lfloor m/2 \rfloor -
\ell)\,f^{w_{m-1,\ell+1}}_\strong
$$
for $m$ odd.  By induction one easily establishes

\begin{prop}
Let $n = 3$.  The numbers $f^{w_{m,\ell}}_\strong$ of standard strong
tableaux depend only on $m$ and are given by
$$
f^{w_{m,\ell}}_\strong = \frac{m!}{2^{\lfloor m/2 \rfloor}}.
$$
\end{prop}
This agrees readily with Proposition~\ref{P:3weak} and
\eqref{E:fact}.

\begin{example}
Suppose $m = 4$.  Then there are three Grassmannian elements
$w_{4,0},w_{4,1},w_{4,2}$ which have cores $(2211),(311),(42)$
respectively.  We have the following weak and strong tableaux,
illustrating the identity \eqref{E:fact}: $24 = 1 \cdot 6 + 2 \cdot
6+1 \cdot 6$.

\setcellsize{11}
\medskip
\begin{center}
\begin{tabular}{||c|c|c||}
\hline  $w_{4,0}$ & $w_{4,1}$ &$w_{4,2}$
\\
\hline $ \tableau{\\{4}\\{3}\\{2}&{4}\\{1}&{3}\\&}$&
$\tableau{{4}\\{3}\\{1}&{2}&{3}}$ \
$\tableau{{3}\\{2}\\{1}&{3}&{4}}$  &
$\tableau{{3}&{4}\\{1}&{2}&{3}&{4}}$
\\
\hline
\end{tabular}
\end{center}

\medskip

\begin{center}
\begin{tabular}{||c|c|c||}
\hline $w_{4,0}$ & $w_{4,1}$ &$w_{4,2}$ \\
\hline $\tableau{\\{4^*}\\{3}\\{3^*}&{4}\\{1^*}&{2^*}}$ \ $
\tableau{\\{4}\\{3}\\{3^*}&{4^*}\\{1^*}&{2^*}}$ &
$\tableau{{3}\\{3^*}\\{1^*}&{2^*}&{4^*}}$ \
$\tableau{{4^*}\\{2^*}\\{1^*}&{3}&{3^*}}$  &
$\tableau{{2^*}&{4}\\{1^*}&{3}&{3^*}&{4^*}}$ \
$\tableau{{2^*}&{4^*}\\{1^*}&{3}&{3^*}&{4}}$
\\
$ \tableau{\\{4}\\{3}\\{2^*}&{4^*}\\{1^*}&{3^*}}$ \ $ \tableau{\\{4^*}\\{3}\\{2^*}&{4}\\{1^*}&{3^*}}$ & $\tableau{{4^*}\\{3}\\{1^*}&{2^*}&{3^*}}$ \ $\tableau{{4^*}\\{3^*}\\{1^*}&{2^*}&{3}}$ & $\tableau{{3^*}&{4^*}\\{1^*}&{2^*}&{3}&{4}}$ \ $\tableau{{3}&{4^*}\\{1^*}&{2^*}&{3^*}&{4}}$ \\
$ \tableau{\\{4}\\{3^*}\\{2^*}&{4^*}\\{1^*}&{3}\\&}$ \ $ \tableau{\\{4^*}\\{3^*}\\{2^*}&{4}\\{1^*}&{3}\\&}$  & $\tableau{{3^*}\\{2^*}\\{1^*}&{3}&{4^*}}$ \ $\tableau{{3}\\{2^*}\\{1^*}&{3^*}&{4^*}}$  & $\tableau{{3^*}&{4}\\{1^*}&{2^*}&{3}&{4^*}}$ \ $\tableau{{3}&{4}\\{1^*}&{2^*}&{3^*}&{4^*}}$\\
 \hline
\end{tabular}
\end{center}

\end{example}

\chapter{Affine Insertion in Terms of Cores}
\label{ch:coresinsertion}

We now translate the affine insertion algorithm for the special case
of Grassmannian elements, into the language of cores. By the proof
of Theorem \ref{T:main} it suffices to describe the local rule of
Chapter \ref{ch:forward} and its reverse in Chapter
\ref{ch:inverse}, in terms of cores.

\section{Internal insertion for cores} \label{SS:intinscore} Let
$C$ be a strong cover $\mu \lessdot t_{r,s} \cdot \mu = \la$ with
$m(C)=s$. By Proposition \ref{P:strongcovercore} this is equivalent
to marking the cell of $\la/\mu$ on the diagonal $s-1$. We shall use
this equivalence freely without further mention. We recall from
Section \ref{SS:comm} the definitions of initial pair, final pair,
and commutation. The general definition for internal insertion is in
Section \ref{SS:intins}.

We first express the condition for commutation in terms of cores.

Let $W=(\ws{\mu}{A}{\nu})$ be a weak strip and $(W,S'_1)$ the input
final pair for the internal insertion at $C$, and let $(W',S')$ be
the output final pair, with $\out(W')=\gamma$.

In every case we shall define $W'=(\ws{\la}{A'}{\gamma})$ in terms
of a set $A' \subsetneq \Z/n\Z$.

\begin{lem} \label{L:corenoncomm} $(W,C)$ does not commute if and
only if the marked component of $\la/\mu$ is contained in $\nu/\mu$.
Moreover, if this occurs then each component of $\la/\mu$ is
contained within a single row.
\end{lem}
\begin{proof} By Lemma \ref{L:Fcyclicobstruction}, $(W,C)$ does not
commute if and only if the marked component of $C$ is contained in
the skew shape $\nu/\mu$ and the tail (northwestmost cell) of the
marked component has residue equal to the minimum of some cyclic
interval $I$ in $A$. Since $\nu/\mu$ corresponds to a weak strip, by
Lemma \ref{L:weakktab} it is a horizontal strip in the usual sense.
This implies that each of the ribbon components of $\la/\mu$, lies
in a single row. In particular the above condition regarding the
tail, automatically holds.
\end{proof}

Using Lemma \ref{L:corenoncomm} we rephrase Cases A, B, and C of
internal insertion (which occur in Sections \ref{sec:bija},
\ref{sec:bijb}, and \ref{sec:bijc} respectively) in terms of cores.

\subsection{Commuting case for cores} \

\medskip\noindent\textbf{Case A (Commuting case)} Suppose
the marked component of $\la/\mu$ is not contained in $\nu/\mu$.
Then let $A'=A$ and $S'=S'_1\cup C'$ where $C'$ is the strong cover
$\nu \lessdot \gamma = t_{c_A(r),c_A(s)} \cdot \nu$, with
$m(C')=c_A(s)$.

\subsection{Noncommuting cases for cores} Let $p_0$ be the diagonal
of the head of the marked component of $C$ and
$\Ah=A\backslash\{\ba{p_0}\}$.

\smallskip \noindent \textbf{Case B (Normal bumping case)}
Suppose the marked component of $\la/\mu$ is contained in $\la/\nu$
and that either $\size(S'_1)=0$, or $\size(S'_1)>0$ and $m(C)\ne
m(\lc(S'_1))$. Let $q_0$ be the diagonal of the tail of the marked
component of $\la/\mu$. Let $q<p$ be the unique pair of consecutive
$\Ah$-nice integers such that $q<q_0$, $\la$ has an addable cell in
the diagonal $q$, and $q$ is maximal. Let $A' = \Ah\cup\{\ba{p-1}\}$
and $S'=S'_1\cup C'$ where $C'$ is the strong cover
$\nu\lessdot\gamma:= t_{c_A(q),c_A(p)}\cdot\nu$ with $m(C')=c_A(p)$.

\smallskip \noindent \textbf{Case C (Replacement Bump)}
Suppose the marked component of $\la/\mu$ is contained in $\la/\nu$,
$\size(S'_1)>0$, and $m(C)=m(\lc(S'_1))$. Let $\nu^-$ be the core
obtained by removing $\lc(S'_1)$ from $\nu$. Let $q<p$ be the unique
pair of consecutive $\Ah$-bad integers such that $q < p_0$, $\nu^-$
has an addable cell in the diagonal $q$, and $q$ is maximal. Set
$A'=\Ah\cup \{\ba{q}\}$ and $S'=(S'_1\backslash \lc(S'_1)) \cup \Cpm
\cup C'$ where $\Cpm$ is the strong cover $\nu^- \lessdot \nu' :=
t_{q,p} \cdot \nu^-$ with $m(\Cpm)=p$, and $C'$ is the strong cover
$\nu' \lessdot \gamma$, with $m(C')=t_{q,p}\,m(C)$.

\section{External Insertion for cores (Case X)}
\label{SS:coreexternalInsertion} Recall the general case of external
insertion given in Subsection \ref{sec:externalInsertion}. Let
$W=(\ws{\mu}{A}{\nu})$ be a weak strip. Let $q=\nu_1$ and $p>q$ the
next larger $A$-bad integer. Let $A' = A \cup \{\ba{q}\}$ and
$W'=(\ws{\mu}{A'}{\gamma})$. Let $C'$ be the strong cover $\nu
\lessdot \gamma := t_{q,p}\,\nu$ with $m(C')=p$. This marks the
component of $\gamma/\nu$ in the first row.

External insertion on $(W,S'_1)$ is given by computing $(W',C')$ as
above and setting $S'=S'_1\cup C'$.

\section{An example}

\begin{example} Let $n=3$ and let
\begin{align*}
m = \begin{pmatrix} 0&1&0 \\0&0&2 \\ 1&0&1
\end{pmatrix}.
\end{align*}
The growth diagram to compute the image $(P,Q)$ of $m$ is given by
Figure \ref{F:growth}. Each row of the diagram defines a strong
tableau, which is depicted to the right of the row. %
\setcellsize{13}%
\Yboxdim{5pt}%
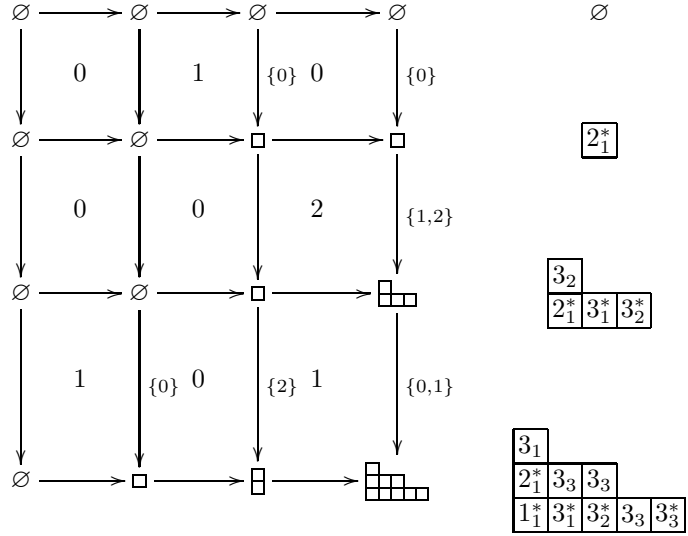
\begin{figure}
\begin{align*}
\xymatrix@=10pt{%
{\vn} \ar[dd] \ar[rr] && {\vn} \ar[dd] \ar[rr]&&
{\vn}\ar[dd]^{\{0\}} \ar[rr] && {\vn}\ar[dd]^{\{0\}} && {\vn} \\
 & 0 && 1 && 0 & \\
{\vn}\ar[dd] \ar[rr] && {\vn} \ar[dd] \ar[rr] && {\yng(1)}\ar[dd]
\ar[rr] &&
{\yng(1)}\ar[dd]^{\{1,2\}} && {\tableau{{2_1^*}}}  \\
 & 0 && 0 && 2 & \\
{\vn}\ar[dd] \ar[rr] && {\vn}\ar[dd]^{\{0\}} \ar[rr] &&
{\yng(1)}\ar[dd]^{\{2\}} \ar[rr]  &&
{\yng(1,3)}\ar[dd] ^{\{0,1\}} && {\tableau{{3_2} \\ {2_1^*}&{3_1^*}&{3_2^*}}} \\
 & 1 && 0 && 1 & \\
{\vn}\ar[rr] && {\yng(1)}\ar[rr] && {\yng(1,1)}\ar[rr] &&
{\yng(1,3,5)} && {\tableau{{3_1} \\ {2_1^*}&{3_3}&{3_3}
\\{1_1^*}&{3_1^*}&{3_2^*}&{3_3}&{3_3^*}}}}
\end{align*}
\label{F:growth} \caption{Growth diagram}
\end{figure}
The tableau $Q$ is given by%
\setcellsize{11}%
\begin{align*}
  Q = \tableau{{3}\\{2}&{3}&{3}\\{1}&{2}&{2}&{3}&{3}}.
\end{align*}
We explain the computation of the last row of $2\times 2$ squares.
The first square consists of a single external insertion where the
input data consists of empty weak and strong strips, all going from
$\vn$ to itself. This external insertion just adds the residue $0$,
giving the new shape $(1)$.

The second square has input strong strip $S$ that consists of a
single marked strong cover $C=(\vn \lessdot (1)=t_{0,1}(\vn))$ with
$m(C)=1$, and input weak strip $W=(\ws{\vn}{A}{(1)})$ with
$A=\{\ba{0}\}$. Since the marked component of $C$ consists of the
cell $(1,1)$ which is contained in $W$, we are in the noncommuting
case. Since the old output strip $S'_1$ is empty, Case B holds. We
have $p_0=q_0=0$, $q=-1$, and $p=0$, so that
$A'=\{\ba{-1}\}=\{\ba{2}\}$ with output weak strip
$\ws{(1)}{A'}{(1,1)}$ and the output strong strip consists of a
single cover $C'=((1) \lessdot (1,1))$ with $m(C')=0$.

The third square has the input strong strip $S=((1) \lessdot (2)
\lessdot (3,1))=C_1\cup C_2$ with $m(C_1)=2$ and $m(C_2)=3$. It has
input weak strip $W=(\ws{(1)}{A}{(1,1)})$ where $A=\{\ba{2}\}$. The
internal insertion at $C_1$ is in Case A: the marked component
consists of the cell $(1,2)$, which is not contained in $W$, which
consists of the cell $(2,1)$. This insertion produces the output
strong strip $(1,1)\lessdot(3,1)$ with mark $3$ and the same set
$A$, which now is associated with the weak strip
$\ws{(2)}{A}{(3,1)}$.

Next we perform the internal insertion at $C_2$. We reindex with
input weak strip $W=(\ws{(2)}{A}{(3,1)})$, and $A=\{\ba{2}\}$, old
output strong strip $S'_1=((1,1)\lessdot(3,1))$, and
$C=C_2=((2)\lessdot(3,1))$ with $m(C)=3$. The marked component of
$C$ is the single cell $(1,3)$ which is contained in $W$. This is
the noncommuting case. Now $m(\lc(S'_1))=3=m(C)$ so we are in Case
C. We have $p_0=2$, $\Ah=\{\}$, $\nu^-=(1,1)$, $q=1$, $p=2$. Thus
$A'=\{\ba{-2}\}=\{\ba{1}\}$, with output weak strip
$(\ws{(3,1)}{A'}{(3,1,1)}$ and output strong strip
$S'=((1,1)\lessdot (2,1,1)\lessdot (3,1,1))=C'_1\cup C'_2$, where
$m(C'_1)=2$ and $m(C'_2)=3$.

To finish up we apply a single external insertion to the output
final pair from the previous step, which are reindexed as
$W=(\ws{(3,1)}{A}{(3,1,1)})$ with $A=\{\ba{1}\}$ and
$S=((1,1)\lessdot (2,1,1)\lessdot(3,1,1))=C'_1\cup C'_2$ with
$m(C'_1)=2$ and $m(C'_2)=3$. We have $q=3$, $p=5$, so that
$A'=\{0,1\}$ and $W'=(\ws{(3,1)}{A'}{(5,3,1)})$ and
$C'=((3,1,1)\lessdot(5,3,1))$ with $m(C')=5$ and $S'=S\cup C'$.
\end{example}

\section{Standard case} The insertion algorithm becomes
particularly simple in the ``standard" case, by which we mean the
restriction of the insertion bijection to the subset of permutation
matrices. All weak strips that occur are empty or a single weak
cover. All strong strips that occur are empty or single marked
strong covers, and Case C never occurs. In particular, in the
standard case, the insertion process is ``context-free"; each
internal and external insertion is independent of previous
computations.

Let us adopt the notation of Section \ref{SS:intinscore} for
internal insertion. We have the weak cover $W=(\mu\wless s_a
\mu=\nu)$ (with singleton weak strip set $A=\{a\}$) and the strong
strip given by a single marked strong cover $S=(\mu\lessdot \mu
t_{rs}=\la)$, marked at the cell on diagonal $s-1$.

Let us consider Lemma \ref{L:corenoncomm} for the standard case: it
says that $(W,C)$ does not commute if and only if the marked
component is contained in $W$, which is a weak cover, whose
components are single cells. Hence the components of the strong
cover are single cells, and it follows that the strong cover and
weak cover coincide. Therefore the noncommuting case is simply
described by the conditions $\mu\ne \la=\nu$, and the marked
component is the single cell of $\la/\mu$ in the diagonal $s-1$
where $\mu \lessdot \la = \mu t_{rs}$.

So Case A occurs if and only if $\la\ne \nu$, and the computation
goes as before: the output weak cover is $\la \wless s_a \cdot
\la=\gamma$ and the output marked strong cover is $\nu \lessdot
\gamma=\nu t_{rs}$ with the cell on the same diagonal $s-1$ marked.

Case B occurs if and only if $\la=\nu$. In that case, one may show
that $A'=\{\overline{q}\}$ where $q$ is the diagonal of the first
$\la$-addable cell to the northwest of the marked cell. The marked
cell for the output marked cover is the cell in diagonal $q$.

Case C does not occur.

For external insertion, we note that it will only occur if the input
strong strip is empty. Therefore in this case, the output marked
strong cover consists of the weak cover with new residue given by
that of the $\mu$-addable cell in the first row, with that cell in
the first row marked.

\begin{example} Let $n=3$ and let
\begin{align*}
m = \begin{pmatrix} 0&0&1&0&0 \\
1&0&0&0&0 \\
0&0&0&0&1 \\
0&1&0&0&0 \\
0&0&0&1&0
\end{pmatrix}.
\end{align*}
The growth diagram is given in Figure \ref{F:growthstd}. A $*$ in a
partition diagram indicates the marking for the strong cover coming
from its left.

\newcommand{\vx}{}%
\setcellsize{6}%
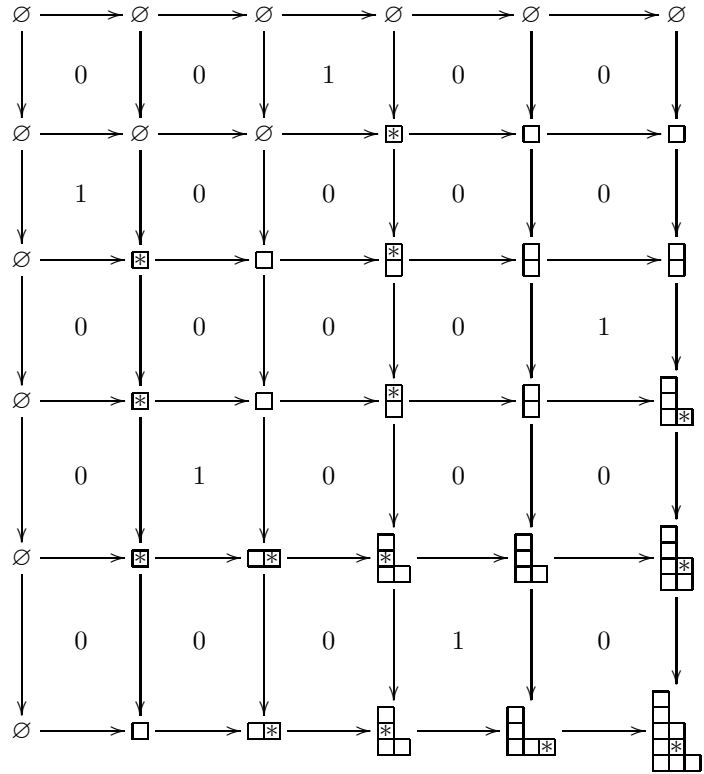
\begin{figure}
$$
\xymatrix@=10pt{%
{\vn} \ar[dd] \ar[rr] && {\vn} \ar[dd] \ar[rr] && {\vn} \ar[dd] \ar[rr] && {\vn}\ar[dd] \ar[rr] && {\vn} \ar[dd] \ar[rr] && {\vn} \ar[dd] \\
 & 0 && 0 && 1 && 0 && 0 & \\
{\vn} \ar[dd] \ar[rr] && {\vn} \ar[dd] \ar[rr] && {\vn} \ar[dd] \ar[rr] && {\tableau{*}} \ar[dd] \ar[rr] && {\tableau{{}}} \ar[dd] \ar[rr] && {\tableau{{}}} \ar[dd]\\
 & 1 && 0 && 0 && 0 && 0 & \\
{\vn} \ar[dd] \ar[rr] && {\tableau{*}} \ar[dd] \ar[rr] &&
{\tableau{{}}} \ar[dd] \ar[rr] && {\tableau{* \\ {}}} \ar[dd]
\ar[rr] &&
{\tableau{{}\\{}}} \ar[dd] \ar[rr] && {\tableau{{}\\{}}} \ar[dd] \\
 & 0 && 0 && 0 && 0 && 1 & \\
{\vn} \ar[dd] \ar[rr] && {\tableau{*}} \ar[dd] \ar[rr] &&
{\tableau{{}}} \ar[dd] \ar[rr] && {\tableau{*\\{}}} \ar[dd] \ar[rr]
&&
{\tableau{{}\\{}}} \ar[dd] \ar[rr] && {\tableau{{}\\{}\\{}&*}} \ar[dd] \\
 & 0 && 1 && 0 && 0 && 0 & \\
{\vn} \ar[dd] \ar[rr] && {\tableau{*}} \ar[dd] \ar[rr] &&
{\tableau{{}&*}} \ar[dd] \ar[rr] && {\tableau{{}\\*\\{}&{}}} \ar[dd]
\ar[rr] && {\tableau{{}\\{}\\{}&{}}} \ar[dd] \ar[rr] &&
{\tableau{{} \\ {} \\ {} & * \\ {}&{} }} \ar[dd] \\
 & 0 && 0 && 0 && 1 && 0 & \\
{\vn}  \ar[rr] && {\tableau{{}}} \ar[rr] && {\tableau{{}&*}} \ar[rr]
&& {\tableau{{}\\*\\{}&{}}} \ar[rr] && {\tableau{{}\\{}\\{}&{}&*}}
\ar[rr] && {\tableau{{}\\{}\\{}&{}\\{}&*\\{}&{}&{}}}
}%
$$
\label{F:growthstd} \caption{Growth diagram of a standard insertion}
\end{figure}
The $P$ and $Q$ tableaux are given by
\setcellsize{11}%
\begin{align*}
P = \tableau{5\\5\\3&5\\3^*&5^*\\1^*&2^*&4^*} \qquad %
Q = \tableau{5\\4\\3&5\\2&4\\1&3&5}
\end{align*}
\end{example}

\section{Coincidence with RSK as $n\to\infty$}

\begin{thm} \label{T:limitRSK} As $n\to\infty$ the bijection of
Theorem \ref{T:affGrRSK} converges to RSK row insertion.
\end{thm}
\begin{proof}
In the limit the set of $n$-bounded matrices becomes the set of all
matrices with finitely many nonzero entries. All residues of corner
cells in a partition will be distinct. Thus a weak strip becomes an
ordinary horizontal strip with one cell per residue and a weak
tableau becomes a semistandard one. A marked strong cover becomes a
skew shape consisting of a single cell, a strong strip becomes a
horizontal strip, and a strong tableau becomes a semistandard one.

We now consider the local rule. The data consists of a horizontal
strip $S=\la/\mu$, a horizontal strip $W=\nu/\mu$, and some
$e\in\Z_{\ge0}$. The algorithm consists of performing internal
insertions on $W$ at the cells of $S$ from left to right. Case A
occurs when the active cell $C$ (the one at which the internal
insertion occurs) is not in $W$, so no bumping takes place. Case B
occurs when the active cell is in $W$, in which case the active cell
$C$ is removed from the ``weak strip" and another cell $C'$ is
added, namely, the next addable one of smaller residue. This is
equivalent to bumping the cell to the end of the next row. Case C
never occurs because cells are always bumped to strictly earlier
diagonals.

Finally, each of the $e$ external insertions is given by adding a
cell to the end of the first row (and putting this cell in the new
weak and strong strips).

Since this coincides with the local rule for ordinary RSK row
insertion, the Theorem follows.
\end{proof}

\section{The bijection for $n  = 3$ and $m = 4$} In
Section~\ref{sec:enu} we described all the strong and weak tableaux
for $n =3$.  Now we give the entire (standard) bijection for $m = 4$
in this case.  For convenience, we have included the outcome of the
usual Robinson-Schensted bijection for comparison. In the following
table, the strong tableau is always drawn first.

\newpage
\setcellsize{11}%
\begin{center}
\begin{tabular}{||c|c|c||}
\hline $w$ & $n = 3$ & $n = \infty$ \\ \hline
1234 & \strongcf \ \wca &  $\tableau{1&2&3&4}$ \ $\tableau{1&2&3&4}$  \\
1243 & \strongcd \ \wca &$\tableau{4 \\ 1 & 2 &3}$ \ $\tableau{{4} \\ {1} & {2} &{3}}$ \\
1324 & \strongce \ \wca &$\tableau{3 \\ 1 & 2 &4}$\ $\tableau{{3} \\ {1} & {2} &{4}}$\\
1342 & \strongcc \ \wca & $\tableau{3 \\ 1 & 2 &4}$\ $\tableau{{4} \\ {1} & {2} &{3}}$\\
1423 & \strongbc \ \wba & $\tableau{4 \\ 1 & 2 &3}$\ $\tableau{{3} \\ {1} & {2} &{4}}$\\
1432 & \strongbd \ \wba& $\tableau{4 \\ 3\\ 1 & 2}$\ $\tableau{4 \\ 3\\ 1 & 2}$\\
2134 & \strongbf \ \wbb & $\tableau{2 \\ 1 & 3 &4}$\ $\tableau{2 \\ 1 & 3 &4}$\\
2143 & \strongac \ \waa & $\tableau{2&4\\1&3}$\ $\tableau{2&4\\1&3}$\\
2314 & \strongca \ \wca & $\tableau{2 \\ 1 & 3 &4}$\ $\tableau{3 \\ 1 & 2 &4}$\\
2341 & \strongcb \ \wca & $\tableau{2 \\ 1 & 3 &4}$ \ $\tableau{4 \\ 1 & 2 &3}$\\
2413 & \strongbf \ \wba & $\tableau{2&4\\1&3}$\ $\tableau{3&4\\1&2}$\\
2431 & \strongbb \ \wba & $\tableau{4 \\ 2\\ 1 & 3}$ \ $\tableau{4
\\ 3\\ 1 & 2}$\\
3124 & \strongba \ \wbb & $\tableau{3 \\ 1 & 2 &4}$\ $\tableau{2 \\ 1 & 3 &4}$\\
&& \\ \hline
\end{tabular}
\end{center}
\newpage
\begin{center}
\begin{tabular}{||c|c|c||}
\hline $w$ & $n = 3$ & $n = \infty$ \\ \hline 3142 & \strongab \
\waa & $\tableau{3&4\\1&2}$\
$\tableau{2&4\\1&3}$\\
3214 & \strongbe \ \wbb & $\tableau{3 \\ 2\\ 1 & 4}$ \ $\tableau{3 \\ 2\\ 1 & 4}$\\
3241 &  \strongae \ \waa & $\tableau{3 \\ 2\\ 1 & 4}$ \ $\tableau{4 \\ 2\\ 1 & 3}$\\
3412 & \strongba \ \wba & $\tableau{3&4\\1&2}$ \ $\tableau{3&4\\1&2}$\\
3421 & \strongbe \ \wba & $\tableau{3 \\ 2\\ 1 & 4}$ \ $\tableau{4 \\ 3\\ 1 & 2}$\\
4123 & \strongbc \ \wbb & $\tableau{4 \\ 1 & 2 &3}$\  $\tableau{{2} \\ {1} & {3} &{4}}$\\
4132 & \strongbd \ \wbb & $\tableau{4 \\ 3\\ 1 & 2}$\ $\tableau{4 \\ 2\\ 1 & 3}$\\
4213 & \strongad \ \waa & $\tableau{4 \\ 2\\ 1 & 3}$\ $\tableau{3 \\ 2\\ 1 & 4}$\\
4231 & \strongbb \ \wbb & $\tableau{4 \\ 2\\ 1 & 3}$\ $\tableau{4 \\ 2\\ 1 & 3}$\\
4312 & \strongaa \ \waa & $\tableau{4 \\ 3\\ 1 & 2}$ \ $\tableau{3 \\ 2\\ 1 & 4}$\\
4321 & \strongaf \ \waa & $\tableau{4\\ 3 \\ 2\\ 1 }$\  $\tableau{4 \\ 3\\ 2\\
1}$\\ && \\ \hline
\end{tabular}
\end{center}

\backmatter

\printindex

\begin{thebibliography}{xx}

\bibitem{BS} {\sc N.~Bergeron and F.~Sottile:}
Schubert polynomials, the Bruhat order, and the geometry of flag
manifolds,  {\sl Duke Math. J.} \textbf{95}, (1998), 373--423.

\bibitem{BB} {\sc A.~Bj\"{o}rner and F.~Brenti:} Affine permutations of type $A$,
{\sl Electron. J. Combin.} \textbf{3}/2 (1996), Research Paper 18.

\bibitem{Bott} {\sc R. Bott:}
The space of loops on a Lie group, {\sl Michigan Math. J.}
\textbf{5} (1958) 35--61.

\bibitem{Fom} {\sc S. Fomin:}
Schensted algorithms for dual graded graphs, {\sl J. Algebraic
Combin.} \textbf{4} (1995), no. 1, 5--45.

\bibitem{Ful} {\sc W.~Fulton:}
Young tableaux. With applications to representation theory and
geometry, {\sl London Mathematical Society Student Texts}
\textbf{35} Cambridge University Press, Cambridge, 1997.

\bibitem{GR} {\sc H.~Garland and M.~S.~Raghunathan:} A Bruhat decomposition for the
loop space of a compact group: a new approach to results of Bott,
{\sl Proc. Nat. Acad. Sci. U.S.A.} \textbf{72} (1975), no. 12,
4716--4717.

\bibitem{[GH]}
{\sc A.~M.~Garsia and M.~Haiman:} A graded representation module
for Macdonald's polynomials,
{\sl Proc. Natl. Acad. Sci. USA} {\bf 90} (1993), 3607-3610.


\bibitem{Gra} {\sc W.~Graham:} Positivity in equivariant Schubert
calculus, {\sl Duke Math. J.} \textbf{109}, no. 3 (2001), 599--614.

\bibitem{[Ha]}
{\sc M.~Haiman:} Hilbert schemes, polygraphs, and the Macdonald
positivity conjecture, {\sl J. Amer. Math Soc.} {\bf{14}} (2001), 941--1006.

\bibitem{H} {\sc J.~Humphreys:} Reflection groups and Coxeter groups,
{\sl Cambridge Studies in Advanced Mathematics}, \textbf{29},
Cambridge University Press, Cambridge, 1990.

\bibitem{[Ji]} {\sc N.~Jing:} Vertex operators and Hall-Littlewood symmetric
functions, {\sl Adv. Math.} {\bf 87} (1991), 226--248.


\bibitem{Kac} {\sc V.~G.~Kac:} Infinite-dimensional Lie algebras,
Cambridge University Press, Cambridge, 1994.

\bibitem{KK} {\sc B.~Kostant and S.~Kumar:} The nil Hecke ring and the
cohomology of $G/P$ for a Kac-Moody group $G$, {\sl Adv. Math.}
\textbf{62} (1986), 187-237.

\bibitem{Kum} {\sc S.~Kumar:}  Kac-Moody groups, their flag
varieties and representation theory, {\sl Progress in Mathematics}
\textbf{204} Birkh\"{a}user Boston, Inc., Boston, MA, 2002.

\bibitem{LamAS} {\sc T.~Lam:} Affine Stanley symmetric functions, {\sl
    Amer. J. Math.} \textbf{128} (2006), no. 6, 1553--1586.

\bibitem{Lam} {\sc T.~Lam:} Schubert polynomials for the affine
Grassmannian, {\sl J. Amer. Math. Soc.} \textbf{21} (2008), 259--281.

\bibitem{LS} {\sc T.~Lam and M.~Shimozono:} Dual graded graphs for
Kac-Moody algebras, {\sl Algebra and Number Theory}, \textbf{1}
(2007), 451--488.

\bibitem{LLM} {\sc L.~Lapointe, A.~Lascoux, and J.~Morse:}
Tableau atoms and a new Macdonald positivity conjecture, {\sl Duke
Math. J.} \textbf{116} (2003),  no. 1, 103--146.


\bibitem{LMfil} {\sc L.~Lapointe and J.~Morse:}
Schur function analogs for a filtration of the symmetric function
space, {\sl J. Combin. Theory Ser. A}  \textbf{101} (2003),  no. 2,
191--224.


\bibitem{LMcore} {\sc L.~Lapointe and J.~Morse:}
Tableaux on $k+1$-cores, reduced words for affine permutations, and
$k$-Schur expansions, {\sl J. Combin. Theory Ser. A} \textbf{112} (2005),
no. 1, 44--81.

\bibitem{LMproofs} {\sc L. Lapointe and J. Morse:}
A $k$-tableaux characterization of $k$-Schur functions,
{\sl Adv. Math.} \textbf{213} (2007), no. 1, 183--204.


\bibitem{LMdual} {\sc L. Lapointe and J. Morse:}
Quantum cohomology and the $k$-Schur basis,
{\sl Trans. Amer. Math. Soc.}, posted on October 5, 2007,
PII S 0002-9947(07)04287-0 (to appear in print).

\bibitem{LMW} {\sc L. Lapointe, J. Morse, and M. Wachs:}
Type A affine Weyl group and the $k$-Schur
functions, unpublished.

\bibitem{L} {\sc A. Lascoux:} Ordering the affine symmetric group, {\sl Algebraic
combinatorics and applications} (G\"{o}{\ss}weinstein, 1999),
219--231, Springer, Berlin, 2001.

\bibitem{LLT} {\sc A.~Lascoux, B.~Leclerc, and J.-Y.~Thibon:}
Ribbon tableaux, Hall-Littlewood symmetric functions, quantum affine
algebras, and unipotent varieties, {\sl J. Math. Phys.}
\textbf{38}(3) (1997), 1041--1068.

\bibitem{LaSc} {\sc A.~Lascoux and M.~P.~Sch\"{u}tzenberger:}
Croissance de polyn\^{o}mes de Foulkes-Green, {\sl C. R. Acad. Sci.
Paris} \textbf{288} (1979), 95--98.



\bibitem{vL} {\sc M.~van Leeuwen:}
Edge sequences, ribbon tableaux, and an action of affine
permutations, {\sl Europ. J. Combinatorics} \textbf{20} (1999),
179--195.

\bibitem{Lus} {\sc G. Lusztig:}  Some examples of square integrable representations of
semisimple $p$-adic groups, {\sl Trans. Amer. Math. Soc.}
\textbf{277} (1983), no. 2, 623--653.

\bibitem{Mac} {\sc I.~G.~Macdonald:}
Symmmetric Functions and Hall Polynomials, Oxford University
Press, 1995.

\bibitem{MM} {\sc K.C. Misra and T. Miwa:} Crystal base for the basic
representation of $U_q(\widehat{\mathfrak{sl}}(n))$, {\sl Comm.
Math. Phys.} \textbf{134} (1990), no. 1, 79--88.

\bibitem{Pak} {\sc I.~Pak:} Periodic permutations and the Robinson-Schensted correspondence,
preprint, 2003.

\bibitem{Pet} {\sc D.~Peterson:} Lecture notes at MIT, 1997.

\bibitem{PS} {\sc A.~Pressley and G.~Segal:} Loop groups, Clarendon
Press, Oxford, 1986.

\bibitem{Shi} {\sc J. Shi:} The Kazhdan-Lusztig cells in certain affine Weyl
  groups, {\sl Lecture
Notes in Mathematics} \textbf{1179}, Springer-Verlag, Berlin, 1986.

\bibitem{SW} {\sc M.~Shimozono and J.~Weyman:} Graded characters of modules
  supported in the closure of a nilpotent conjugacy class, {\sl European J.
Combin.} {\bf 21} (2000), 257--288.

\bibitem{SZ} {\sc M. Shimozono and M. Zabrocki:} Hall-Littlewood vertex
  operators and generalized Kostka polynomials, {\sl Adv. Math.}
{\bf 158} (2001), 66--85.

\bibitem{Sot} {\sc F.~Sottile:}
Pieri's formula for flag manifolds and Schubert polynomials, {\sl
Annales de l'Institut Fourier} \textbf{46} (1996), 89--110.

\end{thebibliography}
\end{document}